\documentclass{article}
\usepackage{float}
\usepackage{subfiles}
\usepackage{amsmath,amsfonts,amscd,amssymb,enumitem,color}
\usepackage{float}
\usepackage{hyperref}
\usepackage{layout}
\usepackage{graphicx}
\usepackage{cite}
\usepackage{color}
\usepackage{appendix}
\usepackage{amssymb,latexsym,amsthm}
\usepackage{amsfonts}
\usepackage{mathtools}

\usepackage{xcolor}

\usepackage{nomencl}
\makenomenclature

\usepackage{tikz-cd}
\usetikzlibrary{matrix, calc, arrows}
\usetikzlibrary{arrows.meta}

\def\newmathop#1{\expandafter\gdef\csname #1\endcsname{\mathop{\rm #1}\nolimits}}
\def\newvmathop#1{\expandafter\gdef\csname v#1\endcsname{\mathop{\rm #1}\nolimits}}

\theoremstyle{plain}
\newcounter{thmcount}[section]

\newtheorem{theorem}[thmcount]{Theorem}
\newtheorem{corollary}[thmcount]{Corollary}
\newtheorem{lemma}[thmcount]{Lemma}
\newtheorem{proposition}[thmcount]{Proposition}
\newtheorem{conjecture}[thmcount]{Conjecture}

\newtheorem*{conjecture*}{Conjecture}
\newtheorem*{theorem*}{Theorem}
\newtheorem*{corollary*}{Corollary}
\newtheorem*{lemma*}{Lemma}

\newtheorem*{rough*}{Rough idea}

\newtheorem*{motivation*}{Motivation}

\newtheorem*{goal*}{Goal}

\theoremstyle{definition}
\newtheorem{remark}[thmcount]{Remark}
\newtheorem{example}[thmcount]{Example}
\newtheorem*{example*}{Example}
\newtheorem*{remark*}{Remark}
\newtheorem{definition}[thmcount]{Definition}
\newtheorem{notation}[thmcount]{Notation}

\numberwithin{equation}{section}

\def\Q{{\mathbb Q}}
\def\F{{\mathbb F}}

\def\Fpb{{\overline{\mathbb F}_p}}

\def\Qp{{{\mathbb Q}_p}}

\def\Qpb{{\overline{\mathbb Q}_p}}
\def\Z{{\mathbb Z}}
\def\Zp{{\mathbb Z}_p}
\def\Zpb{{\overline{\mathbb Z}_p}}

\def\cP{{\mathcal P}}

\def\fS{{\mathfrak S}}

\newmathop{Disc}
\newmathop{Tr}
\newmathop{Norm}
\newmathop{ord}
\newmathop{GL}
\newmathop{SL}
\newmathop{PGL}
\newmathop{Hom}
\newmathop{Ind}
\newmathop{Res}
\newmathop{rk}
\newmathop{corank}
\newmathop{coker}
\newmathop{codim}
\newmathop{cyc}
\newmathop{Reg}
\newmathop{Gal}
\newmathop{Sel}

\newmathop{BSD}
\newmathop{tors}
\newmathop{tr}
\newmathop{Ext}
\newmathop{Mod}

\newmathop{id}
\newmathop{Aut}
\newmathop{Art}

\newmathop{Rep}
\DeclareMathOperator{\HT}{HT}

\newcommand{\ol}[1]{\overline{#1}}

\newcommand{\vr}{\vec{r}}

\def\fM{{\mathfrak M}}
\def\ofM{{\overline{\mathfrak{M}}}}
\def\hofM{{\hat{\overline{\mathfrak{M}}}}}
\def\hM{{\hat{\mathfrak{M}}}}
\def\fN{{\mathfrak N}}
\def\ofN{{\overline{\mathfrak{N}}}}

\def\ofP{{\overline{\mathfrak{P}}}}

\def\fS{{\mathfrak S}}
\def\catKis{{\Mod^{\phi,r}_{\fS}}}
\def\catKistor{{\Mod^{\phi,r}_{\fS, \text{tor}}}}

\def\Tkis{{T_{\fS}}}

\newcommand{\mtwo}[4]{\begin{pmatrix} #1 & #2  \\ #3 & #4 \end{pmatrix}}
\newcommand{\mtwosmall}[4]{\begin{psmallmatrix} #1 & #2  \\ #3 & #4 \end{psmallmatrix}}

\newcommand{\liftrep}{\tilde{\rho}}

\newcommand{\modchar}{\overline{\psi}}
\newcommand{\mcpair}{\overline{\psi}_1,\overline{\psi}_2}

\newcommand{\Frob}{\text{Fr}}

\newcommand{\Block}{\text{Block}}
\newcommand{\wtop}{s}
\newcommand{\wbottom}{t}
\newcommand{\lowparam}{n_1}
\newcommand{\highparam}{n_2}
\newcommand{\singleparam}{n}

\title{Crystalline liftability of irregular weights}
\author{Hanneke Wiersema }
\date{June 2025}
\begin{document}

\maketitle

\begin{abstract}

Let $p$ be an odd prime. Let $K/\Q_p$ be a finite unramified extension. Let $\rho: G_K \to \GL_2(\Fpb)$ be a continuous representation. We prove that $\rho$ has a crystalline lift of small irregular weight if and only if it has multiple crystalline lifts of certain specified regular weights. The inspiration for this result comes from work of Diamond-Sasaki on geometric Serre weight conjectures. Our result provides a way to translate results currently formulated only for regular weights to also include irregular weights. The proof uses results on Kisin and $(\varphi,\hat{G})$-modules obtained from extending recent work of Gee-Liu-Savitt to study crystalline liftability of irregular weights.
\end{abstract}

\tableofcontents
\section{Introduction}
Serre’s modularity conjecture \cite{serre1973}, states that for any prime $p$, any irreducible, continuous, odd representation $\rho: G_{\Q} \to \GL_2(\Fpb)$ is isomorphic to the Galois representation associated to a modular form.  The strong version \cite{serreduke} of the statement includes a prediction for the minimal weight $k \geq 2$ such that a modular form of weight $k$ and level prime to $p$ corresponds to a given Galois representation, commonly referred to as the {weight part} of Serre's conjecture. The full conjecture was proven by Khare and Wintenberger \cite{khareI}, \cite{khare} and Kisin \cite{kisinserre}.

\par
More generally, it is conjectured that such a connection should also exist between mod $p$ representations of the absolute Galois group of a totally real field $F$ and Hilbert modular forms. This is still open. However, if, for a given Galois representation $\rho: G_F \to \GL_2(\Fpb)$ one assumes such a connection exists, i.e. that $\rho$ is modular of some Hilbert modular form, one can still try to predict the weights of the forms such that $\rho$ is modular with respect to these forms. This is called the weight part of Serre's conjecture. The first breakthrough for the weight part in the more general setting  came through the work of Buzzard, Diamond and Jarvis \cite{bdj}.

Buzzard, Diamond and Jarvis reformulated modularity in terms of $\rho$ arising in the Jacobian of Shimura curves associated to quaternion algebras over $F$. In \cite{bdj}, assuming $\rho$ is modular of some weight, the authors predict the set of weights such that a representation $\rho$ will be modular of these weights, now known as the Buzzard--Diamond--Jarvis conjecture. 

In their paper \cite{GLS14}, Gee, Liu and Savitt prove the Buzzard--Diamond--Jarvis conjecture for unitary groups. For one direction of the conjecture their methods are purely local as it consists of studying the crystalline liftability of  representations of the Galois group $G_K$, with $K/\Qp$ finite and unramified. In order to do this, they make extensive use of the theory of Kisin and $(\varphi,\hat{G})$-modules. In particular, the authors use this theory to determine the possible reductions of certain two-dimensional crystalline representations. 

In this paper, we do the same, but we include crystalline representations with non-distinct Hodge--Tate weights. This allows us to study crystalline liftability of Galois representations with respect to non-distinct Hodge--Tate weights, which we shall refer to as irregular weights. In particular, we prove that one can connect crystalline liftability of irregular weights to crystalline liftability of multiple specified regular weights. We shall now state simplified versions of our main results below.

\begin{theorem}[Theorem \ref{twoproof}, Proposition \ref{irreducibletoregular}, Proposition \ref{irreducibletoirregular}]
\label{thmcristwointro}
Let $K$ be an unramified finite extension of $\Qp$. A representation $\rho: G_K \to \GL_2(\Fpb)$ has a crystalline lift of a small irregular weight $({k},{0})$ with $k \neq 1$ if and only if $\rho$ has crystalline lifts of (regular) weights $(k',l')$ and $(k^{\Theta},l^{\Theta})$ where the weights and conditions on $k$ are described in Section \ref{statementsect}.
\end{theorem}

The next statement involves a subset $\tilde{M}$ of the set of embeddings of $K \to \Qpb$, which we will define in the body of the paper (see Definition \ref{deftildeM}).

\begin{theorem}[{Theorem \ref{reduciblethm}, Proposition \ref{irreducibletoregular}, Proposition \ref{irreducibletoirregular}}]
\label{thmcrisintro}
Let $K$ be an unramified finite extension of $\Qp$. A representation $\rho: G_K \to \GL_2(\Fpb)$ has a crystalline lift of a small irregular weight $({k},{0})$ with $k \neq 1$ if and only if
\begin{enumerate}
\item $\rho$ has a crystalline lift of (regular) weight $(k',0)$, and,
\item for each $\mu \in \tilde{M}$, $\rho$ has a crystalline lift of (regular) weight $({k}^{\mu},{l}^{\mu})$,
\end{enumerate}
where the weights and conditions on $k$ are described in Section \ref{statementsect}.
\end{theorem}

Originally this paper just contained Theorem \ref{thmcrisintro}. The inspiration for this theorem was an analogous result regarding modularity with respect to geometric mod $p$ Hilbert modular forms, which we will expand on below. We would like to thank Kalyani Kansal, Brandon Levin and David Savitt for conversations which led to the addition of Theorem \ref{thmcristwointro}. They have proven this result in upcoming work \cite{kls}, using a different, more geometric approach following recent work of Kansal and Savoie \cite{kansal}.

\par
Both versions of the theorem provide a way to translate results formulated only for regular weights to also include irregular weights. For example, consider the geometric Breuil--Mézard conjecture \cite{ge}, due to Emerton and Gee. Their geometric refinement of the conjecture includes an explicit description of the irreducible components of the fibre of a certain Galois deformation ring through the identification of the underlying cycles. These cycles correspond to regular weights. The theorems above suggest a formulation of the conjecture for irregular weights in terms of intersection of these cycles. For further discussion in this direction see also \cite[Remark 5.4]{BBHKLLSW}.

Let us now explain more carefully the reason why we have two versions of this theorem. The first one is obviously simpler but as we shall see the second one has implications for partial weight one \emph{geometric modularity}. Indeed, the original motivation for Theorem \ref{thmcrisintro} in this paper is the (conjectural) relation with crystalline liftability and geometric modularity as defined by Diamond and Sasaki, which we shall briefly discuss.

In \cite{DS}, Diamond and Sasaki define geometric mod $p$ Hilbert modular forms  as sections of automorphic line bundles on Hilbert modular varieties. For a mod $p$ Hilbert modular eigenform of arbitrary weight, they associate a two-dimensional representation of the absolute Galois group of a totally real field $F$, extending existing work on attaching Galois representations to Hilbert modular forms. A (global) Galois representation is then said to be \emph{geometrically modular} of weight $(k,l)$ if it arises from such a mod $p$ Hilbert modular eigenform of weight $(k,l)$. 

We now come to the relation of this to our main results. Diamond and Sasaki formulate a conjecture describing the set of weights of the mod $p$ Hilbert modular forms from which Galois representation arise in terms of crystalline liftability: 

\begin{conjecture}[Conjecture 7.3.2 from \cite{DS}]
\label{conj732}
Let $F$ be a totally real field in which $p$ is unramified. Suppose $\rho: G_F \to \GL_2(\Fpb)$ is geometrically modular of some weight and suppose $k \in  \Xi_{\min}^{+}$. Then $\rho$ is geometrically modular of weight $(k,l)$ if and only if $\rho|_{G_{F_v}}$ has a crystalline lift of weight  $(k_{\tau},l_{\tau})_{\tau \in \Sigma_v}$ for all $v|p$.
\end{conjecture}

Here $\Xi_{\min}^{+}$ is the minimal weight cone defined in \cite[Section 7.3]{DS} and $\Sigma_v$ the set of embeddings $\{\tau: F_v \to \Qpb\}$. 
The work of Diamond and Sasaki allows for the study of geometric modularity with respect to Hilbert modular forms with partial weight one, and, using the above relation therefore allows for the study of crystalline lifts with irregular weights. Naturally, the analogue of Theorem \ref{thmcrisintro} in the geometric modularity setting is a global statement:

\begin{theorem}[Theorem 3.1 from \cite{hw2022}]
\label{geommod}
Let $F$ be a totally real field in which $p$ is unramified. Let $\rho: G_F \to \GL_2(\Fpb)$ be irreducible. Let $(k,l)$ be an irregular weight. Assuming some technical conditions on the weight $(k,l)$ and on the local representations $\rho|_{G_{F_v}}$ for each $v \mid p$ (as in \cite[Section 3]{hw2022}), $\rho$ is geometrically modular of weight $(k,l)$ if and only if 
\begin{enumerate}
\item $\rho$ is geometrically modular of (regular) weight $(k',l')$, and,
\item for each $\mu \in \tilde{M}$, $\rho$ is geometrically modular of (regular) weight $({k}^{\mu},{l}^{\mu})$.
\end{enumerate}
\end{theorem}

This is an extension of Lemma 11.3.1 in \cite{DS}, which dealt with the quadratic case. Let us say a few words about the proof of the above theorem. In the first direction one starts with a mod $p$ Hilbert modular form with irregular or \emph{non-algebraic} (in the terminology of \cite{DS}) weight, meaning we have a partial weight one Hilbert modular form. We then obtain multiple Hilbert modular forms with regular weights through multiplication by partial Hasse invariants and the action of partial Theta operators. The reverse direction is harder, and relies on, amongst other things,  results by Diamond and Sasaki on the ordinariness of Galois representations as well as results on divisibility by Hasse invariants (\cite[Theorem 8.2.2]{DS}).

The analogous result of Theorem \ref{thmcristwointro} is expected to be true in the geometric modularity setting too, but would require either a different method or more general versions of results on divisibility by Hasse invariants. To facilitate the relation with geometric modularity we thus prove both versions. In particular, Theorem \ref{thmcrisintro} allows one to relate partial weight one modularity to the existence of irregular crystalline lifts and vice versa, which is work in progress.
We should note that the two versions are identical in the case where $F$ is quadratic, which is stated and proven by Diamond and Sasaki (\cite[Lemma 11.2.6]{DS}).

Our proof methods are different to those in the aforementioned quadratic result by Diamond and Sasaki. Diamond and Sasaki use Fontaine--Lafaille theory and appeal to results from Chang and Diamond \cite{Chang_Diamond} who use the theory of Wach modules to characterise $(\phi, \Gamma)$-modules corresponding to crystalline representations. We use methods of Gee, Liu and Savitt \cite{GLS14} and \cite{GLS15} instead (as per \cite[Remark 11.2.2]{DS}). 

We extend some of the main theorems in \cite{GLS14} to allow for the study of crystalline lifts with irregular weights. Throughout we take advantage of the fact that many results in \cite{GLS14} are in the generality we need, or exist in a more general form in \cite{GLS15}, and if not we extend the methods in these two papers and adapt definitions and results where necessary. 

We prove Theorem \ref{thmcristwointro} and Theorem \ref{thmcrisintro} separately for reducible and irreducible Galois representations. We explain the steps required. In both cases we shall see it is easier to move from the irregular weight to the set of regular weights than in the other direction. We will spend a lot more time on reducible representations, as this is significantly harder than dealing with irreducible representations, in particular if the representation is not semisimple. 

Let $K/\Qp$ be a finite and unramified extension. Suppose we have a reducible $\rho: G_K \to \GL_2(\Fpb)$ representation such that 
\[
\rho \cong \mtwo{\modchar_1}{\ast}{0}{\modchar_2},
\]
for some characters $\mcpair: G_{K} \to \overline{\F}^{\times}_p$. 
Let $c_{\rho}$ be the corresponding class in $H^1(G_K, \Fpb(\modchar))$ where $\modchar=\modchar_1 \modchar_2^{-1}$. 

For semisimple representations, it suffices to rewrite the characters $\modchar_1$, and $\modchar_2$ and then lift these characters. This is fairly straightforward when going from irregular weights to the regular weights. However, to obtain the irregular weight in this way, we first need to do some `matching': showing that the representation having crystalline lifts of a set of regular weights forces it to have a certain shape (potentially) corresponding to an irregular weight. We note that this is also roughly the strategy we adopt for irreducible representations.

If $\rho$ is not semisimple then we need to show that $c_{\rho}$ is contained in a distinguished subspace of the cohomology space $H^1(G_K, \Fpb(\modchar))$.
We use our extended results from \cite{GLS14} as well as results from \cite{GLS15} to create maps between extensions of rank one Kisin modules. We need to do this multiple times. This is to account for the different versions of our main results as well as the different cases in the proof for each of the theorems. We can then use this to prove equality of (subsets of) the relevant distinguished subspaces. Using the aforementioned `matching', this then completes the reducible case.

\subsubsection*{Outline of paper} 
We conclude the introduction with some notation and conventions.
In the second section we introduce some of the theory of Kisin modules and $(\varphi,\hat{G})$-modules. The third section deals with extending results from \cite{GLS14} to the situation of non-distinct Hodge--Tate weights and we use these results to study certain spaces of extension classes. In the fourth section we introduce the weights involved so we can state our main results. The fifth section contains the `matching' mentioned above, allowing us to prove the reducible semisimple case. The sixth section is the most involved, and contains the main technical results involving extensions of Kisin modules and distinguished subspaces. The final and the seventh section contains the irreducible case. 

\subsection*{Notation}
Let $p$ be an odd prime. Let $K$ be an unramified extension of $\Q_p$ of degree $f$ with residue field $\F$. 
Fix an algebraic closure $\Qpb$ of $\Qp$. Write $\Sigma$ for the embeddings of $K$ into $\Qpb$. We let $L/\Q_p$ be a sufficiently large finite extension of $\Q_p$ containing the image of every embedding in $\Sigma$.  

Let $\Fpb$ be the algebraic closure of $\F$ obtained as the residue field of the ring of integers of $\Qpb$. We identify $\Sigma$ with the set of embeddings of $\F \to \Fpb$.

Write $\Frob$ for absolute Frobenius. For any $\tau_0 \in \Sigma$ we let $\tau_i=\tau_{i+1}^p$ so that $\Sigma=\{ \tau_i \, \vert \, i \in \{ 0, \dots, f-1\} \}$. Note this means $\tau_i=\Frob \circ \tau_{i+1}$.

For each $\tau \in \Hom(\F,\Fpb)$, we write $\omega_{\tau}$ for the fundamental character corresponding to $\tau$. When there is no risk of confusion we will write $\omega_i$ for $\omega_{\tau_i}$. Our conventions then mean that $\omega_{i+1}^p=\omega_i$. 

 Write $K'$ for the quadratic unramified extension of $K$, $\F'$ for the residue field and $\Sigma'=\{ \sigma: \F' \hookrightarrow \Fpb \}$.

Let $\pi$ be an uniformiser of $K$ and write $E(u)$ for the minimal polynomial of $\pi$ over $K$, i.e. $E(u)=u-\pi$. Let $W(\F)$ be the ring of Witt vectors. Define $\mathfrak{S}=W(\F)[[u]]$. This ring is equipped with a Frobenius endomorphism $\varphi$ via $u \mapsto u^p$ along with the natural Frobenius on $W(\F)$.

We let $R=\varprojlim \mathcal{O}_{\overline{K}}/p$ where the transition maps are the $p$th power maps. This is a valuation ring with valuation defined by $v_{R}((x_n)_{n \geq 0})=\lim_{n \to \infty} p^n v_p(x_n)$, where $v_p(p)=1$.

We fix a compatible system of $p^n$th roots of $\pi$: set $\pi_0:=\pi$ and for all $n>0$ we fix some $\pi_n$ such that $\pi_n^p=\pi_{n-1}$. We also fix a compatible system of primitive $p^n$th roots of unity $\zeta_{p^n}$ and we define the following fields
\[
K_{\infty}=\cup_{n=0}^{\infty} K(\pi_n), \quad K_{p^{\infty}}=\cup_{n=1}^{\infty} K(\zeta_{p^n}), \quad \hat{K}=\cup_{n=1}^{\infty} K_{\infty}(\zeta_{p^n}).
\]
Here $\hat{K}$ is the Galois closure of $K_{\infty}$ over $K$. We further write
\[
G_{\infty}=\Gal(\overline{K}/K_{\infty}), \quad \hat{G}=\Gal(\hat{K}/K).
\]
\subsubsection*{$p$-adic Hodge theory}

If $W$ is a de Rham representation of $G_K$ over $\Qpb$ and $\kappa$ is an embedding $K \hookrightarrow \Qpb$ then the multiset $\HT_{\kappa}(W)$ of Hodge--Tate weights of $W$ with respect to $\kappa$ is defined to contain the integer $i$ with multiplicity 
\[
\dim_{\Qpb}(W \otimes_{\kappa,K} \hat{\overline{K}}(-i))^{G_K},
\]
with $\hat{\overline{K}}$ the completion of $\overline{K}$.

We wish to emphasise that with this convention the $p$-adic cyclotomic character $\chi_{\text{cyc}}$ has Hodge--Tate weight one. 

\subsubsection*{Acknowledgements.} I would like to thank Fred Diamond for his continued support, suggestions and comments. I would also like to thank Rong Zhou for comments on earlier versions. Finally I would like to thank Kalyani Kansal, Brandon Levin and David Savitt for helpful conversations. This work was supported by the Herchel Smith Postdoctoral Fellowship Fund, and the Engineering and Physical Sciences Research Council (EPSRC) grant EP/W001683/1.

\section{Kisin modules and Galois representations}
In this section we introduce the techniques from \cite{GLS14} and extend some of the main results in the paper to include irregular weights using results from \cite{GLS15}.

\subsection{Background on $\varphi$- and Kisin modules}

Before introducing Kisin modules, we need to define $\varphi$-modules. The definitions are obtained from \cite{GLS14}.

\begin{definition}
A $\varphi$-module over $\fS$ is an $\fS$-module $\fM$ equipped with a $\varphi$-semilinear map $\varphi_{\fM}: \fM \to \fM$ (we will mostly omit the subscript). A morphism between two $\varphi$-modules $(\fM_1,\varphi_1)$ and $(\fM_2,\varphi_2)$ is a $\fS$-linear morphism compatible with the maps $\varphi_i$. We have a $\fS$-linear map $1 \otimes \varphi: \fS \otimes_{\varphi,\fS} \fM \to \fM$, and we say $(\fM,\varphi)$ has height $r$ if the cokernel of $1 \otimes \varphi$ is killed by $E(u)^r$.
\end{definition}

Now we can introduce Kisin modules.

\begin{definition}
A finite free Kisin module (of height $r$) is a $\varphi$-module (of height $r$) $\fM$ whose underlying $\fS$-module is finite free. A torsion Kisin module is a $\varphi$-module $\fM$ of height $r$ that is  finitely generated over $\fS$ and which is killed by $p^n$ for some $n \geq 0$, and such that the natural map $\fM \to \fM[\tfrac{1}{u}]$ is injective.
\end{definition}

We write Kisin module to mean either a finite free Kisin module or a torsion Kisin module of some height $r$.
We write $\catKis$ and $\catKistor$ for the category of finite free Kisin modules and torsion Kisin modules of height $r$. 

If $\fM$ and $\fN$ are two Kisin modules, then we write $\fM \otimes \fN$ for the Kisin module with underlying $\fS$-module $\fM \otimes_{\fS} \fN$ and with $\varphi$-semilinear map $\varphi_{\fM} \otimes_{\fS} \varphi_{\fN}$. If $\fM$ has height $r$ and $\fN$ has height $s$, then $\fM \otimes \fN$ has height $r+s$.

Often we study Kisin modules with an additional action.  If $A$ is a finite commutative $\Z_p$-algebra we say a Kisin module $\fM$ has a natural $A$-action if $\fM$ is an $A$-module and the $A$-action commutes with the $\fS$ and $\varphi$-action on $\fM$. {In particular, if $\fM$ and $\fN$ are two Kisin modules with natural $\F_L$-action, then $\fM \otimes \fN$ is again a Kisin module with natural $\F_L$-action.}

We introduce one specific kind of rank one Kisin module with $\F_L$-action, as these are the Kisin modules we will mostly be working with.

\begin{definition}
\label{def61}
Suppose $r_0, \dots, r_{f-1}$ are non-negative integers and $a \in \F_L^{\times}$. Let $\ofM(r_0,\dots, r_{f-1};a)$ be the Kisin module with natural $\F_L$-action that is rank one over $\fS \otimes_{\Zp} \F_L$ and satisfies
\begin{itemize}
\item $\ofM(r_0,\dots,r_{f-1};a)_i$ is generated by $e_i$, and
\item $\varphi(e_{i-1})=(a)_i u^{r_i} e_i$,
\end{itemize}
where $(a)_i=a$ if $i \equiv 0 \mod f$ and $(a)_i=1$ otherwise. 
\end{definition}

Note that in this definition the $\ofM_i$ are as in \cite[Section 4.4]{GLS14}. By Lemma 6.2 of \cite{GLS14} any rank $\varphi$-module over $\fS \otimes_{\Zp} \F_L$ is isomorphic to one of the form $\ofM(r_0,\dots, r_{f-1};a)$, and we will study these frequently in later sections. 

Finally we define $(\varphi,\hat{G})$-modules. We recall enough of the definition for our purposes; see \cite[Section 5]{GLS14} for more details.

\begin{definition}
A $(\varphi,\hat{G})$-module of height $r$ is a triple $(\fM,\varphi_{\fM},\hat{\mathcal{G}})$ in which $(\fM,\varphi_{\fM})$ is a Kisin module of height $r$, and where $\hat{\mathcal{G}}$ is a $\hat{G}$-action on $\hat{\fM}:=\hat{\mathcal{R}} \otimes_{\varphi,\fS} \fM$, where $\hat{\mathcal{R}} = W(R) \cap \mathcal{R}_K$. The ring $\mathcal{R}_K$ is a certain subring inside $B_{\text{cris}}^+$, and both $\mathcal{R}_K$ and $\hat{\mathcal{R}}$ are stable under the action of $G_K$, and this action factors through $\hat{G}$ (see \cite[Definition 5.1]{GLS14} for more details).
\end{definition}

We will use $\hat{\fM}$ to denote any $(\varphi,\hat{G})$-module. For any $(\varphi,\hat{G})$-module $\hat{\fM}=(\fM,\varphi_{\fM},\hat{\mathcal{G}})$, we say $(\fM,\varphi)$ is the ambient Kisin module of $\hat{\fM}$ and we say a sequence of $(\varphi,\hat{G})$-modules is exact if the sequence of ambient Kisin modules is exact.

By Lemma 6.3 of \cite{GLS14}, for all $\ofM(r_0,\dots, r_{f-1};a)$ as above there exists a $(\varphi,\hat{G})$-module $\hofM(r_0,\dots,r_{f-1};a)$ such that $\ofM(r_0,\dots, r_{f-1};a)$ is the ambient Kisin module.

\subsection{Functors to Galois modules}

We briefly introduce some functors from Kisin modules and $(\varphi,\hat{G})$-modules to Galois representations and summarise some of their properties. Again we follow \cite{GLS14}.

\subsubsection{From Kisin modules to $G_{\infty}$-representations}

\begin{definition}
We define contravariant functors $T_{\fS}$ from $\Mod_{\fS}^{\varphi,r}$ and $\Mod_{\fS,\text{tor}}^{\varphi,r}$ to the category $\Rep_{\Zp}(G_{\infty})$ of $\Z_p[G_{\infty}]$-modules by setting:
\[
T_{\fS}(\fM):=\Hom_{\fS,\varphi}(\fM,W(R)) \text{ if $\fM$ is a finite free Kisin module},
\]
and
\[
T_{\fS}(\fM):=\Hom_{\fS,\varphi}(\fM,\Qp/\Zp \otimes_{\Zp} W(R)) \text{ if $\fM$ is a torsion Kisin module}.
\]
\end{definition}
If $\fM$ has a natural $A$-action, then $T_{\fS}(\fM)$ is an $A[G_{\infty}]$-module. We summarise some of the functor's properties:

\begin{theorem}
\label{thm32}
The functor $\Tkis$ from $\catKis$ to $\Rep_{\Zp}(G_{\infty})$ is exact and fully faithful. If $V$ is a semi-stable representation of $G_K$ with non-negative Hodge--Tate weights in some range $[0,r]$, and $L \subset V$ is a $G_K$-stable $\Zp$-lattice, there exists $\fM \in \catKis$ such that $\Tkis(\fM) \cong L|_{G_{\infty}}$. If $L$ is an $A$-module such that the $A$-action commutes with the action of $G_K$, then the Kisin module $\fM$ has a natural $A$-action. Moreover, the functor $\Tkis$ is compatible with the formation of tensor products over $\F_L$.
\end{theorem}
\begin{proof}
This is Theorem 3.2(1)-(3) and Proposition 3.4(1) of \cite{GLS14}.  The compatibility with tensor products follows from the fact that $\Tkis(\fM)$ is also the representation associated to an \'{e}tale $\varphi$-module (as in the proof of \cite[Theorem 3.3.2]{cegm}) and the statement for \'{e}tale $\varphi$-modules is Theorem 3.1.8 in \cite{brinonconrad}.
\end{proof}

We use this result to attach Kisin modules to (lattices in) our Galois representations.

\begin{definition}
Let $V$ be a semi-stable representation of $G_K$ with non-negative Hodge--Tate weights in some range $[0,r]$ with lattice $L$ as in Theorem \ref{thm32}. We say the Kisin module $\fM$ obtained from Theorem \ref{thm32} is the Kisin module attached to the lattice $L$. This is well-defined up to isomorphism.
\end{definition}


\subsubsection{From $(\varphi, \hat{G})$-modules to $\Z_p[G_{K}]$-modules}
\label{thatfunctor}
We similarly define a functor from $(\varphi,\hat{G})$-modules to Galois representations, this time to $\Z_p[G_{K}]$-modules. We do this by setting:
\[
\hat{T}(\hM):=\Hom_{\hat{\mathcal{R}},\varphi}(\hat{\mathcal{R}} \otimes_{\varphi,\fS} \fM,W(R)) \text{ if $\fM$ is a finite free Kisin module},
\]
and
\[
\hat{T}(\hM):=\Hom_{\hat{\mathcal{R}},\varphi}(\hat{\mathcal{R}} \otimes_{\varphi,\fS} \fM,\Qp/\Zp \otimes_{\Zp} W(R)) \text{ if $\fM$ is a torsion Kisin module}.
\]

We let $A$ be a finite commutative $\Z_p$-algebra. We say a $(\varphi,\hat{G})$-module $\hat{\fM}$ has a natural $A$-action if the ambient Kisin module $\fM$ has a natural $A$-action that also commutes with the $\hat{G}$-action on $\hat{\mathcal{R}}\otimes_{\varphi,\fS} \fM$. 
If $\hat{\fM}$ has a natural $A$-action, then $\hat{T}(\hat{\fM})$ is an $A[G_{K}]$-module.

We highlight some of the functor's properties:

\begin{theorem}[{\cite[Theorem 5.2]{GLS14}}]
\label{thm52}
There is a natural isomorphism $\theta: T_{\fS}(\fM) \to \hat{T}(\hM)|_{G_{\infty}}$.The functor $\hat{T}$ is exact.
\end{theorem}

\subsection{Reductions of rank one Kisin modules}

We use the theory of Kisin and $(\varphi,\hat{G})$-modules to study Galois representations. In particular, we often want to show two Galois representations are isomorphic. For rank one Kisin modules and one-dimensional representations we have the following result:

\begin{proposition}[{Proposition 6.7 from \cite{GLS14}}]
\label{prop67}
Write $\hofM=\hofM(r_0,\dots, r_{f-1};a)$ and $\hat{\overline{\mathfrak{M}'}}=\hofM(r_0',\dots, r_{f-1}';a')$. Let $\ofM$, $\overline{\mathfrak{M}'}$ denote the ambient Kisin modules of $\hofM, \hat{\overline{\mathfrak{M}'}}$.
\begin{enumerate}
\item We have $\hat{T}(\hofM)|_{I_K} \cong \omega_0^{r_0} \dots \omega_{f-1}^{r_{f-1}}$.
\item We have $\hat{T}(\hofM) \cong \hat{T}(\hat{\overline{\mathfrak{M}'}})$ if and only if $T_{\fS}(\ofM) \cong T_{\fS}(\overline{\mathfrak{M}'})$.
\item The isomorphism in (2) occurs if and only if $a=a'$ and $\sum_{i=0}^{f-1} p^{f-i-1}r_i \equiv \sum_{i=0}^{f-1} p^{f-i-1} r_i' \mod p^f-1$.
\end{enumerate}
\end{proposition}

More generally we have the following result:
\begin{lemma}[{Lemma 5.8 from \cite{GLS14}}]
\label{lem58}
Suppose that $\hat{f}:\hofM \to \hat{\ol{\mathfrak{M}'}}$ is a map of torsion $(\varphi,\hat{G})$-modules with natural $\F_L$-action, and let $\ofM, \ol{\mathfrak{M}'}$ be the ambient Kisin modules of $\hofM$ and $\hat{\ol{\mathfrak{M}'}}$. Then $\hat{T}(\hat{f})$ is injective (resp. surjective, an isomorphism) if and only if the induced map $\ofM[\tfrac{1}{u}] \to \ol{\mathfrak{M}'}[\tfrac{1}{u}]$ is injective (resp. surjective, an isomorphism).
\end{lemma}

\section{Reductions of crystalline representations}
\label{reductionsofcrysreps}
In this section we study reductions of crystalline Galois representations via Kisin modules. We use results of Gee-Liu-Savitt in \cite{GLS15} and adapt some definitions in \cite{GLS14} to extend the main results therein. In particular, our generalisation involves extending the range of Hodge--Tate weights from $[1,p]$ to $[0,p]$. This might seem like a small step, but it will allow us to deal with the case of non-distinct Hodge--Tate weights.

\subsection{Extensions of rank one $\varphi$-modules}
We are interested in studying extensions of Kisin modules. In particular, let $r_0, \dots, r_{f-1}$ be non-negative integers and fix a subset $J \subset \{0,\dots, f-1\}$. We study extensions of $\ofN=\overline{\fM}(h_0,\dots, h_{f-1};a)$ by $\overline{\mathfrak{P}}=\overline{\fM}(r_0-h_0,\dots,r_{f-1}-h_{f-1};b)$ where 
\[
h_i=\begin{cases}
r_i & \text{ if } i \in J, \\
0 & \text{ otherwise.}
\end{cases}
\]
We use a structure theorem of extensions of Kisin modules (\cite[Proposition 5.1.3]{GLS15}). By this result we can choose bases for these extensions such that we have explicit descriptions for the action of $\varphi$ on these bases. This will allow us later to pin down the extensions, and to count the possible extensions. 

The results of the structure theorem depend on whether or not there exist non-zero maps from $\ofN$ to $\ofP$. In \cite{GLS15}, the authors give a very explicit condition for this existence, using the following expression.

\begin{definition}
\label{defalpha}
Write $\ofN=\ofM(s_0, \dots, s_{f-1};a)$. Then let 
\[
\alpha_i(\ofN) := \frac{1}{p^f-1} \sum_{j=1}^f p^{f-j} (s_{j+i}). 
\]
\end{definition}

These numbers satisfy the following relation: 
\begin{align}
    \label{eq: alphaidentity}
\alpha_i(\ofN) + s_i = p \alpha_{i-1}(\ofN).
\end{align}

\begin{lemma}[Lemma 5.1.2 from \cite{GLS15}]
\label{lemma512}
Write $\ofN_1=\ofM(s_0,\dots, s_{f-1};a)$ and $\ofN_2=\ofM(s_0',\dots, s_{f-1}';a')$.
There exists a non-zero map $\ofN_1 \to \ofN_2$ if and only if 
\begin{equation}
\alpha_i(\ofN_1,\ofN_2):= \alpha_i(\ofN_1) - \alpha_i(\ofN_2) = \frac{1}{p^f-1} \sum_{j=1}^f p^{f-j} (s_{j+i}-s'_{j+i}) \in \Z_{\geq 0},
\end{equation}
for all $i$, and $a=a'$.
\end{lemma}

 We use this lemma in conjunction with the following lemma:

\begin{lemma}[Lemma 7.1 from \cite{GLS14}]
\label{lem71}
Suppose that $r_0, \dots r_{f-1}$ are integers in the range $[-p,p]$ that satisfy $\sum_{i=0}^{f-1} p^{f-1-i}r_i \equiv 0 \mod (p^f-1)$ for some integer $f$. Then either
\begin{enumerate}
\item $(r_0, \dots, r_{f-1}) = \pm (p-1, \dots, p-1)$, or,
\item the numbers $r_0,\dots r_{f-1}$, considered as a cyclic list, can be broken up into strings of the form $\pm (-1,p-1,\dots,p-1,p)$ (where there may not be any occurrences of $p-1$) and strings of the form $(0,\dots, 0)$.
\end{enumerate}
\end{lemma}

We define a set $\cP'$ of $f$-tuples. We show that if there exist non-zero maps between $\ofN \to \ofP$, that then it is necessary that $(r_0,\dots,r_{f-1}) \in \cP'$. The set $\cP'$ is a variant of the set $\mathcal{P}$ defined in \cite[Definition 7.2]{GLS14}. We make this adjustment to accommodate for possible non-distinct Hodge--Tate weights. 

Before we do, let us stress that we view any $f$-tuple $(r_0, \dots, r_{f-1})$ as indexed by elements of $\Z/f\Z$ so that $r_f=r_0$, although we usually pick representatives in $[0,f-1]$.

\begin{definition}
\label{adjustedset}
Let $\cP'$ be the set of $f$-tuples $(r_0, \dots, r_{f-1})$ with $r_i \in \{0,1,p-1,p\}$ for all $i$ such that
\begin{itemize}
\item if $r_i=p$, then $r_{i+1} \in \{0,1\}$,
\item if $r_i \in \{1, p-1 \}$, then $r_{i+1} \in \{p-1,p\}$,
\item if $r_i=0$, then either $r_i=0$ for all $i$, or there exist $s,t \geq 1$ such that 
\[
(r_i,r_{i+1},\dots, r_{i+s})=(0,\dots, 0,1), \quad (r_{i-t},\dots, r_{i-1},r_i)=(p,0, \dots, 0).
\]
\end{itemize}
\end{definition}

\begin{lemma}
\label{auxLemma2}
Let $\ofN$ and $\ofP$ be the two Kisin modules defined at the beginning of this section and suppose $0 \leq r_i \leq p$ for all $i$. Suppose there exists a non-zero map $\ofN \to \ofP$. Then we must have $(r_0,\dots, r_{f-1}) \in \mathcal{P}'$ and $\{ i \mid r_i=p-1,p \} \subset J$ and $i \not \in J$ if $r_i=1$.
\end{lemma}
\begin{proof}
We note first that by Lemma \ref{lemma512} we must have
\[
\alpha_i(\ofN,\ofP)=\frac{1}{p^f-1} \sum_{j=1}^f p^{f-j} (2h_{j+i}-r_{j+i}) \in \Z_{\geq 0}.
\]
So that in particular
\[
 \sum_{j=1}^f p^{f-j} (2h_{j+i}-r_{j+i})  \equiv 0 \mod (p^f-1).
\]
The claim that $(r_0,\dots, r_{f-1}) \in \mathcal{P}'$ follows directly from Lemma \ref{lem71}. To ease notation, let $m_i:=2h_i-r_i$. We note that
\[
m_i= \begin{cases} r_i & \text{ if } i \in J, \\
-r_i & \text{ otherwise.}
\end{cases}
\]
In particular, the same Lemma \ref{lem71} tells us that either $(m_0,\dots, m_{f-1})=\pm(p-1,\dots, p-1)$, or $(m_0,\dots, m_{f-1})$ consists of strings of the form  $\pm(-1,p-1,\dots,p-1,p)$ and $(0, \dots, 0)$. The possibility $-(p-1,\dots, p-1)$ obviously violates the conditions in Lemma \ref{lemma512}. The possibility $(p-1,\dots,p-1)$ satisfies the conditions.

The possibility of having a string of the form $-(-1,p-1,\dots,p-1,p)$ requires the existence of some $k \not \in J$ such that $r_k=p$ and thus $m_k=2h_k-r_k=-p$. Then we find:
\[
\alpha_{k-1}(\ofN,\ofP):= \frac{1}{p^f-1}\left(-p^f + \sum_{j=2}^f p^{f-j} m_{j+k-1} \right).
\]
Since $-p^f + \sum_{j=2}^f p^{f-j} m_{j+k-1}<0$, this gives a contradiction with the existence of the map by Lemma \ref{lemma512}.

So indeed we find that all remaining possibilities, including $(-1,p-1,\dots,p-1,p)$, are such that $\{ i \mid r_i=p-1,p \} \subset J$ and $i \not \in J$ if $r_i=1$.
\end{proof}

Now we are ready to prove a theorem studying extensions of Kisin modules coming from crystalline representations. This is a generalised version of \cite[Theorem 7.9]{GLS14}.

\begin{theorem}
\label{thm79}
Suppose $K/\Q_p$ is unramified and $p>2$. Let $T$ be a $G_K$-stable $\mathcal{O}_L$-lattice in a crystalline representation $V$ of $L$-dimension 2 whose $\kappa_s$-labelled Hodge--Tate weights are $\{0,r_s\}$ with $r_s \in [0,p]$ for all $s$. Let $\fM$ be the Kisin module associated to $T$ and let $\overline{\fM}:=\fM \otimes_{\mathcal{O}_L} \F_L$.
\par
Assume that the $\F_L[G_K]$-module $\overline{T}:= T/ \mathfrak{m}_L T$ is reducible. Then $\overline{\fM}$ is an extension of two $\varphi$-modules of rank one, and there exist $a,b \in \F_L^{\times}$ and a subset $J \subset \{ 0, \dots, f-1 \}$ so that $\overline{\fM}$ is as follows:
\par
Set $h_i=r_i$ if $i \in J$ and $h_i=0$ if $i \not \in J$. Then $\overline{\fM}$ is an extension of $\overline{\fM}(h_0,\dots, h_{f-1};a)$ by $\overline{\fM}(r_0-h_0,\dots,r_{f-1}-h_{f-1};b)$, and we can choose bases $e_i,f_i$ of the $\overline{\fM}_i$ so that $\varphi$ has the form
\begin{align*}
    \varphi(e_{i-1})&=(b)_i u^{r_i-h_i}e_i, \\
     \varphi(f_{i-1})&=(a)_i u^{h_i}f_i+x_ie_i, 
\end{align*}
with $x_i=0$ if $i \not \in J$ or if $r_i=0$, and $x_i \in \F_L$ constant if $i \in J$, except in the following case
\begin{itemize}
    \item $(r_0,\dots, r_{f-1}) \in \mathcal{P}'$,
    \item $\{ i \mid r_i=p-1,p \} \subset J$ and $i \not \in J$ if $r_i=1$, and,
    \item $a=b$.
\end{itemize}
If this is the case, and if all $r_i=0$, then we need to allow for a constant term for any one choice of $j \in J$. Otherwise there exists $i_0 \in J$ such that $x_i$ may be taken to be 0 if $i \not \in J$ or if $r_i=0$, to be a constant for all $i \in J$ except $i=i_0$, and to be the sum of a constant and a term of degree $p$ if $i=i_0$.
\par
Finally, $\overline{T}|_{I_K} \cong \mtwo{\prod_{i \in J} \omega_i^{r_i}}{\ast}{0}{\prod_{i \not \in J} \omega_i^{r_i}}$.
\end{theorem}

\begin{remark}
Completely analogously to Remark 7.10 from \cite{GLS14}, it follows from Proposition \ref{prop67} and Lemma \ref{auxLemma2} that the exceptional cases above can only occur if $\overline{T}$ is an extension of a character by itself.
\label{rem710}
\end{remark}

We follow the proof of \cite[Theorem 7.9]{GLS14}. For clarity we give the entire argument, which is identical to that in \cite{GLS14}, except for the first and fourth paragraph where we adjust for the extended range of Hodge--Tate weights and use the results above.

\begin{proof}
It follows from \cite[Lemma 5.5]{GLS14} that $\ofM$ is an extension of two rank one $\varphi$-modules. Then \cite[Proposition 7.8]{GLS14} guarantees that if $\ofM$ is an extension of $\overline{\mathfrak{M}'}$ by $\overline{\mathfrak{M}''}$, then $\overline{\mathfrak{M}''}$ has the form $\ofM(r_0', \dots, r_{f-1}';b)$ with $r_i' \in \{ 0, r_i \}$ for all $i$. Now taking $i \in J$ if $r_i'=0$ and $i \not \in J$ if $r_i' \geq 1$ and $r_i'=r_i $ puts $\overline{\mathfrak{M}''}$ in the correct form. Considering the determinant of $\varphi$ in \cite[Theorem 4.22]{GLS14} one finds that $\overline{\mathfrak{M}'}$ then also has the correct form. Note that if $r_i=0$ for some $i$ then there exist several choices for the set $J$ that give $\overline{\mathfrak{M}''}$ and $\overline{\mathfrak{M}'}$ the correct form.
\par
Now $\ofM$ can be taken to have the form given by Proposition 5.1.3 of \cite{GLS15}. We obtain $\deg(x_i) < h_i$, except when there exists a non-zero map $\ofM' \to \ofM''$. 
\par
Suppose a non-zero map $\ofM' \to \ofM''$ exists, then we want to show that the result still holds, upon allowing one term to be the sum of a degree $p$ term and possibly a constant (or just a constant term if all $r_i=0$). If $r_i \geq 1$ for all $i$, we note that in this case $\cP'=\cP$. It follows from \cite[Proposition 7.4]{GLS14} that we can pick $i=i_0$ such that $x_{i_0}$ is the sum of a polynomial of degree less than $h_{i_0}$ and a term of degree $p$ and $\deg(x_i)<h_i$ for all the other terms.
\par
Now suppose for some $i$ we have $r_i=0$ and that we have a non-zero map $\overline{\mathfrak{M}'} \to \overline{\mathfrak{M}''}$.  If all $r_i$ are zero, then for one choice of $j \in J$ we allow $x_j$ to be a constant term and $x_i=0$ otherwise, as follows from Proposition 5.1.3 of \cite{GLS15}. Now suppose not all $r_i=0$, but just some. Then by Definition \ref{adjustedset} and Lemma \ref{auxLemma2} there must be $i_0 \in J$ such that $r_{i_0}=p$ and $r_{i_0+1}=0$. Then by Proposition 5.1.3 of \cite{GLS15} it suffices to allow for the polynomial $x_{i_0}$ of degree $<p$ to also have an added term of degree $p$. That these are all the possible exceptional cases (which have $r_i=0$ for some $i$) follows from Lemma \ref{auxLemma2}.
\par
It remains to show that each $x_i$ with $i \in J$ cannot have any non-zero terms of degree between $1$ and $r_i-1$. But \cite[Theorem 4.22]{GLS14} implies that the image $\varphi(\ofM_{i-1}) \subset \ofM_i$ is spanned over $\F_L[[u^p]]$ by an element divisible exactly by $u^0$ and an element divisible exactly by $u^{r_i}$. On the other hand, if $x_i$ were to have a term of degree between $1$ and $r_i-1$ then neither $(b)_ie_i + \varphi(c) ((a)_i u^{r_i} f_i + x_ie_i)$ nor $(a)_i u^{r_i} f_i+x_i e_i+ \varphi(c)(b)_i e_i$ would be divisible exactly by $u^{r_i}$ for any $c \in \F_L[[u]]$. This is a contradiction.
\par
Finally, that $\overline{T}|_{I_K}$ is as claimed follows from parts (1) and (2) of Proposition Proposition \ref{prop67}, together with the fact that two mod $p$ characters of $G_K$ that are equal on $G_{\infty}$ must be equal.
\end{proof}

Just like in \cite{GLS14}, the following follows immediately from Theorem \ref{thm79}.

\begin{corollary}[Generalised Corollary 7.11]
\label{cor711}
Suppose that $K/\Q_p$ is unramified and $p>2$. Let $\rho:G_K \to \GL_2(\Fpb)$ be the reduction mod $p$ of a $G_K$-stable $\overline{\Z}_p$-lattice in a crystalline $\Qpb$-representation of dimension 2 whose $\kappa$-labelled Hodge--Tate weights are $\{0, r_{\kappa} \}$ with $r_{\kappa} \in [0,p]$ for all $\kappa$.
\par
Assume that $\rho$ is reducible. Let $\Sigma=\Hom(\F,\Fpb)$, and identify the set $\Sigma$ with $\Hom_{\Qp}(K,\Qpb)$. Then there is a subset $J \subset \Sigma$ such that
\[
\rho|_{I_K} \cong \mtwo{\prod_{\tau \in J} \omega_{\tau}^{r_{\tau}}}{\ast}{0}{\prod_{\tau \not \in J} \omega_{\tau}^{r_{\tau}}}.
\]
\end{corollary}

This corollary also is a consequence of Theorem 1.0.1 by Bartlett \cite{bartlett}, which precedes this work. This result relates Hodge--Tate weights to inertial weights for crystalline representations $\rho:G_K \to \GL_n(\Zpb)$.

\subsection{Reduction of crystalline representations and their shape}
In this subsection we use the results from the previous section to study the possible reductions of crystalline representations. This generalises two main theorems of \cite{GLS14} which serve as important input for their proof of the Buzzard--Diamond--Jarvis conjecture for unitary groups.

As well as extending the results, we need to adapt some of the definitions occurring in Section 8 and 9 of \cite{GLS14}.

\subsubsection{Types of Kisin modules}
We need a slightly more general notion of types of Kisin modules than in \cite{GLS14} as follows:

\begin{definition}[Adapted Definition 8.6 from \cite{GLS14}] 
\label{def86}
Write $\vr:=(r_0,\dots,r_{f-1})$ with $r_i \in [0,p]$ for all $i$. We say that a Kisin module $\ofM$ is an extension of type $(\vr,a,b,J)$ if it has the same shape as the Kisin modules described by Theorem \ref{thm79}; that is, $\ofM$ sits in a short exact sequence
\[
0 \to \ofM(r_0-h_0,\dots,r_{f-1}-h_{f-1};b) \to \ofM \to \ofM(h_0,\dots,h_{f-1};a) \to 0
\]
in which the extension parameters $x_i$ satisfy $x_i=0$ if $i \not \in J$ or $r_i=0$, and $x_i \in \F_L$ if $i \in J$ (except that in exceptional cases for some $i_0 \in J$,  $x_{i_0}$ is allowed to have a term of degree $p$, and if $r_i=0$ for all $i$, for one choice of $j \in J$ we have $x_j \in \F_L$ and $x_i=0$ otherwise).

\begin{remark}
\label{zerodoesnotmatter}
We note here that it follows from Theorem $\ref{thm79}$ that if there exists some $j \not \in J$ such that $r_j=0$, that then an extension is of type $(\vr,a,b,J)$ if and only if it is of type $(\vr,a,b,J \cup \{j\})$.
\end{remark}

We say a $(\varphi,\hat{G})$-module $\hofM$ with natural $\F_L$-action is of type $(\vr,a,b,J)$ if it is an extension
\[
0 \to \hat{\overline{\fM''}} \to \hat{\overline{\fM}} \to \hat{\overline{\fM'}}  \to 0
\]
such that the ambient short exact sequence of Kisin modules is an extension of type $(\vr,a,b,J)$, and if for all $x \in \fM$ there exist $\alpha \in R$ and $y \in R \otimes_{\varphi,\mathfrak{S}} \ofM$ such that $\tau(x)-x=\alpha y$ and $v_R(\alpha) \geq \tfrac{p^2}{p-1}$ where $R=\varprojlim \mathcal{O}_{\overline{K}}/p$ is a valuation ring with valuation $v_R$ and $\tau$ is as in (4.8) in \cite{GLS14}.
\end{definition}

This last condition is motivated by the following lemma:
\begin{lemma}[Adapted Lemma 8.1 from \cite{GLS14}]
\label{lem81}
Except possibly for the case that $r_i=h_i=p$ for all $i=0, \dots, f-1$, there is at most one way to extend the exact sequence
\[
0 \to \ofM(r_0-h_0,\dots,r_{f-1}-h_{f-1};b) \to \ofM \to \ofM(h_0,\dots,h_{f-1};a) \to 0
\]
to an exact sequence of $(\varphi,\hat{G})$-modules with natural $\F_L$-action such that for any $x \in \fM$ there exist $\alpha \in R$ and $y \in R \otimes_{\varphi,\mathfrak{S}} \ofM$ such that $\tau(x)-x=\alpha y$ and $v_R(\alpha) \geq \tfrac{p^2}{p-1}$. In particular, the $\hat{G}$-action on $\hofM$ is uniquely determined by $\ofM$, except possibly for the case that $r_i=h_i=p$ for all $i=0, \dots, f-1$.
\end{lemma}
\begin{proof}
We observe that the right-hand-side of equation (8.4) in the proof of Lemma 8.1 in \cite{GLS14} is still at most $\tfrac{p^2}{p-1}$ upon extending the range so that $r_i \in [0,p]$. The result then follows from the same proof.
\end{proof}


Since we have expanded the definition of the types, we need to also adjust the set $J_{\max}$ defined in \cite{GLS14}. The role of this set is as follows: if $\hofM$ is a $(\varphi,\hat{G})$-module of some type $(\vr,a,b,J)$, then there exists a $(\varphi,\hat{G})$-module $\hat{\ofN}$ of type $(\vr,a,b,J_{\max})$ whose associated Galois representation is isomorphic to that of $\hofM$.

Again fix $r_0,\dots, r_{f-1} \in [0,p]$, $J \subset \{0, \dots, f-1 \}$ and let $h_i$ be defined as in Theorem $\ref{thm79}$. We write 
\[
h(J)=\sum_{i=0}^{f-1} p^{f-1-i} h_i \mod p^f-1.
\]
Fix $h \in \Z/(p^f-1)\Z$ and suppose for some $J \subset \{0, \dots, f-1 \}$ we have $h=h(J)$. There may be multiple choices for $J$. Let $i_1, \dots, i_{\delta}$ be the distinct integers in the range $\{0,\dots, f-1\}$ such that
\begin{itemize}
    \item $(r_{i_j},\dots, r_{i_j+s_j})=(1,p-1,\dots,p-1,p)$ for some $s_j>0$, and
    \item either $i_j \in J$ and $i_j+1,\dots i_j+s_j \not \in J$ or vice versa
\end{itemize}
for all $1 \leq j \leq \delta$.

We define $J_{\max}$ as in \cite[p. 39]{GLS14} with a slight adjustment.

\begin{definition}[Adjusted definition of $J_{\max}$]
\label{jmax}
For each $h$ such that $h=h(J)$ for at least one $J$, we define $J_{\max}$ to be the unique subset of $\{0,\dots,f-1\}$ such that
\begin{itemize}
    \item $h=h(J_{\max})$, and
    \item $i_j \not \in J$ and  $i_j+1,\dots i_j+s_j \in J$ for all $1 \leq j \leq \delta$
\end{itemize}
and such that (this is our adjustment):
\[
i  \not \in J_{\max} \text{ if } r_i=0.
\]
\end{definition}

The adaptation of our definition ensures the set $J_{\max}$ is unique. In particular, with this choice and with Remark \ref{zerodoesnotmatter}, the following proposition (and its proof!) generalises to our setting. The result means that for our purposes it suffices to look at $(\varphi,\hat{G})$-modules of type $(\vr,a,b,J_{\max})$.

\begin{proposition}[Adapted Proposition 8.8 from \cite{GLS14}]
\label{prop88}
Let $\hofM$ be a $(\varphi,\hat{G})$-module of type $(\vr,a,b,J)$ and set $h=h(J)$. Then there exists a $(\varphi,\hat{G})$-module $\hat{\ofN}$ of type $(\vr,a,b,J_{\max})$ such that $\hat{T}(\hat{\ofN}) \cong \hat{T}(\hofM)$.
\end{proposition}

\subsubsection{A cohomological subspace}
In this section we define a cohomological subspace consisting of elements such that the Galois representation associated to the extension class has a crystalline lift of a certain shape. We will then study this space as we will use its properties later. Our definition of the space differs slightly from Definition 9.2 in \cite{GLS14}, as we allow for an extended range of Hodge--Tate weights. We first need the following lemma. 

\begin{lemma}
\label{cryslift}
Let $A=\{a_{\tau}\}_{\tau \in \Sigma}$ be a collection of integers.
\begin{enumerate}
    \item There is a crystalline character $\psi_A: G_K \to \overline{\Z}_p^{\times}$ such that for each $\tau \in \Sigma$ we have $\HT_{\tau} (\psi_A)=a_{\tau}$ and this character is unique up to unramified twists.
    \item We have $\overline{\psi}_A|_{I_K}=\prod_{\tau \in \Sigma} \omega_{\tau}^{a_{\tau}}$
\end{enumerate}
\end{lemma}
\begin{proof}
    Follows directly from Lemma 5.1.6 of \cite{herziggee}.
\end{proof}

Let $\liftrep: G_K \to \GL_2(\Zpb)$ be a continuous representation such that its reduction $\rho: G_K \to \GL_2(\Fpb)$ is reducible. Suppose that $\liftrep$ is crystalline with $\kappa$-Hodge--Tate weights $\{b_{\kappa,1},b_{\kappa,2} \}$ for each $\kappa \in \Hom_{\Qp}(K, \Qpb)$, and suppose further that $0 \leq b_{\kappa,1} - b_{\kappa,2} \leq p$ for each $\kappa$.
\par 
Write $\rho \cong \mtwo{\overline{\psi}_1}{\ast}{0}{\overline{\psi}_2}$, by Corollary \ref{cor711} we find that there is a decomposition $\Hom_{\Qp}(K,\Qpb)=J \prod J^c$ such that
\[
\overline{\psi}_1|_{I_K}=\prod_{\kappa \in J} \omega_{\overline{\kappa}}^{b_{\kappa,1}} \prod_{\kappa \in J^c} \omega_{\overline{\kappa}}^{b_{\kappa,2}}, \quad \overline{\psi}_2|_{I_K}=\prod_{\kappa \in J^c} \omega_{\overline{\kappa}}^{b_{\kappa,1}} \prod_{\kappa \in J} \omega_{\overline{\kappa}}^{b_{\kappa,2}},
\]
there may be several such $J$ but we temporarily fix one choice.
\par
Let $\psi_1,\psi_2: G_K \to \overline{\Z}_p^{\times}$ be crystalline lifts of $\overline{\psi}_1,\overline{\psi}_2$, with the properties that
\begin{align}
\label{eq: charpsi}
\HT_{\kappa}(\psi_1)=\begin{cases} b_{\kappa,1} & \kappa \in J, \\
b_{\kappa,2} & \text{otherwise, } \end{cases} \quad
\HT_{\kappa}(\psi_2)=\begin{cases} b_{\kappa,2} & \kappa \in J, \\
b_{\kappa,1} & \text{otherwise. } \end{cases} 
\end{align}
These exist by Lemma \ref{cryslift} and are unique up to an unramified twist.

\begin{definition}
\label{extspace}
We define $L_{\psi_1,\psi_2}$ to be the subset of $H^1(G_K, \ol{\psi}_1\ol{\psi}_2^{-1})$ consisting of all elements such that the corresponding representation has a crystalline lift of the form
\[
\mtwo{\psi_1}{\ast}{0}{\psi_2}.
\]
\end{definition}


We define $J_0=\{ \kappa \in \Hom_{\Qp}(K,\Qpb) \mid b_{\kappa,1}-b_{\kappa,2}=0\}$.

We are interested in the dimension of the subspaces $L_{\psi_1,\psi_2}$, just like in Lemma 9.3 of \cite{GLS14}. We need to adapt both the statement and the proof. The arguments are essentially identical, the only observation we need is that if $i \in J_0$, this does not contribute to the dimension of the subspace.

We first recall the definition of the space classifying crystalline extensions.

\begin{definition}
We let $V$ be a $p$-adic representation of $G_K$, then as in \cite[p. 353]{blochkato} we set
\[
H_f^1(G_K,V):= \ker(H^1(G_K,V) \to H^1(G_K, V \otimes B_{\text{cris}})).
\]
\end{definition}

Any $\alpha \in H^1(G_K,V)$ corresponds to an extension
\[
0 \to V \to E \to \Qp \to 0.
\]
Assume $V$ is crystalline, then any such $E$ is crystalline if and only if $\alpha \in H_f^1(G_K,V)$. So the space $H_f^1(G_K,V)$ classifies those extensions in which $E$ is again a crystalline representation.

\begin{lemma}[Adjusted Lemma 9.3 from \cite{GLS14}]
\label{lemma93}
$L_{\psi_1,\psi_2}$ is an $\Fpb$-vector subspace of $H^1(G_K, \overline{\psi}_1 \overline{\psi}_2^{-1})$ of dimension $|J\setminus (J \cap J_0) |$, unless $\overline{\psi}_1=\overline{\psi_2}$ in which case it has dimension $|J \setminus (J \cap J_0) |+1$.
\end{lemma}
\begin{proof}
Let $\psi=\psi_1 \psi_2^{-1}$. By \cite[Proposition 1.24(2)]{nekovar}, we see that, since $b_{\kappa,1} > b_{\kappa,2}$ if and only if $\kappa \not \in J_0$, we have 
\[
\dim_{\Qpb}(H^1_f(G_K,\Qpb(\psi))=|J \setminus (J \cap J_0)|
\]
if $\psi \neq 1$. If $\psi=1$, the dimension is $|J \setminus (J \cap J_0) |+1$.

Now write $\eta$ for the natural map $H^1(G_K,\Zpb(\psi)) \to H^1(G_K,\Qpb(\psi))$ and recall that $H^1_f(G_K,\Zpb(\psi))$ is the preimage of $H_f^1(G_K,\Qpb(\psi))$ under $\eta$. The kernel of $\eta$ is precisely the torsion part of $H^1(G_K, \Zpb(\psi))$, which we will denote $H^1(G_K, \Zpb(\psi))_{\text{tor}}$. So $H^1_f(G_K,\Zpb(\psi))$ contains this torsion part. If $\psi \neq 1$, this torsion is non-zero if and only if $\overline{\psi}=1$ in which case it has the form $\lambda^{-1}\Zpb/\Zpb$ for some $\lambda \in \mathfrak{m}_{\Zpb}$. 
\par
 Note $L_{\psi_1,\psi_2}$ is the image of $H^1_f(G_K,\Zpb(\psi))$ in $H^1(G_K,\overline{\psi})$, so that $\dim L_{\psi_1,\psi_2}=|J \setminus (J \cap J_0) |$  unless $\overline{\psi}_1=\overline{\psi}_2$. If $\overline{\psi}_1=\overline{\psi}_2$ and $\psi \neq 1$, then by the above the torsion part is one-dimensional so that in this case we obtain dimension $|J \setminus (J \cap J_0)| +1$. If $\psi=1$, then the torsion part is zero so that we obtain dimension $|J \setminus (J \cap J_0) | +1$.
\end{proof}

\subsubsection{The reducible case}

We are ready to prove a generalised version of one of the main results in \cite{GLS14}. 

\begin{theorem}[Generalised Theorem 9.1]
\label{thm91}
Suppose $p>2$. Let $\liftrep: G_K \to \GL_2(\Zpb)$ be a continuous representation such that its reduction $\rho: G_K \to \GL_2(\Fpb)$ is reducible. Suppose that $\liftrep$ is crystalline with $\kappa$-Hodge--Tate weights $\{b_{\kappa,1},b_{\kappa,2} \}$ for each $\kappa \in \Hom_{\Qp}(K, \Qpb)$, and suppose further that $0 \leq b_{\kappa,1} - b_{\kappa,2} \leq p$ for each $\kappa$.
\par
Then there is a reducible crystalline representation $\tilde{\rho}': G_K \to \GL_2(\Zpb)$ with the same $\kappa$-Hodge--Tate weights as $\liftrep$ for each $\kappa$ such that ${\rho} \cong {\rho}'$, where $\rho'$ is the reduction of $\tilde{\rho}'$.
\end{theorem}
The proof is almost identical to the proof of \cite[Theorem 9.1]{GLS14} using our generalised results and our adjusted notions of types of Kisin modules in Definition \ref{def86} and the set $J_{\max}$ from Definition \ref{jmax}.  We give it in full as we will refer back to the proof later.

\begin{proof}
As before we write $\rho \cong \mtwosmall{\overline{\psi}_1}{\ast}{0}{\overline{\psi}_2}$. If $b_{\kappa,1}-b_{\kappa,2}=p$ for each $\kappa$ and $(\overline{\psi}_1 \overline{\psi}_2^{-1})|_{I_K}=\omega$ then it follows from \cite[Lemma 9.4]{GLS14}. So assume that for some $\kappa$ either one of these conditions does not hold. After twisting we can and do assume in addition that $b_{\kappa,2}=0$ for each $\kappa$. Write $r_{\kappa}:=b_{\kappa,1}$ for each $\kappa$. 
\par
We choose a sufficiently large finite extension $L/\Q_p$ such that $\rho$ is defined over $\mathcal{O}_L$, and for each tuple of integers $\{s_{\kappa}\}$ in the range $[0,p]$ such that if $\overline{\psi}_i$ ($i=1,2$) has a crystalline lift $\psi_i$ with $\HT_{\kappa}(\psi_i)=s_{\kappa}$ for all $\kappa$, it has such a lift defined over $\mathcal{O}_L$. We fix one choice for each possible $\psi_i$ (for each possible choice of Hodge--Tate weights) and then further enlarge $L$ so that each space $H^1_f(G_K,\overline{\Z}_p(\psi_1 \psi_2^{-1}))$ is defined over $\mathcal{O}_L$.
\par
Now we will allow $\liftrep$ (and thus also $\rho$) to vary over all crystalline representations $G_K \to \GL_2(\mathcal{O}_L)$ which have $\rho \cong \mtwo{\overline{\psi}_1}{\ast}{0}{\overline{\psi}_2}$ (where $\ast$ is allowed to vary) and have $\kappa$-labelled Hodge--Tate weights $\{0, r_{\kappa}\}$ for each $\kappa$.
\par
Now by Theorem \ref{thm79} we obtain a Kisin module $\overline{\mathfrak{M}}$ which is an extension of type $(\vec{r},a,b,J)$. We let $T$ be the $G_K$-stable $\mathcal{O}_L$-lattice as in the theorem. We know $\overline{T}$ is reducible, and then by Remark 8.7(2) in \cite{GLS14}, which extends to our setting, the $(\varphi,\hat{G})$-module $\hat{\overline{\mathfrak{M}}}$ associated to $\overline{T}$ via \cite[Theorem 5.2(4)]{GLS14} is of type $(\vec{r},a,b,J)$ for some $a,b$ and $J$. Now by Proposition \ref{prop88} there exists a $(\varphi,\hat{G})$-module $\hat{\overline{\mathfrak{N}}}$ of type $(\vec{r},a,b,J_{\text{max}})$ such that $\hat{T}(\hat{\overline{\mathfrak{N}}}) \cong \rho$. It follows from Proposition \ref{prop67} that $a,b$ are uniquely determined.
\par
Again by Theorem \ref{thm79} and by the assumption that we are not in the case that $(\overline{\psi}_1 \overline{\psi}_2^{-1})|_{I_K}=\omega$ and each $r_{\kappa}=p$, we see we are not in the exceptional case of Lemma \ref{lem81}. This means we can uniquely extend our exact sequences of Kisin modules to one of $(\varphi,\hat{G})$-modules. Now by looking at the extension parameters $x_i$ we obtain from Theorem \ref{thm79}, this implies that there are at most $(|\F_L|)^{|J_{\max}|}$ isomorphism classes of $(\varphi,\hat{G})$-modules $\hat{\overline{\mathfrak{N}}}$ of type $(\vec{r},a,b,J_{\text{max}})$ and so we have (by Theorem \ref{thm79} and Remark \ref{rem710}) at most $(|\F_L|)^{|J_{\max}|}$ elements of $H^1(G_K,\overline{\psi}_1 \overline{\psi}_2^{-1})$ corresponding to representations $\rho$, unless $\overline{\psi}_1=\overline{\psi_2}$, in which case $(|\F_L|)^{|J_{\max}|}$ must be replaced with $(|\F_L|)^{|J_{\max}|+1}$.
\par
Now by Lemma \ref{cryslift} we may choose crystalline characters $\psi_1,\psi_2$ lifting $\overline{\psi}_1,\overline{\psi}_2$ respectively such that 
\[
\HT_{\kappa}(\psi_1)=\begin{cases} r_{\kappa} & \kappa \in J_{\max}, \\
0 & \text{otherwise, } \end{cases} \quad
\HT_{\kappa}(\psi_2)=\begin{cases} 0 & \kappa \in J_{\max}, \\
r_{\kappa} & \text{otherwise. } \end{cases}
\]
By our choice of $L$ we may further suppose $\psi_1,\psi_2$ and $H_f^1(G_K,\Zpb(\psi_1 \psi_2^{-1}))$ are all defined over $\mathcal{O}_L$. 
\par
Recall that by construction $J_{\max}$ contains no element of $J_0$. By Lemma \ref{lemma93} there are $(\# \F_L)^{|J_{\max}|}$ extension classes which arise as the reductions of crystalline representations which are extensions of $\psi_2$ by $\psi_1$, unless $\overline{\psi}_1=\overline{\psi}_2$ in which case there are $(\# \F_L)^{|J_{\max}|+1}$  extension classes. Since we have already shown that there are at most $(\# \F_L)^{|J_{\max}|}$ (or $(\# \F_L)^{|J_{\max}|+1}$ if $\overline{\psi}_1=\overline{\psi}_2$) extension classes arising from the reduction of crystalline representations with $\kappa$-labelled Hodge--Tate weights $\{0, r_{\kappa} \}$, the result follows. 
\end{proof}

We have the following corollary:

\begin{corollary}[Corollary about extension classes]
\label{corextclass}
Suppose that $K/\Q_p$ is unramified and $p>2$. Let $\rho: G_K \to \GL_2(\Fpb)$ be the reduction mod $p$ of a $G_K$-stable $\overline{\Z}_p$-lattice in a crystalline $\Qpb$-representation of dimension 2 whose $\kappa$-labelled Hodge--Tate weights are $\{b_{\kappa,1}, b_{\kappa,2} \}$ with $0 \leq b_{\kappa,1} - b_{\kappa,2} \leq p$ for all $\kappa$.
\par
Assume that $\rho$ is reducible. Let $\Sigma=\Hom(\F,\Fpb)$, then there is a subset $J \subset \Sigma$ such that
\[
\rho|_{I_K} \cong \mtwo{\prod_{\tau \in J} \omega_{\tau}^{b_{\tau,1}} \prod_{\tau \not \in J} \omega_{\tau}^{b_{\tau,2}}}{\ast}{0}{\prod_{\tau \in J} \omega_{\tau}^{b_{\tau,2}} \prod_{\tau \not \in J} \omega_{\tau}^{b_{\tau,1}}},
\]
with $\ast \in L_{\psi_1,\psi_2}$.
\end{corollary}

\begin{proof}
The decomposition follows from Corollary \ref{cor711} (after twisting). The statement about the extension classes follows from the definition of $L_{\psi_1,\psi_2}$ and the proof of Theorem \ref{thm91}.
\end{proof}

\begin{remark}
\label{diffxi}
We can see from the proof of Theorem \ref{thm91} that different choices of the extension parameters $x_i$ in Theorem \ref{thm79} give rise to different Galois representations.
\end{remark}

\subsubsection{The irreducible case}
We now prove the analogous result in the case that $\rho$ is irreducible, which we will deduce from the reducible case. We do this by simply using the fact that an irreducible $\rho$ becomes reducible after restriction to an unramified quadratic extension.

Recall that we write $K'$ for the quadratic unramified extension of $K$, $\F'$ for the residue field and $\Sigma'=\{ \tau: \F' \hookrightarrow \Fpb \}$. We say that $J \subset \Sigma'$ is a balanced subset if it consists of precisely one element of $\Sigma'$ extending each element of $\Sigma$.

\begin{theorem}
\label{thm101}
Let $\liftrep:G_K \to \GL_2(\Qpb)$ be a continuous irreducible representation such that its reduction $\rho: G_K \to \GL_2(\Fpb)$ is also irreducible. Suppose that  $\liftrep$ is crystalline with $\kappa$-Hodge--Tate weights $\{b_{\kappa,1},b_{\kappa,2}\}$ for each $\kappa \in \Hom(K,\Qpb)$ such that $0 \leq b_{\kappa,1}-b_{\kappa,2} \leq p$. Then there is a balanced subset $J \subset \Sigma'$ such that
\[
\rho|_{I_K} \cong \mtwo{\prod_{\sigma \in J} \omega_{\sigma}^{b_{\sigma|_{\F,1}}} \prod_{\sigma \not \in J} \omega_{\sigma}^{b_{\sigma|_{\F,2}}}}{0}{0}{\prod_{\sigma \in J} \omega_{\sigma}^{b_{\sigma|_{\F,2}}} \prod_{\sigma \not \in J} \omega_{\sigma}^{b_{\sigma|_{\F,1}}}},
\]
where $b_{\sigma,i}=b_{\kappa,i}$ if $\sigma \in \Sigma$ is the reduction of $\kappa$.
\end{theorem}

This is the analogue of Theorem 10.1 in \cite{GLS14}. Like in the previous sections we follow the original proof in \cite{GLS14}, but in this case we need to make some more adjustments to account for the extended range.

\begin{proof}
 By Corollary \ref{cor711} we have a decomposition as in the statement for some $J \subset \Sigma'$ . The point is showing that we can find such a decomposition with $J$ balanced. Now since $\rho|_{I_K}$ is irreducible, we must have
\[
\prod_{\sigma \in J} \omega_{\sigma}^{b_{\sigma|_{\F,1}}} \prod_{\sigma \not \in J} \omega_{\sigma}^{b_{\sigma|_{\F,2}}}=\prod_{\sigma \in J} \omega_{\sigma  \circ \Frob^f}^{b_{\sigma|_{\F,2}}}\prod_{\sigma \not \in J} \omega_{\sigma  \circ \Frob^f}^{b_{\sigma|_{\F,1}}}.
\]
Write $J_1$ for the set of places in $\Sigma$ both of whose extensions to $\Sigma'$ are in $J$, and $J_2$ for the set of places in $\Sigma$ neither of whose extensions to $\Sigma'$ are in $J$. Then we see that we have
\[
\prod_{\sigma \in J_1} \omega_{\sigma}^{b_{\sigma,1}-b_{\sigma,2}}=\prod_{\sigma \in J_2} \omega_{\sigma}^{b_{\sigma,1}-b_{\sigma,2}}.
\]
If both $J_1,J_2$ are empty, then $J$ is balanced and we are done. So assume this is not the case.
\par
 Next we let $J_0'=\{ \sigma \in \Sigma' \, \vert \, b_{\sigma|_{\F},1} - b_{\sigma|_{\F},2} =0 \}$. Then we first redefine $J$ to be $J \setminus J_0'$, it is straightforward to see we can do this. 
\par
Now we set $x_{\sigma}=b_{\sigma,1}-b_{\sigma,2}$ if $\sigma \in J_1$, $x_{\sigma}=b_{\sigma,2}-b_{\sigma,1}$ if $\sigma \in J_2$ and $x_{\sigma}=0$ otherwise. Note that for any $\sigma \in J_0$ we automatically have $x_{\sigma}=0$. Note also that since $\rho$ is irreducible, there is at least one embedding $\sigma$ with $x_{\sigma}=0$. By construction we have $\prod_{\sigma \in \Sigma} \omega_{\sigma}^{x_{\sigma}}=1$ and each $x_{\sigma} \in [-p,p]$.  
\par
Choose an element $\sigma_0 \in \Sigma$ and recursively define $\sigma_i=\sigma_{i+1}^p$, we write $\omega_{\sigma_{i+1}}=\omega_{i+1}$ and have $\omega_{i+1}^p=\omega_i$. Now identify $\Sigma$ with $\{0, \dots, f-1\}$ by identifying $\sigma_i$ with $i$. Then by Lemma \ref{lem71} the cyclic set of those $i$ with $x_i \neq 0$ must break up as a disjoint union of sets of the form $(i,i+1,\dots, i+j)$ with $(x_i,x_{i+1},\dots,x_{i+j})=\pm(-1,p-1,p-1,\dots,p-1,p)$. For each such interval $(i,i+j)$, we may choose a lift of $i$ to $\Sigma'$, say $\sigma_i'$. We replace $J$ with $J \Delta \{i,\dots, i+j\}$, i.e. we replace $J$ with $J \Delta \{ \Frob^{-k} \circ \sigma_i' \}$ for $k=0,\dots j$.
\par
Now we have ended up with $J \subset \Sigma'$ containing exactly one element of $\Sigma'$ extending each element $\sigma \in \Sigma$ such that $b_{\sigma,1} - b_{\sigma,2} \neq 0$. 
\par
If for some $\sigma$ we have $b_{\sigma,1} - b_{\sigma,2} = 0$, then any choice of its (two) extensions $\sigma'$ will satisfy $b_{\sigma'|_{\F},1} - b_{\sigma'|_{\F},2} =0$. So for each such $\sigma$ all we have to do is to just pick one such extension to add into $J$. This way we end up with a balanced subset.
\end{proof}

\subsection{Extensions of Kisin modules and Galois representations}
The goal of this subsection is to establish some results which will allow us to relate different spaces of extension classes $L_{\psi_1,\psi_2}$ and $L_{\psi'_1,\psi'_2}$ later. We do this via the theory of Kisin modules.

Let $\ol{\psi}=\ol{\psi}_1 \ol{\psi}_2^{-1}$ for some characters $\ol{\psi}_1, \ol{\psi}_2: G_K \to \F_L^{\times}$. Throughout we will identify $H^1(G_K, \Fpb(\ol{\psi}))$ with $\Ext_{G_K}(\ol{\psi}_2,\ol{\psi}_1)$.

We will need the following lemma: 

\begin{lemma}[{\cite[Lemma 5.4.2]{GLS15}}]
\label{lem542}
Let $\ol{\psi}: G_K \to \F_L^{\times}$ be a continuous character. If $\ol{\psi} \neq \omega$, then the restriction map
\[
H^1(G_K,\Fpb(\ol{\psi})) \to H^1(G_{\infty},\Fpb(\ol{\psi}))
\]
is injective.
\end{lemma}

\subsubsection{Twists of extension classes}
\label{twistsubsection}
We will next study what happens to extension classes after twisting, first in the setting of Galois representations.

Suppose we have a representation $\rho: G_K \to \GL_2(\Fpb)$, and suppose
\[
\rho \sim \mtwo{\ol{\psi}_1}{\ast}{0}{\ol{\psi}_2},
\]
is such that $c_{\rho} \in L_{\psi_1,\psi_2}$. This means that the representation corresponding to $c_{\rho}$ has a crystalline lift of the form
\[
\mtwo{{\psi}_1}{\ast}{0}{{\psi}_2}.
\]
Suppose $\ol{\psi}: G_K \to \overline{\F}_p^{\times}$ is such that $\ol{\psi}_{I_K}=\prod_{i} \omega_i^{\delta_i}$ for some integers $\delta_i$. Then $\ol{\psi} \otimes \rho$ has a crystalline lift of the form
\[
\mtwo{{\psi}'_1}{\ast}{0}{{\psi}'_2},
\]
with Hodge--Tate weights 
\[
\HT_i(\psi'_1)=\HT_i(\psi_1)+\delta_i, \quad \HT_i(\psi'_2)=\HT_i(\psi_2)+\delta_i.
\]
We can achieve this by twisting by a crystalline lift of $\ol{\psi}$, which exists uniquely up to an unramified twist by Lemma \ref{cryslift}. 

More generally, suppose we have a pair of characters $M_1,M_2: G_K \to \overline{\F}_p^{\times}$ and an extension $M \in \Ext_{G_{\infty}}^1(M_1|_{G_{\infty}},M_2|_{G_{\infty}})$. Let $M'$ be another $G_K$-character and let $M_1'=M_1 \otimes M'$ and $M_2'=M_2 \otimes M'$. Then we obtain a isomorphism between $\Ext_{G_{\infty}}^1(M_1|_{G_{\infty}},M_2|_{G_{\infty}})$ and $\Ext_{G_{\infty}}^1(M_1'|_{G_{\infty}},M_2'|_{G_{\infty}})$ by twisting any extension by $M'|_{G_{\infty}}$.

Next let $\fM_1, \fM_2$ be a pair of rank one Kisin modules and let $\fM \in \Ext^1(\fM_1,\fM_2)$ be an extension of $\fM_1$ by $\fM_2$. We then obtain a short exact sequence
\[
0 \to \fM_1 \to \fM \to \fM_2 \to 0.
\]
Now let $\fM'$ be any rank one Kisin module as in Definition \ref{def61} and let $\fM_1'=\fM_1 \otimes \fM'$ and $\fM_2'=\fM_2 \otimes \fM'$. Then we obtain an homomorphism from $\Ext^1(\fM_1,\fM_2)$ to $\Ext^1(\fM'_1,\fM'_2)$ by twisting any $\fM \in \Ext^1(\fM_1,\fM_2)$ by $\fM'$.

\subsubsection{Twists of Kisin modules and Galois representations}
We now wish to relate twisting our Galois extensions to twists of Kisin modules. First recall we have a functor $T_{\fS}$ from Kisin modules to $G_{\infty}$-modules. By Theorem \ref{thm32} the functor $T_{\fS}$ is exact and we obtain a commutative diagram as in Figure 1,
\par
\begin{figure}[h!]
\label{comdiag}
\center
\begin{tikzcd}
\Ext^1(\fM_1,\fM_2) \arrow[r]\arrow[d]
& \Ext^1(\fM'_1,\fM'_2) \arrow[d] \\
\Ext^1_{G_{\infty}}(T_{\fS}(\fM_1),T_{\fS}(\fM_2)) \arrow[r,"\sim"]
& \Ext^1_{G_{\infty}}(T_{\fS}(\fM'_1),T_{\fS}(\fM'_2))
\end{tikzcd}
\caption{Commutative diagram of extension spaces.}
\end{figure}
\noindent where the isomorphism on the bottom is obtained by twisting any extension by $T_{\fS}( \fM')$.

\section{Statements of the main results}

Before stating the main theorems, we need some terminology. The following is taken from \cite{DS}.

\begin{definition}
We say a representation $\rho: G_K \to \GL_2(\Fpb)$ has a crystalline lift of weights $(k,l) \in \Z_{\geq 1}^{\Sigma} \times \Z^{\Sigma}$ if there exists a continuous representation $\tilde{\rho}: G_K \to \GL_2(\mathcal{O}_L)$ such that $\tilde{\rho} \otimes_{\mathcal{O}_L} \F_L \cong \rho$ and $\tilde{\rho} \otimes_{\mathcal{O}_L} L$ is crystalline with Hodge--Tate weights $\{k_{\kappa}+l_{\kappa}-1,l_{\kappa}\}_{\kappa}$ where $\kappa \in \Hom(K,\Qpb)$.
\end{definition}

We will mostly have $l_{\kappa}=0$ (which can be achieved by twisting).

\begin{definition}
We say a crystalline lift of weight $(k,l)$ is \emph{regular} if $k_{\kappa} \geq 2$ for all $\kappa \in \Hom(K, \Qpb)$. Otherwise we say the lift is \emph{irregular}.
\end{definition}

We will next assign multiple regular weights to a fixed irregular weight. 

\subsection{The regular weights}
\label{algweights}
Suppose $(k,0) \in \Z_{\geq 1}^{\Sigma} \times \Z^{\Sigma}$ is an irregular weight. From this irregular weight we obtain regular weights by increasing the weight using the analogues of weight shifts by partial Hasse invariants and partial Theta operators invariants as in \cite{DS} and \cite{hw2022}. For any  $\tau: \F \to \Fpb$ we let $h_{\tau}=p e_{\Frob^{-1} \circ \tau} - e_{\tau}$ and $\theta_{\tau}=p e_{\Frob^{-1} \circ \tau} + e_{\tau}$ where $e_{\tau}$ is the basis vector of $\Z^{\Sigma}$ corresponding to $\tau$.

Before we describe our regular weights we define a subset of the set of embeddings. Recall that we have identified the set of embeddings $\Sigma = \{\kappa: K \hookrightarrow \Qpb\}$ with the set of embeddings $\{\tau: \F \to \Fpb\}$. 

\begin{definition}
We let $M$ be the subset of $\Sigma$ consisting of $\tau \in \Sigma$ such that
\[
 k_{\Frob^{-1} \circ \tau}=\dots =k_{\Frob^{(1-s)} \circ \tau}=2 \text{ and } k_{\Frob^{-s} \circ \tau}=1,
\]
for some $s \geq 1$ (the first condition is vacuous if $s=1$).
\end{definition}

\subsubsection*{The definition of $(k',l')$}

\begin{definition}
    
Let $k \in \Z_{\geq 1}^{\Sigma}$, then we define $k' \in \Z_{\geq 2}^{\Sigma}$ as follows
\begin{align}
\label{defkprime}
k':=k+\sum_{\tau \in M} h_{\tau}.
\end{align}
We set $l'=l=(0, \dots,0)$. Using the definition of $h_{\tau}$ one can easily check $(k',l')$ is indeed regular.
\end{definition}

Before defining the next regular weights we define a subset $\tilde{M}$ of $\Sigma$. 

\begin{definition}
\label{deftildeM}
Let $\tilde{M}$ be the set of $\tau \in \Sigma$ such that one of the following conditions holds:
\begin{itemize}
\item  $k_{\tau} \geq 3$ and $\tau \in M$,
\item $k_{\tau}=2, k_{\Frob^{-1} \circ \tau}=1, k_{\Frob^{-2} \circ \tau}=1$ or $2$ and if $k_{\Frob^{-2} \circ \tau}=2$ then $\Frob^{-2} \circ \tau \in M$.
\end{itemize}
\end{definition}

In the case that there is no $\tau$ such that $k_{\tau}=2$ and $k_{\Frob^{-1} \circ \tau}=1$, then $\tilde{M}$ is just the set $\tau \in \Sigma$ such that $k_{\tau} \neq 1, k_{\Frob^{-1} \circ \tau}=1$.

\subsubsection*{The definition of $(k^{\mu}, l^{\mu})$}

\begin{definition}
Next let $\mu \in \tilde{M}$, then we define $k^{\mu}$ by setting
\begin{align}
\label{defkmu}
k^{\mu}=k+ \theta_{\mu} + \sum_{\tau \in M \setminus \{\mu\}} h_{\tau}.
\end{align}
We set $l^{\mu}= - e_{\mu}$. Again it is straightforward to see that the newly defined weight $(k^{\mu}, l^{\mu})$ is regular. 
\end{definition}

\subsubsection*{The definition of $(k^{\Theta}, l^{\Theta})$}

\begin{definition}
We define $k^{\Theta}$ by setting
\begin{align}
\label{defktheta}
k^{\Theta}=k+ \sum_{\tau \in \tilde{M}} \theta_{\tau} + \sum_{\tau \in M \setminus \tilde{M}} h_{\tau}.
\end{align}
We set $l^{\Theta}= \sum_{\tau \in \tilde{M}} - e_{\tau}$. Again it is straightforward to see that the newly defined weight $(k^{\Theta}, l^{\Theta})$ is regular. 
\end{definition}

\begin{example}
\label{explicitweight}
Suppose $f=4$ and that we have an irregular weight $(k,l)=((1,k_1,1,k_3),(0,0,0,0))$. If $k_1, k_3 \geq 3$, then we obtain weights 
\begin{align*}   
 (k',l')&=((p+1,k_1-1,p+1,k_3-1),(0,0,0,0)), \\(k^{\mu_1},l^{\mu_1})&=((p+1,k_1+1,p+1,k_3-1),(0,-1,0,0)), \\
 (k^{\mu_2},l^{\mu_2})&=((p+1,k_1-1,p+1,k_3+1),(0,0,0,-1)), \\(k^{\Theta},l^{\Theta})&=((p+1,k_1+1,p+1,k_3+1),(0,-1,0,-1)). 
\end{align*}
Now if say $k_1=2$ and $k_3 \geq 3$, then we would have
\begin{align*}   
 (k',l')&=((p,p+1,p+1,k_3-1),(0,0,0,0)), \\ (k^{\mu_1},l^{\mu_1})&=((p,p+3,p+1,k_3-1),(0,-1,0,0)), \\
 (k^{\mu_2},l^{\mu_2})&=((p,p+1,p+1,k_3+1),(0,0,0,-1)), \\(k^{\Theta},l^{\Theta})&=((p,p+3,p+1,k_3+1),(0,-1,0,-1)),
\end{align*}
in particular, the gap between the Hodge--Tate weights at the embedding $\mu_1$ for the second and fourth weight would be larger than $p$.
\end{example}

\subsection{The statement}
We have two versions of the theorem. The first one is simpler as it involves fewer weights, the second is more complicated but, as alluded to before, has as its main benefit that there is a direct analogous statement concerning geometric modularity of Galois representations (\cite{hw2022}).
\label{statementsect}

\begin{theorem}
\label{twoweightversion}
Let $p$ be an odd prime. Let $(k,0) \in \Z^{\Sigma} \times \Z^{\Sigma}$ be an irregular weight satisfying $1 \leq k_{\tau} \leq p$ for all $\tau \in \Sigma$. Suppose there is no $\tau \in \Sigma$ such that $k_{\tau}=2$ and $k_{\Frob^{-1} \circ \tau}=1$. Suppose further $k \neq 1$. Then a representation $\rho: G_K \to \GL_2(\Fpb)$ has a crystalline lift of weight $({k},{0})$ if and only if $\rho$ has  crystalline lifts of weights $(k',l')$ and $(k^{\Theta},l^{\Theta})$.
\end{theorem}
\begin{proof}
Theorem \ref{twoproof}, Proposition \ref{irreducibletoregular} and Proposition \ref{irreducibletoirregular}.
\end{proof}

\begin{theorem}
\label{conjcris}
Let $p$ be an odd prime. Let $(k,0) \in \Z^{\Sigma} \times \Z^{\Sigma}$ be an irregular weight satisfying $1 \leq k_{\tau} \leq p$ for all $\tau \in \Sigma$. Suppose there is no $\tau \in \Sigma$ such that $k_{\tau}=2$ and $k_{\Frob^{-1} \circ \tau}=1$. Suppose further $k \neq 1$. Then a representation $\rho: G_K \to \GL_2(\Fpb)$ has a crystalline lift of weight $({k},{0})$ if and only if
\begin{enumerate}
\item $\rho$ has a crystalline lift of weight $(k',l')$, and,
\item for each $\mu \in \tilde{M}$, $\rho$ has a crystalline lift of weight $({k}^{\mu},{l}^{\mu})$.
\end{enumerate}
\end{theorem}
\begin{proof}
Theorem \ref{reduciblethm}, Proposition \ref{irreducibletoregular} and Proposition \ref{irreducibletoirregular}.
\end{proof}

Let us note that both theorems hold just as well for $(k, l)$, after twisting, with the regular weights again obtained from the formulae \eqref{defkprime}, \eqref{defkmu} and \eqref{defktheta}, respectively.

We explain the reason for the `(2,1) condition', i.e. the condition that there is no $\tau$ such that $k_{\tau}=2$ and $k_{\Frob^{-1} \circ \tau}=1$. In Example \ref{explicitweight} above one can see that in this case gaps between the Hodge--Tate weights corresponding to the regular weights are larger than $p$. In order to use the techniques of \cite{GLS14} we need them to be at most $p$. 

However, the condition that there is no $\tau$ such that $k_{\tau}=2$ and $k_{\Frob^{-1} \circ \tau}=1$ could be removed by slightly adapting the weights  $(k^{\Theta},l^{\Theta})$ by replacing the set $M \setminus \tilde{M}$ in the definition by the set  $M \setminus \tilde{M}_{2}$ where $\tilde{M}_{2}$ is the set of embeddings $\tau \in \Sigma$ such that 
\begin{enumerate}
    \item $\tau \in \tilde{M}$, or, 
    \item $\Frob^{-1} \circ \tau \in \tilde{M}$ with $k_{\Frob^{-1} \circ \tau}=2$. 
\end{enumerate}
Write $(k^{\Theta, \text{alt}},l^{\Theta, \text{alt}})$ for this alternative weight. Note that we have $(k^{\Theta, \text{alt}},l^{\Theta, \text{alt}}) = (k^{\Theta},l^{\Theta})$ if there is no $\tau$ such that $k_{\tau}=2$ and $k_{\Frob^{-1} \circ \tau}=1$. 

We continue with our earlier example:

\begin{example}[Example \ref{explicitweight} continued]
Again taking $f=4$ and irregular weight $(k,l)=((1,2,1,k_3),(0,0,0,0))$ with $k_3 \geq 3$, we would have
\begin{align*}   
(k',l')&=((p,p+1,p+1,k_3-1),(0,0,0,0)), \\
(k^{\Theta, \text{alt}},l^{\Theta, \text{alt}})&=((p+1,3,p+1,k_3+1),(0,-1,0,-1)),
\end{align*}    
so the gaps between the Hodge--Tate weights would not exceed $p$.
\end{example}

We then formulate the analogous statement as a conjecture:

\begin{conjecture}
Let $\rho: G_K \to \GL_2(\Fpb)$ be a Galois representation. Let $(k,0)$ be an irregular weight such that $1 \leq k_{\tau} \leq p$ for all $\tau \in \Sigma$.  Suppose further $k \neq 1$. Then a representation $\rho: G_K \to \GL_2(\Fpb)$ has a crystalline lift of weight $({k},{0})$ if and only if $\rho$ has a crystalline lift of weight $(k',l')$ and $\rho$ has a crystalline lift of weight $({k}^{\Theta,\text{alt}},{l}^{\Theta, \text{alt}})$.
\end{conjecture}

\begin{remark}
\label{condition2-1}    
For these adapted weights, the methods of proof are the same but including these cases would require many more cases and many more pages because the weight combinatorics would be less uniform. We can make similar adjustments to the weights $(k^{\mu},l^{\mu})$ to remove the (2,1)-condition in that version of the theorem too. However, with these adjustments we would not have the direct analogy with the geometric modularity version of our theorem. So we will leave the above as a conjecture and will go ahead with our earlier defined weights and will leave the (2,1)-condition as is. 
\end{remark}

\subsection{Notation for the weights}
\label{weightnotation}
We continue with a fixed irregular weight $(k,0)$. We assume that $k \neq 1$ and that there is no $\tau$ such that $k_{\tau}=2$ and $k_{\Frob^{-1} \circ \tau} \neq 1$. We write $\{b_{\tau,1},b_{\tau,2}\}_{\tau \in \Sigma}$ for the corresponding Hodge--Tate weights. For a fixed subset $J \in \Sigma$, write
\[
\wtop_{\tau}= \begin{cases}
    b_{\tau,1} & \tau \in J, \\
    b_{\tau,2} & \tau \not \in J, \\
\end{cases}
\quad
\wbottom_{\tau}= \begin{cases}
    b_{\tau,2} & \tau \in J, \\
    b_{\tau,1} & \tau \not \in J. \\
\end{cases}
\]
Consider $(k',l')$ as in Section \ref{algweights} and write $\{b'_{\tau,1},b'_{\tau,2}\}_{\tau \in \Sigma}$ for the corresponding Hodge--Tate weights. Similarly, for each $\mu \in \tilde{M}$, consider $(k^{\mu},l^{\mu})$ and write $\{b^{\mu}_{\tau,1},b^{\mu}_{\tau,2}\}_{\tau \in \Sigma}$ for the corresponding weights. Finally we do the same for the weight $(k^{\Theta}, l^{\Theta})$. Recall that $J_0 = \{ \tau \in \Sigma \, | \, k_{\tau}=1 \}$.

This gives 
\begin{align}
\label{eq: bprime}
b'_{\tau,1}= \begin{cases} k_{\tau} - 1, & \tau \not \in M \cup J_0, \\
                          k_{\tau} - 2, & \tau \in \tilde{M}, \\
                          p-1,         & \tau \in J_0 \cap M, \\
                          p,           & \tau \in J_0, \tau \not \in M,    
\end{cases}
\quad \quad
b'_{\tau,2}=0.
\end{align}
We further have 
\begin{align}
\label{eq: bmu}
b^{\mu}_{\tau,1} = \begin{cases} k_{\tau} - 1, & \tau \not \in M \cup J_0, \\
                          k_{\tau} - 1, & \tau \in \tilde{M}, \tau = \mu, \\
                          k_{\tau} - 2, & \tau \in \tilde{M}, \tau \neq \mu, \\
                          p-1,          & \tau \in J_0 \cap M, \\    
                          p,            & \tau \in J_0, \tau \not \in M,    
\end{cases}
\quad \quad
b^{\mu}_{\tau,2}= \begin{cases}
    -1, & \tau = \mu, \\
     \,0, & \tau \neq \mu,
\end{cases}
\end{align}
and
\begin{align}
\label{eq: btheta}
b^{\Theta}_{\tau,1} = \begin{cases} k_{\tau} - 1, & \tau \not \in J_0, \\
                          p-1,          & \tau \in J_0 \cap M, \\    
                          p,           & \tau \in J_0, \tau \not \in M,    
\end{cases}
\quad \quad
b^{\Theta}_{\tau,2}= \begin{cases}
    -1, & \tau \in \tilde{M}, \\
     \,0, & \tau \not \in \tilde{M}.
\end{cases}
\end{align}

For fixed subsets $J'$, $J^{\mu}$ and $J^{\Theta}$ of $\Sigma$, we then set
\begin{align}
\label{eq: rexpressions}    
\wtop'_{\tau}=\begin{cases}
    b'_{\tau,1} & \tau \in J', \\
    b'_{\tau,2} & \tau \not \in J',     
\end{cases}
\quad \quad
\wtop^{\mu}_{\tau}=\begin{cases}
    b^{\mu}_{\tau,1} & \tau \in J^{\mu}, \\
    b^{\mu}_{\tau,2} & \tau \not \in J^{\mu},  
\end{cases}
\quad \quad
\wtop^{\Theta}_{\tau}=\begin{cases}
    b^{\Theta}_{\tau,1} & \tau \in J^{\Theta}, \\
    b^{\Theta}_{\tau,2} & \tau \not \in J^{\Theta},    
\end{cases}
\end{align}
and similarly
\begin{align}
\label{eq: sexpressions}    
\wbottom'_{\tau}=\begin{cases}
    b'_{\tau,2} & \tau \in J', \\
    b'_{\tau,1} & \tau \not \in J',     
\end{cases}
\quad \quad
\wbottom^{\mu}_{\tau}=\begin{cases}
    b^{\mu}_{\tau,2} & \tau \in J^{\mu}, \\
    b^{\mu}_{\tau,1} & \tau \not \in J^{\mu},    
\end{cases}
\quad \quad
\wbottom^{\Theta}_{\tau}=\begin{cases}
    b^{\Theta}_{\tau,2} & \tau \in J^{\Theta}, \\
    b^{\Theta}_{\tau,1} & \tau \not \in J^{\Theta}.    
\end{cases}
\end{align}

\section{Matching for reducible representations}
In this section we will assume a reducible representation $\rho: G_K \to \GL_2(\Fpb)$ has a crystalline lift of some weights, and will investigate the implications on the shape of $\rho|_{I_K}$. At the end of the section we can prove our main results for $\rho$ semisimple.

Throughout the remaining sections we fix an irregular weight $(k,0)$. We assume that $1 \leq k_{\tau} \leq p$ for all $\tau$ and moreover that there is no $\tau$ such that $k_{\tau}=2$ and $k_{\Frob^{-1} \circ \tau}=1$. Finally we assume $k \neq 1$. 

To ease notation, in this section we sometimes use integers as indices, e.g. $\wtop_i = \wtop_{\tau_i}$, and, as before,  sometimes we write $i \in J$ to mean $\tau_i \in J$. Recall that we have labelled our embeddings $\tau_0, \dots, \tau_{f-1}$ with the convention that $\tau_i=\tau_{i+1}^{p}$.


\subsection{Consequences of having regular lifts}
Suppose a reducible representation $\rho: G_K \to \GL_2(\Fpb)$ has a crystalline lift of regular weights $(k',l')$ and $\{(k^{\mu},l^{\mu})\}_{\mu \in \tilde{M}}$. By Corollary \ref{corextclass} we have subsets $J'$ and $J^{\mu}$ such that
    \begin{align*}    
    \rho|_{I_{K}} & \sim \mtwo{\prod_{\tau \in J'} \omega_{\tau}^{b'_{\tau,1}} \prod_{\tau \not \in J'} \omega_{\tau}^{b'_{\tau,2}}}{\ast}{0}{\prod_{\tau \in J'} \omega_{\tau}^{b'_{\tau,2}} \prod_{\tau \not \in J'} \omega_{\tau}^{b'_{\tau,1}}} \\
    & \sim \mtwo{\prod_{\tau \in J^{\mu}} \omega_{\tau}^{b^{\mu}_{\tau,1}} \prod_{\tau \not \in J^{\mu}} \omega_{\tau}^{b^{\mu}_{\tau,2}}}{\ast}{0}{\prod_{\tau \in J^{\mu}} \omega_{\tau}^{b^{\mu}_{\tau,2}} \prod_{\tau \not \in J^{\mu}} \omega_{\tau}^{b^{\mu}_{\tau,1}}},
    \end{align*}
for each $\mu \in \tilde{M}$.

Then, using the notation from the previous section, this gives us
    \begin{align}    
    \label{eq: rhoinertprime}
    \rho|_{I_{K}} & \sim \mtwo{\prod_{\tau} \omega_{\tau}^{\wtop'_{\tau}}}{\ast}{0}{\prod_{\tau} \omega_{\tau}^{\wbottom'_{\tau}}} \\
    & \sim \mtwo{\prod_{\tau} \omega_{\tau}^{\wtop^{\mu}_{\tau}}}{\ast}{0}{\prod_{\tau} \omega_{\tau}^{\wbottom^{\mu}_{\tau}}}, 
      \label{eq: rhoinertmu}
    \end{align}
for each $\mu \in \tilde{M}$. So that 
\begin{align}
\label{eq: congstring}
    \sum_{i \in \{0, \dots f-1\}} \wtop'_{i} \cdot p^{f-i-1} \equiv     \sum_{i \in \{0, \dots f-1\}} \wtop^{\mu}_{i} \cdot p^{f-i-1} \mod (p^f-1),
\end{align}
and
\begin{align}
\label{eq: hcongstring}
    \sum_{i \in \{0, \dots f-1\}} \wbottom'_{i} \cdot p^{f-i-1} \equiv     \sum_{i \in \{0, \dots f-1\}} \wbottom^{\mu}_{i} \cdot p^{f-i-1} \mod (p^f-1),
\end{align}
for each $\mu \in \tilde{M}$.

Similarly, if $\rho$ has crystalline lifts of weights $(k',0)$ and $(k^{\Theta},l^{\Theta})$, then 
    \begin{align}    
    \label{eq: rhoinertprime}
    \rho|_{I_{K}} & \sim \mtwo{\prod_{\tau} \omega_{\tau}^{\wtop'_{\tau}}}{\ast}{0}{\prod_{\tau} \omega_{\tau}^{\wbottom'_{\tau}}} \\
    & \sim \mtwo{\prod_{\tau} \omega_{\tau}^{\wtop^{\Theta}_{\tau}}}{\ast}{0}{\prod_{\tau} \omega_{\tau}^{\wbottom^{\Theta}_{\tau}}}, 
      \label{eq: rhoinerttheta}
    \end{align}
so we obtain
\begin{align}
\label{eq: congthetatopstring}
    \sum_{i \in \{0, \dots f-1\}} \wtop'_{i} \cdot p^{f-i-1} \equiv     \sum_{i \in \{0, \dots f-1\}} \wtop^{\Theta}_{i} \cdot p^{f-i-1} \mod (p^f-1),
\end{align}
and
\begin{align}
\label{eq: hcongthetabottomstring}
    \sum_{i \in \{0, \dots f-1\}} \wbottom'_{i} \cdot p^{f-i-1} \equiv     \sum_{i \in \{0, \dots f-1\}} \wbottom^{\Theta}_{i} \cdot p^{f-i-1} \mod (p^f-1).
\end{align}

\begin{remark}
    \label{symmrandh}
We have $\wtop'_{\tau}-\wtop^{\mu}_{\tau}=-(\wbottom'_{\tau}-\wbottom^{\mu}_{\tau})$ for all $\tau \in \Sigma$ for each $\mu \in \tilde{M}$, and $\wtop'_{\tau}-\wtop^{\Theta}_{\tau}=-(\wbottom'_{\tau}-\wbottom^{\Theta}_{\tau})$ for all $\tau \in \Sigma$.
\end{remark}

\subsection{Reducible matching}
In this section we show that if a reducible representation $\rho: G_K \to \GL_2(\Fpb)$ has a crystalline lift of irregular weights, we can rewrite $\rho|_{I_K}$ so the shape corresponds to a representation of the regular weights. We also show the converse of this statement. This will be sufficient for the semisimple case. 

\subsubsection{Blocks}

We work in \emph{blocks}. This will be a maximal length string of embeddings such that the first part of the block satisfies $k_{\tau} \neq 1$ and the second part satisfies $k_{\tau}=1$. Each block will contain exactly one embedding in $\tilde{M}$, i.e. some embedding $\mu$ such that $k_{\mu} \neq 1$ and $k_{\Frob^{-1} \circ \mu} =1$.

\begin{definition}
\label{defblock}
Consider the irregular weight $(k,0)$.
We say $\tau_{\lowparam}, \dots, \tau_{\highparam}$ for some $\highparam > \lowparam$ is a block if 
\[
\tau_{\lowparam-1} \in J_0 \text{ and } \tau_{\lowparam} \not \in J_0,
\]
and if $\highparam$ is the smallest integer $i>\lowparam$ such that
\[
\tau_i \in J_0 \text{ and } \tau_{i+1} \not \in J_0.
\]
\end{definition}

\begin{example}
Suppose we have a weight $(k,0)$ with $k=(k_0,k_1,1,k_3,1,1)$ where $k_0,k_1,k_3 \geq 2$. Then we have two blocks with embeddings $\tau_0,\tau_1,\tau_2$ and $\tau_3, \tau_4, \tau_5$ corresponding to
\[
    (k_0,k_1,1) \quad \text{ and } \quad (k_3,1,1).
\]
\end{example}

We introduce some more notation:

\begin{definition}
    Let $\mu \in \tilde{M}$. We let $\text{Block}_{\mu}$ be the set of embeddings $\tau \in \Sigma$ that are in the same block as $\mu$.
\end{definition}

\subsubsection{Shape at inertia}

\begin{proposition}[Possible regular shape]
\label{regularshape}
Suppose $\rho: G_K \to \GL_2(\Fpb)$ is a reducible representation with a crystalline lift of weights $(k,0)$ and let $J$ be a subset as in Corollary \ref{corextclass}. Then there are subsets $J'$, $J^{\Theta}$ and $\{J^{\mu}\}_{\mu \in \tilde{M}}$ so that we have 
    \begin{align*}    
    \rho|_{I_{K}} & \sim \mtwo{\prod_{\tau \in J'} \omega_{\tau}^{b'_{\tau,1}} \prod_{\tau \not \in J'} \omega_{\tau}^{b'_{\tau,2}}}{\ast}{0}{\prod_{\tau \in J'} \omega_{\tau}^{b'_{\tau,2}} \prod_{\tau \not \in J'} \omega_{\tau}^{b'_{\tau,1}}}, \\
    & \sim \mtwo{\prod_{\tau \in J^{\Theta}} \omega_{\tau}^{b^{\Theta}_{\tau,1}} \prod_{\tau \not \in J^{\Theta}} \omega_{\tau}^{b^{\Theta}_{\tau,2}}}{\ast}{0}{\prod_{\tau \in J^{\Theta}} \omega_{\tau}^{b^{\Theta}_{\tau,2}} \prod_{\tau \not \in J^{\Theta}} \omega_{\tau}^{b^{\Theta}_{\tau,1}}}, \\
    & \sim \mtwo{\prod_{\tau \in J^{\mu}} \omega_{\tau}^{b^{\mu}_{\tau,1}} \prod_{\tau \not \in J^{\mu}} \omega_{\tau}^{b^{\mu}_{\tau,2}}}{\ast}{0}{\prod_{\tau \in J^{\mu}} \omega_{\tau}^{b^{\mu}_{\tau,2}} \prod_{\tau \not \in J^{\mu}} \omega_{\tau}^{b^{\mu}_{\tau,1}}},
    \end{align*}
for each $\mu \in \tilde{M}$.
The subsets $J'$, $J^{\Theta}$ and $J^{\mu}$ are such that 
\[
\tau \in J^{\mu} \iff \tau \in J' \iff \tau \in J^{\Theta} \iff \tau \in J,
\]
for all $\tau \in \Sigma \setminus J_0$.
\end{proposition}
\begin{proof}
We use the set $J$ to construct sets $J'$, $J^{\Theta}$ and $\{J^{\mu}\}_{\mu \in \tilde{M}}$ so that the above holds.

Similar to before, we write
\[
\wtop_{i}= \begin{cases}
    b_{\tau_i,1} & \tau_i \in J, \\
    b_{\tau_i,2} & \tau_i \not \in J. \\
\end{cases}
\]

For all $\tau$ such that $\tau \in \Sigma \setminus(J_0 \cup \tilde{M})$ we have
\[
{b_{\tau,1}}={b'_{\tau,1}}={b^{\mu}_{\tau,1}}={b^{\Theta}_{\tau,1}}, \quad
{b_{\tau,2}}={b'_{\tau,2}}={b^{\mu}_{\tau,2}}={b^{\Theta}_{\tau,1}},
\]
for each $\mu \in \tilde{M}$. Now for all embeddings $\tau \in \Sigma \setminus J_0$, we let $\tau \in J', \tau \in J^{\mu}, \tau \in J^{\Theta} \iff \tau \in J$, for each $\mu \in \tilde{M}$.

For each set $J', J^{\mu}, J^{\Theta}$ it thus remains to make a choice for embeddings in $J_0$. Once we have made such a choice, we prove congruences per block. Without loss of generality, say $\tau_0$ is the starting embedding of a block and $\tau_E$ is the end of this block. Let $\nu$ be the unique embedding in this block such that $\nu \in \tilde{M}$. Let $i_{\nu}$ be the integer such that $\tau_{i_{\nu}} = \nu$.

We first choose the set $J'$ so that the desired congruence holds. If $\nu \in J'$, for all embeddings $\tau \in J_0 \cap \text{Block}_{\nu}$ we set $\tau \in J'$. 
We obtain the following congruences
\[
\sum_{i \in \{0, \dots, E \}} \wtop_i p^{f-i-1} \equiv \sum_{i \in \{0, \dots, E \}} \wtop'_i p^{f-i-1}  \mod (p^f-1),
\]
since $\wtop_i=\wtop'_i$ for $i=0,  \dots, i_{\nu}-1$ and $\wtop_{i_{\nu}}=k_{\nu}-1$ and $\wtop'_{i_{\nu}}=k_{\nu}-2$ with $\wtop'_{i_{\nu}+1}, \dots, \wtop'_E = p-1, \dots, p-1,p$ (with $\wtop_{i_{\nu}+1}, \dots, \wtop_E = 0$).

If $\nu \not \in J'$, then for all embeddings $\tau \in J_0 \cap \text{Block}_{\nu}$ we set $\tau \not \in J'$. We have
\[
\sum_{i \in \{0, \dots, E \}} \wtop_i p^{f-i-1} \equiv \sum_{i \in \{0, \dots, E \}} \wtop'_i p^{f-i-1}  \mod (p^f-1),
\]
since $\wtop_i=\wtop'_i$ for $i=0,  \dots, E$.

This argument works exactly the same for any other block. We put the congruences together and we indeed find:
\[
\sum_{i \in \{0, \dots, f-1 \}} \wtop_i p^{f-i-1} \equiv \sum_{i \in \{0, \dots, f-1 \}} \wtop'_i p^{f-i-1}  \mod (p^f-1),
\]
and similarly, by symmetry,
\[
\sum_{i \in \{0, \dots, f-1 \}} \wbottom_i p^{f-i-1} \equiv \sum_{i \in \{0, \dots, f-1 \}} \wbottom'_i p^{f-i-1} \mod (p^f-1).
\]
This gives the result for the weights $\{b'_{\tau,1}, b'_{\tau,2}\}_{\tau \in \Sigma}$ corresponding to $(k',l')$.

We fix $\mu \in \tilde{M}$. We define the set $J^{\mu}$, again proceeding per block. If we are dealing with a block not containing $\mu$, then we proceed exactly as above for $J'$. This works because the blocks in $(k^{\mu},l^{\mu})$ not containing $\mu$ are identical to the corresponding blocks $(k',l')$, so for all embeddings in these blocks we have $\{b'_{\tau,1}, b'_{\tau,2}\}=\{b^{\mu}_{\tau,1}, b^{\mu}_{\tau,2}\}$. 

Now consider the block containing $\mu$. If $\mu \in J^{\mu}$, then for all embeddings $\tau \in J_0 \cap \text{Block}_{\mu}$ we set $\tau \not \in J^{\mu}$. We obtain
\[
\sum_{i \in \{0, \dots, E \}} \wtop_i p^{f-i-1} \equiv \sum_{i \in \{0, \dots, E \}} \wtop^{\mu}_i p^{f-i-1}  \mod (p^f-1),
\]
since $\wtop_i=\wtop^{\mu}_i$ for $i=0,  \dots, E$.

Suppose $\mu \not \in J^{\mu}$, then for all embeddings $\tau \in J_0 \cap \text{Block}_{\mu}$ we set $\tau \in J^{\mu}$. We obtain
\[
\sum_{i \in \{0, \dots, E \}} \wtop_i p^{f-i-1} \equiv \sum_{i \in \{0, \dots, E \}} \wtop^{\mu}_i p^{f-i-1}  \mod (p^f-1),
\]
since $\wtop_i=\wtop^{\mu}_i$ for $i=0,  \dots, i_{\mu}-1$ and $\wtop_{i_{\mu}}=0$ and $\wtop^{\mu}_{i_{\mu}}=-1$ with $\wtop^{\mu}_{i_{\mu}+1}, \dots, \wtop^{\mu}_E = p-1, \dots, p-1,p$ (with $\wtop_{i_{\mu}+1}, \dots, \wtop_E = 0$).

Similar to above, this gives the result for $(k^{\mu},l^{\mu})$ for each $\mu \in \tilde{M}$.

Now finally for the set $J^{\Theta}$, we also work by blocks. Suppose $\nu \in \tilde{M}$ is the unique embedding of $\tilde{M}$ in our block, then we find
\[
\{b_{\tau,1}^{\Theta},b_{\tau,2}^{\Theta}\}=\{b_{\tau,1}^{\nu},b_{\tau,2}^{\nu}\}
\]
for all $\tau \in \text{Block}_\nu$. So we proceed as in the case for $J^{\mu}$ with $\nu=\mu$. Doing this for all elements in $\tilde{M}$ then gives the result for $(k^{\Theta},l^{\Theta})$.
\end{proof}

Before we do the converse, we need some auxiliary lemmas. We will prove stronger versions of these later, but these require a longer proof for some special cases. Since we do not need these for the semisimple cases we postpone the stronger statements.

\begin{lemma}
\label{thetasemimatch}
Suppose $\rho: G_K \to \GL_2(\Fpb)$ is a reducible representation with crystalline lifts of weights $(k',l')$ and $(k^{\Theta},l^{\Theta})$. Let $J'$ and $J^{\Theta}$ be sets obtained from Corollary \ref{corextclass}. Let $\nu \in \tilde{M}$ and write $\tau_{i_{\nu}}, \dots, \tau_E$ for the part of $\text{Block}_{\nu}$ starting at $\tau_{i_{\nu}}=\nu$ with $\tau_E$ the end of the block. 
Then either
\[
\tau_{i_{\nu}}, \tau_{i_{\nu}+1}, \dots, \tau_E\in J' \text{and }\tau_{i_\nu+1},\dots, \tau_{i_E} \not \in J^{\Theta},
\]
or
\[
\tau_{i_{\nu}}, \tau_{i_{\nu}+1}, \dots, \tau_E \not \in J' \text{ and }\tau_{i_\nu+1},\dots, \tau_{i_E}  \in J^{\Theta}.
\]
\end{lemma}
\begin{proof}
By our assumptions, the congruence \eqref{eq: congthetatopstring} must be satisfied.
Now suppose $\nu \in J'$, then we have
\[
\wtop'_{i_\nu}-\wtop^{\Theta}_{i_\nu}=
\begin{cases}
    -1, & \text{if } \nu \in J^{\Theta}, \\
    k_{\nu}-1, & \text{if } \nu \not \in J^{\Theta}.
\end{cases}
\]
If $\nu \in J^{\Theta}$, then by Lemma \ref{lem71}, \eqref{eq: bprime} and \eqref{eq: btheta} we must have $\tau_{i_\nu+1},\dots, \tau_{i_E} \in J'$ and $\tau_{i_\nu+1},\dots, \tau_{i_E} \not \in J^{\Theta}$. If $\nu \not \in J^{\Theta}$, noting that by assumption we have $k_{\nu} \neq 2$ (since $\nu \in \tilde{M}$) we conclude again by Lemma \ref{lem71}, \eqref{eq: bprime} and \eqref{eq: btheta} that we must have $\tau_{i_\nu+1},\dots, \tau_{i_E} \in J'$ and $\tau_{i_\nu+1},\dots, \tau_{i_E} \not \in J^{\Theta}$.
Now if $\nu \not \in J'$, then we have 
\[
\wbottom'_{i_\nu}-\wbottom^{\Theta}_{i_\nu}=
\begin{cases}
    -1, & \text{if } \nu \not\in J^{\Theta}, \\
    k_{\nu}-1, & \text{if } \nu \in J^{\Theta}, 
\end{cases}
\]
and since the congruence \eqref{eq: hcongthetabottomstring} must hold it follows as above that we must have $\tau_{i_\nu+1},\dots, \tau_{i_E} \not \in J'$ and $\tau_{i_\nu+1},\dots, \tau_{i_E} \in J^{\Theta}$ by Lemma \ref{lem71} in both cases.
\end{proof}

\begin{proposition}
\label{thetatoirregshape}
Suppose  $\rho: G_K \to \GL_2(\Fpb)$ is a reducible representation with crystalline lifts of weights $(k',0)$ and $(k^{\Theta},l^{\Theta})$. Let $J'$ and $J^{\Theta}$ be as in Corollary \ref{corextclass}. Then there is a subset $J$ such that
  \begin{align*}    
    \rho|_{I_{K}} & \sim \mtwo{\prod_{\tau \in J} \omega_{\tau}^{b_{\tau,1}} \prod_{\tau \not \in J} \omega_{\tau}^{b_{\tau,2}}}{\ast}{0}{\prod_{\tau \in J} \omega_{\tau}^{b_{\tau,2}} \prod_{\tau \not \in J} \omega_{\tau}^{b_{\tau,1}}}.
    \end{align*}
The subset $J$ is such that 
\[
\tau \in J \iff \tau \in J',
\]
for all $\tau \in \Sigma \setminus J_0$.
\end{proposition}

\begin{proof}
We will choose our set $J$ so that the above holds. Note that for the weight $(k,0)$ we have $b_{\tau,1},b_{\tau,2}=0$ for all $\tau \in J_0$. So we only need to decide for $\Sigma \setminus J_0$.
Recall that for all $\tau$ such that $\tau \in \Sigma \setminus(J_0 \cup \tilde{M})$ we have
\[
{b_{\tau,1}}={b'_{\tau,1}}, \quad \text{and} \quad {b_{\tau,2}}={b'_{\tau,2}}.
\]
For such embeddings we set 
\[
\tau \in J \iff \tau \in J'.
\]
Moreover, we set 
\[
\mu \in J \iff \mu\in J',
\]
for all $\mu \in \tilde{M}$.

We prove congruences per block. Say $\tau_0$ is the start of a block and $\tau_E$ is the end of this block. Let $\nu$ be the unique embedding in this block such that $\nu \in \tilde{M}$ and let $i_{\nu}$ be such that $\tau_{i_{\nu}} = \nu$.

Now suppose $\nu \in J'$, by Lemma \ref{thetasemimatch} we have $\tau_{i_\nu+1},\dots, \tau_{i_E} \in J'$. Then we obtain the following congruences
\[
\sum_{i \in \{0, \dots, E \}} \wtop_i p^{f-i-1} \equiv \sum_{i \in \{0, \dots, E \}} \wtop'_i p^{f-i-1}  \mod (p^f-1),
\]
since $\wtop_i=\wtop'_i$ for $i=0,  \dots, i_{\nu}-1$ and $\wtop_{i_{\nu}}-\wtop'_{i_{\nu}}=(k_{\nu}-1) - (k_{\nu}-2)=-1$ with $\wtop'_{i_{\nu}+1}, \dots, \wtop'_E = p-1, \dots, p-1,p$ (with $\wtop_{i_{\nu}+1}, \dots, \wtop_E = 0$). Note this also gives 
\[
\sum_{i \in \{0, \dots, E \}} \wbottom_i p^{f-i-1} \equiv \sum_{i \in \{0, \dots, E \}} \wbottom'_i p^{f-i-1}  \mod (p^f-1),
\]
since $\wbottom_i=\wbottom'_i$ for $i=0,  \dots, E$. 

The case where $\nu \not \in J'$ follows by symmetry.

We put the congruences from each block together and we indeed find:
\[
\sum_{i \in \{0, \dots, f-1 \}} \wtop_i p^{f-i-1} \equiv \sum_{i \in \{0, \dots, f-1 \}} \wtop'_i p^{f-i-1}  \mod (p^f-1),
\]
and by symmetry,
\[
\sum_{i \in \{0, \dots, f-1 \}} \wbottom_i p^{f-i-1} \equiv \sum_{i \in \{0, \dots, f-1 \}} \wbottom'_i p^{f-i-1}  \mod (p^f-1).\qedhere
\]
\end{proof}

\begin{lemma}
\label{nusemimatch}
Suppose $\rho: G_K \to \GL_2(\Fpb)$ is a reducible representation with crystalline lifts of weights $(k',0)$ and $\{(k^{\mu},l^{\mu})\}_{\mu \in \tilde{M}}$.  Let $J'$ and $\{J^{\mu}\}_{\mu \in \tilde{M}}$ be as in Corollary \ref{corextclass}.  Let $\nu \in \tilde{M}$ and write $\tau_{i_{\nu}} = \nu$ and let $\tau_E$ be the last embedding in $\text{Block}_{\nu}$. Then either
\[
\tau_{i_{\nu}}, \tau_{i_{\nu}+1},\dots, \tau_E \in J' \text{ and } \tau_{i_{\nu}+1},\dots, \tau_E \not \in J^{\nu},
\]
or
\[
\tau_{i_{\nu}}, \tau_{i_{\nu}+1},\dots, \tau_E \not \in J' \text{ and } \tau_{i_{\nu}+1},\dots, \tau_E \in J^{\nu}.
\]
\end{lemma}
\begin{proof}
This follows immediately from Lemma \ref{thetasemimatch} and its proof given that the weights $(k^{\nu},l^{\nu})$ and $(k^{\Theta},l^{\Theta})$ are identical for all embeddings in $\text{Block}_{\nu}$.
\end{proof}

\begin{remark}
\label{remarknusemimatch}
Note that the above also holds upon replacing $J'$ by $J^{\mu}$ for any $\mu \in \tilde{M}$ such that $\mu \neq \nu$. This follows from the proof and the fact that the weights $(k^{\mu},l^{\mu})$ and $(k',l')$ are identical at all embeddings outside of $\text{Block}_{\mu}$.
\end{remark}

\begin{proposition}
\label{weightoneshape}
Suppose $\rho: G_K \to \GL_2(\Fpb)$ is a reducible representation with crystalline lifts of weights $(k',l')$ and $\{(k^{\mu},l^{\mu})\}_{\mu \in \tilde{M}}$.  Let $J'$ and $\{J^{\mu}\}_{\mu \in \tilde{M}}$ be as in Corollary \ref{corextclass}, then there is a subset $J$ so that we have 
    \begin{align*}    
    \rho|_{I_{K}} & \sim \mtwo{\prod_{\tau \in J} \omega_{\tau}^{b_{\tau,1}} \prod_{\tau \not \in J} \omega_{\tau}^{b_{\tau,2}}}{\ast}{0}{\prod_{\tau \in J} \omega_{\tau}^{b_{\tau,2}} \prod_{\tau \not \in J} \omega_{\tau}^{b_{\tau,1}}}.
    \end{align*}
The subset $J$ is such that 
\[
\tau \in J \iff \tau \in J',
\]
for all $\tau \in \Sigma \setminus J_0$.
\end{proposition}
\begin{proof}
The proof is the same as the proof of Proposition \ref{thetatoirregshape}, with Lemma \ref{nusemimatch} used instead of Lemma \ref{thetasemimatch}. 
\end{proof}

\subsection{Proof in the semisimple case}
The above results allow us to prove our main theorems in the semisimple case.

\begin{proposition}
\label{semisimpleproof}
Let $p$ be an odd prime. Let $(k,0) \in \Z^{\Sigma} \times \Z^{\Sigma}$ be an irregular weight satisfying $1 \leq k_{\tau} \leq p$ for all $\tau \in \Sigma$ and such that $k_{\tau} \neq 2$ if $k_{\Frob^{-1} \circ \tau}=1$.  Suppose further $k \neq 1$. Suppose $\rho: G_K \to \GL_2(\Fpb)$ is a semisimple reducible representation. Then we have the following equivalence
\begin{align*}
&\text{
$\rho$ has a crystalline lift of weight $({k},{0})$} \\
&\iff \text{$\rho$ has crystalline lifts of weights $(k',l')$ and $(k^{\Theta},l^{\Theta})$} \\
&\iff \text{$\rho$ has crystalline lifts of weights $(k',l')$ and $\{(k^{\mu},l^{\mu})\}_{\mu \in \tilde{M}}$.}
\end{align*}
\end{proposition}

\begin{proof}
First assume $\rho$ has a lift of weight $(k,0)$. By Lemma \ref{regularshape} we know that for each regular weight $(k',l')$, $(k^{\Theta},l^{\Theta})$ and $\{(k^{\mu},l^{\mu})\}_{\mu \in \tilde{M}}$ there exists sets $J', J^{\Theta}$ and $\{J^{\mu}\}_{\mu \in \tilde{M}}$ such that $\rho|_{{I_K}}$ can be rewritten to have the shape corresponding to these weights. Since $\rho$ is semisimple, it follows from Lemma \ref{cryslift} that $\rho$ has a crystalline lift of the required regular weights. 

If $\rho$ has lifts of weights $(k',l')$ and $(k^{\Theta},l^{\Theta})$, then by Proposition \ref{thetatoirregshape} we know that there exist a set $J$ such that $\rho|_{{I_K}}$ can be written into the shape corresponding to $(k,0)$. Since $\rho$ is semisimple, it follows from Lemma \ref{cryslift} that $\rho$ has a crystalline lift of weight $(k,0)$. If $\rho$ instead has crystalline lifts of weights $(k',l')$ and $\{(k^{\mu},l^{\mu})\}_{\mu \in \tilde{M}}$, then it follows similarly from Proposition \ref{weightoneshape}. 
\end{proof}

\section{Non-semisimple case}
\label{nonsemisimplesection}
The non-semisimple case is more complex. We need to study various subspaces of $H^1(G_K,\Fpb(\overline{\psi}))$ as in Definition \ref{extspace}, and prove these are the same under some conditions. This makes extensive use of Kisin modules and our results in Section \ref{reductionsofcrysreps}. In the final subsection we use these results to finish the reducible case.

\subsection{Descriptions of the relevant Kisin modules}
If a reducible representation $\rho: G_K \to \GL_2(\Fpb)$ has a crystalline lift of certain Hodge--Tate weights, under some conditions on the Hodge--Tate weights, we can use Theorem \ref{thm79} to obtain an extension of two rank one Kisin modules. Below we describe the extensions for the sets of Hodge--Tate weights we are interested in. 

\begin{remark}
\label{remarknotexceptional}
We want to use Theorem \ref{thm79} to describe the extensions using the parameters $x_i$. These are constants in $\F_L$, outside of some exceptional cases. We assume we are not in these cases so that we will have $x_i=0$ if $i \in J_0$ or $i \not \in J$ and $x_i \in \F_L$ otherwise. It follows from Lemma \ref{notexceptional} in the Appendix that we can make this assumption for the irregular weight. For the regular weights we justify why we can have this assumption later in the proofs when we need them.
\end{remark}

\subsubsection{Description of extension corresponding to the irregular weight}
\label{irregularweightdescription}
\textbf{The weight $(k,0)$.}
If a reducible representation $\rho: G_K \to \GL_2(\Fpb)$ has a crystalline lift of irregular weight $(k,0)$, let $J$ be a subset given by Corollary \ref{corextclass}. We write $\psi_1,\psi_2$ for the corresponding characters as in \eqref{eq: charpsi} and $L_{\psi_1,\psi_2}$ for the space of extension classes. 

For each element in the space $L_{\psi_1,\psi_2}$ we obtain an extension $\ofM$ of Kisin modules as in Theorem \ref{thm79}:
\begin{align}
\label{eq: kisinirreg}
0 \to \ofM(\wbottom_0,\dots,\wbottom_{f-1};b) \to \ofM \to \ofM(\wtop_0,\dots,\wtop_{f-1};a) \to 0, 
\end{align}
where $s_i$ and $t_i$ are as in Section \ref{weightnotation} 
and where $a, b \in \F_L$. We obtain bases $\{e_i,f_i\}$ for the $\ofM_i$ such that
\begin{align*}    
\varphi(e_{i-1})=(b)_i u^{\wbottom_i} e_i, \quad \varphi(f_{i-1})=(a)_i u^{\wtop_i} f_i + x_i e_i,
\end{align*}
where $x_i \in \F_L$ and $x_i=0$ if $i \in J_0$ or $i \not \in J$.

\subsubsection{Description of extension corresponding to the regular weights}
\label{regularweightdescription}
\textbf{The weight $(k',l')$.} If $\rho: G_K \to \GL_2(\Fpb)$ has a crystalline lift of weight $(k',l')$, this gives a subset $J'$ as in Corollary \ref{corextclass}. We write $\psi'_1,\psi'_2$ for the corresponding characters as in \eqref{eq: charpsi} and $L_{\psi'_1,\psi'_2}$ for the space of extension classes. 

For any element in $L_{\psi'_1,\psi'_2}$ we find by Theorem \ref{thm79}:
\begin{align}
\label{eq: kisinprime}
0 \to \ofM(\wbottom'_0,\dots,\wbottom'_{f-1};b') \to \ofM' \to \ofM(\wtop'_0,\dots,\wtop'_{f-1};a') \to 0, 
\end{align}
with $\wtop'_i$ and $\wbottom'_i$ as in \eqref{eq: rexpressions} and \eqref{eq: sexpressions} and where $a',b' \in \F_L$. We obtain bases $\{e'_i,f'_i\}$ for the $\fM'_i$ such that
\begin{align*}
\varphi(e'_{i-1})=(b')_i u^{\wbottom'_i} e'_i, \quad \varphi(f'_{i-1})=(a')_i u^{\wtop'_i} f'_i + x'_i e'_i,
\end{align*}
where $x'_i \in \F_L$ and $x'_i=0$ if $\tau_i \not \in J'$.

\noindent \textbf{The weight $(k^{\mu},l^{\mu})$.} Fix $\mu \in \tilde{M}$.  If $\rho: G_K \to \GL_2(\Fpb)$ has a crystalline lift of weight $(k^{\mu},l^{\mu})$ this gives a subset $J^{\mu}$ as in Corollary \ref{corextclass}. We write $\psi^{\mu}_1,\psi^{\mu}_2$ for the corresponding characters as in \eqref{eq: charpsi} and $L_{\psi^{\mu}_1,\psi^{\mu}_2}$ for the space of extension classes. Before moving to Kisin modules, we twist any extension corresponding to an element in $L_{\psi^{\mu}_1,\psi^{\mu}_2}$ by $\omega_{\mu}$. This removes the one negative Hodge--Tate weight at $\mu$.

Now we can apply Theorem \ref{thm79} and we obtain an extension $\ofM^{\mu}$ 
\begin{align}
\label{eq: kisinmu}
0 \to \ofM(\tilde{t}^{\mu}_0,\dots,\tilde{t}^{\mu}_{f-1};b^{\mu}) \to \ofM^{\mu} \to \ofM(\tilde{s}^{\mu}_0,\dots,\tilde{s}^{\mu}_{f-1};a^{\mu}) \to 0,
\end{align}
where
\begin{align*}
    \tilde{s}^{\mu}_i = \begin{cases}
        \wtop^{\mu}_i, & \tau_i \neq \mu, \\
        \wtop^{\mu}_i+1, & \tau_i = \mu,
    \end{cases} \quad
    \tilde{t}^{\mu}_i = \begin{cases}
        \wbottom^{\mu}_i, & \tau_i \neq \mu, \\
        \wbottom^{\mu}_i+1, & \tau_i = \mu,
    \end{cases}
\end{align*}
with $\wtop^{\mu}_i$ and $\wbottom^{\mu}_i$ as in \eqref{eq: rexpressions} and \eqref{eq: sexpressions} and $a^{\mu}, b^{\mu} \in \F_L$. We obtain bases $\{e^{\mu}_i,f^{\mu}_i\}$ for $\ofM^{\mu}_i$ such that
\begin{align*}
    \varphi(e^{\mu}_{i-1}) = (b^{\mu})_i u^{\tilde{t}^{\mu}_i} e^{\mu}_i, \quad    \varphi(f^{\mu}_{i-1}) = (a^{\mu})_i u^{\tilde{s}^{\mu}_i} f^{\mu}_i + x^{\mu}_i e^{\mu}_i,
\end{align*}
where $x^{\mu}_i \in \F_L$ and $x^{\mu}_i=0$ if $\tau_i \not \in J^{\mu}$.

\noindent \textbf{The weight $(k^{\Theta},l^{\Theta})$.} If $\rho: G_K \to \GL_2(\Fpb)$ has a crystalline lift of weight $(k^{\Theta},l^{\Theta})$ this gives a subset $J^{\Theta}$ as in Corollary \ref{corextclass}. We write $\psi^{\Theta}_1,\psi^{\Theta}_2$ for the corresponding characters as in \eqref{eq: charpsi} and $L_{\psi^{\Theta}_1,\psi^{\Theta}_2}$ for the space of extension classes. Before moving to Kisin modules, we twist any extension corresponding to an element in $L_{\psi^{\Theta}_1,\psi^{\Theta}_2}$ by $\prod_{\mu \in \tilde{M}} \omega_{\mu}$. This removes the negative Hodge--Tate weights at these embeddings.

Now we can apply Theorem \ref{thm79} and we obtain an extension $\ofM^{\Theta}$ 
\begin{align}
\label{eq: kisintheta}
0 \to \ofM(\tilde{t}^{\Theta}_0,\dots,\tilde{t}^{\Theta}_{f-1};b^{\Theta}) \to \ofM^{\Theta} \to \ofM(\tilde{s}^{\Theta}_0,\dots,\tilde{s}^{\Theta}_{f-1};a^{\Theta}) \to 0,
\end{align}
where
\begin{align*}
    \tilde{s}^{\Theta}_i = \begin{cases}
        \wtop^{\Theta}_i, & \tau_i \not \in \tilde{M}, \\
        \wtop^{\Theta}_i+1, & \tau_i \in \tilde{M},
    \end{cases} \quad
    \tilde{t}^{\Theta}_i = \begin{cases}
        \wbottom^{\Theta}_i, & \tau_i \not \in \tilde{M}, \\
        \wbottom^{\Theta}_i+1, & \tau_i \in \tilde{M},
    \end{cases}
\end{align*}
with $\wtop^{\Theta}_i$ and $\wbottom^{\Theta}_i$ as in \eqref{eq: rexpressions} and \eqref{eq: sexpressions} and where $a^{\Theta}, b^{\Theta} \in \F_L$. We obtain bases $\{e^{\Theta}_i,f^{\Theta}_i\}$ for $\ofM^{\Theta}_i$ such that
\begin{align*}
    \varphi(e^{\Theta}_{i-1}) = (b^{\Theta})_i u^{\tilde{t}^{\Theta}_i} e^{\Theta}_i, \quad    \varphi(f^{\Theta}_{i-1}) = (a^{\Theta})_i u^{\tilde{s}^{\Theta}_i} f^{\Theta}_i + x^{\Theta}_i e^{\Theta}_i,
\end{align*}
where $x^{\Theta}_i \in \F_L$ and $x^{\Theta}_i=0$ if $\tau_i \not \in J^{\Theta}$.

\subsection{Identifying subspaces of $H^1(G_K,\Fpb(\overline{\psi}))$}
In this section we use the Kisin modules introduced above as well as the machinery from earlier sections to identify certain subspaces of $H^1(G_K,\Fpb(\overline{\psi}))$.
We need some auxiliary results. Some of these have been moved to the appendix, as to not distract too much from the main arguments.

\subsubsection{Auxiliary results}
We start with an auxiliary result which is a simpler variant of \cite[Proposition 5.2.3]{GLS15}. 

\begin{proposition}
\label{prop523irreg}
Let $s_0\dots, s_{f-1}, s'_0, \dots, s'_{f-1},t_0,\dots, t_{f-1},t'_0 \dots, t'_{f-1}$ be non-negative integers. Let $a,b \in \F_L$.
Suppose we are given Kisin modules
\[
\ofN := \ofM(s_0,\dots,s_{f-1};a) \text{ and } \ofP:=\ofM(t_0,\dots,t_{f-1};b),
\]
as well as
\[
\ofN' := \ofM(s'_0,\dots,s'_{f-1};a) \text{ and } \ofP':=\ofM(t'_0,\dots,t'_{f-1};b),
\]
such that there exist non-zero maps $\ofP \to \ofP'$ and $\ofN \to \ofN'$.
Consider $\ofM \in \Ext^1(\ofN,\ofP)$ such that we have bases $e_i,f_i$ of the $\ofM_i$ so that $\varphi$ has the form:
\[
\varphi(e_{i-1})=(b)_i u^{t_i} e_i, \quad \varphi(f_{i-1})=(a)_i u^{s_i} f_i + x_i e_i,
\]
with $x_i \in \F_L$. Suppose that if $x_i \neq 0$ that then $\alpha_{i-1}(\ofN,\ofN')=0$.
Then there exists $\ofM' \in \Ext^1(\ofN',\ofP')$ such that 
\[
T_{\fS}(\ofM) \cong T_{\fS}(\ofM'),
\]
with basis $e'_i,f'_i$ such that 
\[
\varphi(e'_{i-1})=(b)_i u^{t'_i} e'_i, \quad \varphi(f'_{i-1})=(a)_i u^{s'_i} f'_i +  x'_i e'_i,
\]
with $x'_i$ satisfying $u^{\alpha_i(\ofP,\ofP')} x_i= x'_i$.
\end{proposition}
\begin{proof}
By our assumptions and Lemma \ref{lemma512} we have $\alpha_i(\ofP,\ofP') \in \Z_{\geq 0}$ and $\alpha_i(\ofN,\ofN') \in \Z_{\geq 0}$.

We define $\ofM' \in \Ext^1(\ofN',\ofP')$ by the formulae:
\[
\varphi(e'_{i-1})=(b)_i u^{t'_i} e'_i, \quad \varphi(f'_{i-1})=(a)_i u^{s'_i} f'_i +  x'_i e'_i,
\]
with $x'_i$ such that $u^{\alpha_i(\ofP,\ofP')} x_i= x'_i$.
We next pick a basis for $\ofM$ as in the statement. We define a morphism $g: \ofM \to \ofM'$ by:
\begin{align*}
    e_i \mapsto u^{\alpha_i(\ofP,\ofP')} e'_i, \quad   f_i \mapsto u^{\alpha_i(\ofN,\ofN')} f'_i.
\end{align*}
We verify this map is compatible with $\varphi$, i.e. that we have $g \circ \varphi= \varphi \circ g$. We have: 
\begin{align*}
  g(\varphi(e_{i-1})) & = & g( (b)_i \cdot u^{\wbottom_i} e_i) &= &  (b)_i \cdot u^{\alpha_i(\ofP,\ofP')} u^{\wbottom_i} e_i, \\
    \varphi(g(e_{i-1})) & = & \varphi(u^{\alpha_{i-1}(\ofP,\ofP')} e'_{i-1}) & = &  u^{p(\alpha_{i-1}(\ofP,\ofP'))} \cdot (b)_i \cdot u^{\wbottom'_i} e'_i, 
\end{align*}
and
\begin{align*}
  g(\varphi(f_{i-1})) & = & g((a)_i \cdot u^{\wtop_i} f_i + x_i e_i) &= & (a)_i \cdot u^{\wtop_i} u^{\alpha_i(\ofN,\ofN')} f'_i + u^{\alpha_i(\ofP,\ofP')} x_i e'_i,  \\
      \varphi(g(f_{i-1})) &= & \varphi(u^{\alpha_{i-1}(\ofN,\ofN')} f'_{i-1}) &= &  u^{p(\alpha_{i-1}(\ofN,\ofN'))} \cdot ( (a)_i \cdot u^{\wtop'_i} f'_i + x'_i e'_i)).
\end{align*}
and the compatibility follows since the equalities
\begin{align*}  
\alpha_i(\ofP,\ofP') +\wbottom_i &= p(\alpha_{i-1}(\ofP,\ofP')) + \wbottom'_i, \\
\alpha_i(\ofN,\ofN') + \wtop_i &= p(\alpha_{i-1}(\ofN,\ofN')) + \wtop'_i,
\end{align*}
follow from the identity \eqref{eq: alphaidentity}. The rest follows from the definition of $x'_i$ and our assumptions on $\alpha_{i-1}(\ofN,\ofN')$.

Now it follows from Lemma \ref{lem58} and Theorem \ref{thm52} that we obtain an isomorphism $\Tkis(\ofM') \xrightarrow[]{\sim} \Tkis(\ofM)$ since the morphism $g: \ofM \to \ofM'$ induces an isomorphism $\ofM[1/u] \xrightarrow[]{\sim} \ofM'[1/u]$ of étale $\varphi$-modules.
\end{proof}

Next we write down a more direct adaptation of \cite[Proposition 5.2.3]{GLS15}. We use this result to relate different spaces of extensions of Kisin modules. The proof follows directly from that of \cite[Proposition 5.2.3]{GLS15}. The key difference with Proposition \ref{prop523irreg} is that the map between the rank one Kisin modules are in the opposite direction, as we make precise below.

\begin{proposition}
\label{prop523reg}
Let $s_0\dots, s_{f-1}, s'_0, \dots, s'_{f-1},t_0,\dots, t_{f-1},t'_0 \dots, t'_{f-1}$ be non-negative integers. Let $a,b \in \F_L$.
Suppose we are given Kisin modules
\[
\ofN := \ofM(s_0,\dots,s_{f-1};a) \text{ and } \ofP:=\ofM(t_0,\dots,t_{f-1};b),
\]
as well as
\[
\ofN' := \ofM(s'_0,\dots,s'_{f-1};a) \text{ and } \ofP':=\ofM(t'_0,\dots,t'_{f-1};b),
\]
such that there exist non-zero maps $\ofP \to \ofP'$ and $\ofN' \to \ofN$. 

Consider $\ofM \in \Ext^1(\ofN,\ofP)$ such that we have bases $e_i,f_i$ of the $\ofM_i$ so that $\varphi$ has the form:
\[
\varphi(e_{i-1})=(b)_i u^{t_i} e_i, \quad \varphi(f_{i-1})=(a)_i u^{s_i} f_i + y_i e_i,
\]
with $y_i \in \F_L$ or $y_i=u x_i$ for some $x_i \in \F_L$.
Then there exists $\ofM' \in \Ext^1(\ofN',\ofP')$ with basis $e'_i,f'_i$ such that 
\[
\varphi(e'_{i-1})=(b)_i u^{t'_i} e'_i, \quad \varphi(f'_{i-1})=(a)_i u^{s'_i} f'_i + y'_i e'_i,
\]
with $y'_i=u^{\alpha_i(\ofP,\ofP')+p(\alpha_{i-1}(\ofN',\ofN))} y_i$, such that $T_{\fS}(\ofM) \cong T_{\fS}(\ofM')$.
\end{proposition}
\begin{proof}
This follows directly from the proof of \cite[Proposition 5.2.3]{GLS15}.    
\end{proof}

\subsubsection{Equalities of irregular and regular subspaces}

If $\rho: G_K \to \GL_2(\Fpb)$ has a crystalline lift of irregular weight $(k,0)$, then by Corollary \ref{corextclass} we have $c_{\rho} \in L_{\psi_1,\psi_2}$ for some subset $J$. We study this space and show it equals certain subspaces of the regular spaces $L_{\psi'_1,\psi'_2}$, $L_{\psi^{\mu}_1, \psi^{\mu}_2}$ and $L_{\psi^{\Theta}_1, \psi^{\Theta}_2}$ under some assumptions on the characters.

\begin{notation}
In the following we use $\ofN,\ofP$ for the rank one modules corresponding to $\ofM$ as in \eqref{eq: kisinirreg}, so that $\ofM \in \Ext^1(\ofN,\ofP)$. We do the same for the regular extensions, i.e. $\ofM' \in \Ext^1(\ofN',\ofP'), \ofM^\mu \in \Ext^1(\ofN^{\mu},\ofP^{\mu})$ and $\ofM^{\Theta} \in \Ext^1(\ofN^{\Theta},\ofP^{\Theta})$ corresponding to \eqref{eq: kisinprime}, \eqref{eq: kisinmu} and \eqref{eq: kisintheta}, respectively.
\end{notation}

\begin{proposition}
\label{isoregprime}
Consider $L_{\psi_1,\psi_2}$ as in Section \ref{irregularweightdescription} for a fixed $J$. Let $J'$ be as in the proof of Proposition \ref{regularshape} and consider $L_{\psi'_1,\psi'_2}$ as in Section \ref{regularweightdescription}. Suppose $\overline{\psi_1 \psi_2^{-1}}=\overline{\psi'_1} \overline{\psi_2^{'-1}}$. 
Consider the subspace $L'_\text{sub}$ of elements in $L_{\psi'_1,\psi'_2}$ such that the corresponding Kisin modules as in \eqref{eq: kisinprime} additionally satisfy $x'_i=0$ if $i \in J_0$.  
Then $L'_\text{sub} = L_{\psi_1,\psi_2}$.
\end{proposition}
\begin{proof}
By Lemma \ref{primenotexceptional} we are not in the exceptional case of Theorem \ref{thm79}. Now for an element in $L'_{\text{sub}}$ consider the corresponding extension $\ofM'$ as in \eqref{eq: kisinprime} for some $a', b' \in \F_L$ and with $x'_i=0$ if $i \in J_0$ or if $i \not \in J_0$ and $i \not \in J'$.

For any $\ofM$ corresponding to an element in $L_{\psi_1,\psi_2}$ as in \eqref{eq: kisinirreg} we can assume $a=a'$ and $b=b'$ (possibly after twisting by $\ofM(0,\dots,0;c)$ where $c$ is such that $\Tkis(\ofM(0,\dots,0;c)) = (\overline{\psi'_1}\overline{\psi^{}_1}^{-1})|_{G_{\infty}}=(\overline{\psi'_2}\overline{\psi_2}^{-1})|_{G_{\infty}}$).

Let $\ofN, \ofP$ be the two rank one modules as in \eqref{eq: kisinirreg} with $a=a'$ and $b=b'$ (for the set $J$ as above). We check the conditions of Proposition \ref{prop523irreg}. The conditions on the basis for the $\ofM'_i$ follow from \eqref{eq: kisinprime}. By Lemma \ref{alpharegtoirreg} we find that there exist non-zero maps $\ofN' \to \ofN$ and $\ofP' \to \ofP$.  By our assumptions we have $x'_i=0$ if $i \not \in J'$ or $i \in J_0$. By Lemma \ref{alpharegtoirreg} we then see that $\alpha_{i-1}(\ofN',\ofN)=0$ if $x'_i \neq 0$. So we have $\Tkis(\ofM') \cong \Tkis(\ofM)$ for some $\ofM \in \Ext^1(\ofN,\ofP)$ and we obtain $x_i$ such that
\[
u^{\alpha_i(\ofP',\ofP)} x'_i=  x_i.
\] 
When we have $x'_i=0$, we find $x_i=0$. This happens if $i \not \in J'$ or $i \in J_0$. By Lemma \ref{alpharegtoirreg}, it is easy to see that in the remaining cases with $x'_i \neq 0$ we have $\alpha_i(\ofP',\ofP)=0$. So we have $x'_i=x_i$ for all $i$, and by our assumptions on the set $J'$, we see that $\ofM$ is as in Section \ref{irregularweightdescription}. By Remark \ref{diffxi} and Lemma \ref{lemma93} we find that $\ofM$ is the Kisin module corresponding to an element in $L_{\psi_1,\psi_2}$.

By Lemma \ref{lem542}, we find that we have an isomorphism of $G_K$-representations $\hat{T}(\hat{\ofM'}) \cong \hat{T}(\hofM)$. By Remark \ref{diffxi}, we know that different choices of the $x_i$ give different Galois representations. By Lemma \ref{lemma93} we find that $L_{\psi_1,\psi_2}$ has $\F_L^{\{ \tau \in J \, | \, \tau \not \in J_0\}}$ elements and we find $L'_\text{sub}$ does too. So the spaces are the same.
\end{proof}

\begin{proposition}
\label{isoregmu}
Consider $L_{\psi_1,\psi_2}$ as in Section \ref{irregularweightdescription} for a fixed set $J$. Let $\mu \in \tilde{M}$. Consider $L_{\psi^{\mu}_1,\psi^{\mu}_2}$ as in Section \ref{regularweightdescription} for the set $J^{\mu}$ defined in the proof of Proposition \ref{regularshape}. Suppose $\overline{\psi_1 \psi_2^{-1}}=\overline{\psi^{\mu}_1} (\overline{\psi^{\mu}_2)^{-1}}$.
Consider the subspace $L^{\mu}_\text{sub}$ of extensions  in $L_{\psi^{\mu}_1,\psi^{\mu}_2}$  such that the corresponding Kisin modules as in \eqref{eq: kisinmu} additionally satisfy $x^{\mu}_i=0$ if $i \in J_0$. 

Then $L^{\mu}_\text{sub} = L_{\psi_1,\psi_2}$.
\end{proposition}
\begin{proof}
By Lemma \ref{munotexceptional} we are not in the exceptional case of Theorem \ref{thm79}.
Now for an element in $L^{\mu}_{\text{sub}}$ consider the corresponding extension $\ofM^{\mu}$ as described in \eqref{eq: kisinmu} for $J^{\mu}$ as above.

For any $\ofM$ corresponding to an element in $L_{\psi_1,\psi_2}$ as in \eqref{eq: kisinirreg} we can assume $a=a^{\mu}$ and $b=b^{\mu}$ (possibly after twisting by $\ofM(0,\dots,0;c)$ where $c$ is such that $\Tkis(\ofM(0,\dots,0;c)) = (\overline{\psi^{\mu}_1}\overline{\psi^{}_1}^{-1})|_{G_{\infty}}=(\overline{\psi^{\mu}_2}\overline{\psi_2}^{-1})|_{G_{\infty}}$).

Next let
\[
\ofN_{twist} := \ofN \otimes \ofM(d^{\mu}_0,\dots,d^{\mu}_{f-1};1), \quad \ofP_{twist}:= \ofP  \otimes \ofM(d^{\mu}_0,\dots,d^{\mu}_{f-1};1).
\]
where $d^{\mu}_i=1$ if $\tau_i=\mu$ and $0$ otherwise. Note these twists add $1$ to $s_i,t_i$ for $i$ such that $\tau_i=\mu$. In Section \ref{regularweightdescription} we added 1 to $s^{\mu}_i,t^{\mu}_i$ for $i$ such that $\tau_i=\mu$ in order to apply Theorem \ref{thm79} (by twisting by a fundamental character). 
Now by Lemma \ref{lemma512} and Lemma \ref{alpharegtoirreg} we find that indeed there exist non-zero maps $\ofP^{\mu} \to \ofP_{twist}$ and $\ofN^{\mu} \to \ofN_{twist}$.

We check the remaining conditions of Proposition \ref{prop523irreg}. We have $x_i^{\mu}=0$ if $i \in J_0$ or if $i \not \in J^{\mu}$. By Lemma \ref{alpharegtoirreg} we see that $\alpha_{i-1}(\ofN^{\mu},\ofN_{twist})=0$ if $x^{\mu}_i \neq 0$.

The conditions on the basis for the $\ofM^{\mu}_i$ in Proposition \ref{prop523irreg} follow from \eqref{eq: kisinmu}, so we find some $\ofM_{\text{twist}} \in \Ext^1(\ofN_{twist},\ofP_{twist})$ such that 
\[
\Tkis(\ofM^\mu) \cong \Tkis(\ofM_{twist}),
\]
with $x_{twist,i}$ such that 
\[
u^{\alpha_i(\ofP^{\mu},\ofP_{twist})} x^{\mu}_i= x_{twist,i}.
\] 
When we have $x_i^{\mu}=0$, i.e. if $i \in J_0$ or $i \not \in J^{\mu}$, then we obtain $x_{twist,i}=0$.

Now suppose that $i \not \in J_0$ and $i \in J^{\mu}$. Then by Lemma \ref{alpharegtoirreg} we find 
\[
\alpha_i(\ofP^{\mu},\ofP_{twist})= \begin{cases}
1, & i=i_{\mu},\\
0, & \text{otherwise}.    
\end{cases}
\]
So
\[
x_{twist, i} =
\begin{cases}
u x^{\mu}_i, & i=i_{\mu}, \\
x^{\mu}_i, & \text{otherwise.}
\end{cases}
\]
It is easy to see we have $\Tkis(\ofM_{twist}) \cong \Tkis(\ofM \otimes \ofM(d^{\mu}_0,\dots,d^{\mu}_{f-1};1)$ for some $\ofM \in \Ext^1(\ofN,\ofP)$ as in \eqref{eq: kisinirreg} with $a=a^{\mu}, b=b^{\mu}, x_i=x_i^{\mu}$. By our assumptions on the set $J^{\mu}$, we see that $\ofM$ is as in Section \ref{irregularweightdescription}. By Remark \ref{diffxi} and Lemma \ref{lemma93} we find that $\ofM$ is the Kisin module corresponding to an element in $L_{\psi_1,\psi_2}$.

We have $\Tkis(\ofM^\mu) \cong \Tkis(\ofM \otimes \ofM(d^{\mu}_0,\dots,d^{\mu}_{f-1};1))$, and by Lemma \ref{lem542}, we obtain an isomorphism of $G_K$-representations.
By Remark \ref{diffxi}, we know that different choices of the $x_i$ give different Galois representations. By Lemma \ref{lemma93} we find that both $L_{\psi_1,\psi_2}$ and $L^{\mu}_{\text{sub}}$ have $|\F_L|^{| J \setminus (J \cap J_0)|}$ elements and the spaces are the same (the twists are compatible by Section \ref{twistsubsection}).
\end{proof}

\begin{proposition}
\label{isoregtheta}
Consider $L_{\psi_1,\psi_2}$ as in Section \ref{irregularweightdescription} for a fixed set $J$. Consider $L_{\psi^{\Theta}_1,\psi^{\Theta}_2}$ as in Section \ref{regularweightdescription} for the set $J^{\Theta}$ defined as in the proof of Proposition \ref{regularshape}. Suppose $\overline{\psi_1 \psi_2^{-1}}=\overline{\psi^{\Theta}_1} \overline{\psi^{\Theta -1}_2}$.
Consider the subspace $L^{\Theta}_\text{sub}$ of extensions  in $L_{\psi^{\Theta}_1,\psi^{\Theta}_2}$ such that the corresponding Kisin modules as in \eqref{eq: kisintheta} additionally satisfy $x^{\Theta}_i=0$ if $i \in J_0$. 
Then $L^{\Theta}_\text{sub} = L_{\psi_1,\psi_2}$.
\end{proposition}

\begin{proof}
By Lemma \ref{thetanotexceptional} we are not in the exceptional case of Theorem \ref{thm79}.
Now for each element in $L^{\Theta}_{\text{sub}}$ we consider the corresponding extension $\ofM^{\Theta}$ with constants $a^{\Theta}$, $b^{\Theta} \in \F_L$ as described in \eqref{eq: kisintheta} for $J^{\Theta}$ and $x_i^{\Theta}$ as above. 

For any $\ofM$ corresponding to an element in $L_{\psi_1,\psi_2}$ as in \eqref{eq: kisinirreg} we can assume $a=a^{\Theta}$ and $b=b^{\Theta}$ (possibly after twisting by $\ofM(0,\dots,0;c)$ where $c$ is such that $\Tkis(\ofM(0,\dots,0;c)) = (\overline{\psi^{\Theta}_1}\overline{\psi^{}_1}^{-1})|_{G_{\infty}}=(\overline{\psi^{\Theta}_2}\overline{\psi_2}^{-1})|_{G_{\infty}}$).

Next write 
\begin{align*}   
\ofN_{twist} := \ofN \otimes \ofM(d_0,\dots,d_{f-1};1), \quad 
\ofP_{twist}:= \ofP  \otimes \ofM(d_0,\dots,d_{f-1};1),
\end{align*}
with $d_i=1$ if $\tau_i \in \tilde{M}$ and $0$ otherwise. Note these twists add $1$ to $s_i,t_i$ for $i$ such that $\tau_i \in \tilde{M}$.  In Section \ref{regularweightdescription} we added 1 to $s^{\Theta}_i,t^{\Theta}_i$ for $i$ such that $\tau_i \in \tilde{M}$ in order to apply Theorem \ref{thm79} (by twisting by a product of fundamental characters). 
By Lemma \ref{lemma512} and Lemma \ref{alpharegtoirreg} we find that there exist non-zero maps $\ofP^{\Theta} \to \ofP_{twist}$ and $\ofN^{\Theta} \to \ofN_{twist}$. We check the remaining conditions in Proposition \ref{prop523irreg}. We have $x_i^{\Theta}=0$ if $i \in J_0$ or if $i \not \in J^{\Theta}$. By Lemma \ref{alpharegtoirreg} we see that $\alpha_{i-1}(\ofN^{\Theta},\ofN_{twist})=0$ if $x^{\Theta}_i \neq 0$.

The conditions on the basis for the $\ofM^{\Theta}_i$ in Proposition \ref{prop523irreg} follow from \eqref{eq: kisintheta}, so by that result we have $\Tkis(\ofM^\Theta) \cong \Tkis(\ofM_{twist})$ for some $\ofM_{twist} \in \Ext^1(\ofN_{twist},\ofP_{twist})$ as above. We obtain parameters $x_{twist,i}$ such that
\[
u^{\alpha_i(\ofP^{\Theta},\ofP_{twist})} x^{\Theta}_i= x_{twist,i}.
\] 
When we have $x_i^{\Theta}=0$, i.e. if $i \in J_0$ or if $i \not \in J^{\Theta}$, then the above gives us $x_{twist,i}=0$. 

Suppose $x_i^{\Theta} \neq 0$, so that $i \in J^{\Theta}$ and $i \not \in J_0$. We use Lemma \ref{alpharegtoirreg} to find that for such $i$ we have
\[
\alpha_i(\ofP^{\Theta},\ofP_{twist})=
\begin{cases}
 1, & i \in \tilde{M}, \\
 0, & \text{otherwise.}  
\end{cases} 
\]
So this gives
\[
x_{twist, i} =
\begin{cases}
u x^{\Theta}_{i}, & i \in \tilde{M}, \\
x^{\Theta}_{i}, & i \not  \in \tilde{M}.
\end{cases}
\]
It is not hard to see that $\Tkis(\ofM_{twist}) \cong \Tkis(\ofM \otimes \ofM(d_0,\dots,d_{f-1};1)$ for some $\ofM$ as in \eqref{eq: kisinirreg} with $a=a^{\Theta}, b=b^{\Theta}$ and $x_i=x_i^{\Theta}$. By our assumptions on the set $J^{\mu}$, we see that $\ofM$ is as in Section \ref{irregularweightdescription}. By Remark \ref{diffxi} and Lemma \ref{lemma93} we find that $\ofM$ is the Kisin module corresponding to an element in $L_{\psi_1,\psi_2}$.

We have $\Tkis(\ofM^\Theta) \cong \Tkis(\ofM \otimes \ofM(d_0,\dots,d_{f-1};1))$, and by Lemma \ref{lem542} we find that we have an isomorphism of $G_K$-representations.
By Remark \ref{diffxi}, we know that different choices of the $x_i$ give different Galois representations. So by Lemma \ref{lemma93} both $L_{\psi_1,\psi_2}$ and $L^{\Theta}_{\text{sub}}$ have $|\F_L|^{| J \setminus (J \cap J_0)|}$ elements and it follows the spaces are the same (again using that we have compatible twists).
\end{proof}

\subsubsection{Consequences of $\rho$ having regular lifts}

In this section we assume that $\rho: G_K \to \GL_2(\Fpb)$ has crystalline lifts of weights $(k',l')$ and $\{k^{\mu},l^{\mu}\}_{\mu \in \tilde{M}}$, or $\rho: G_K \to \GL_2(\Fpb)$ has crystalline lifts of weights $(k',l')$ and $(k^{\Theta},l^{\Theta})$. We show that this implies that $c_{\rho}$ is an element of the subspace $L_{\psi_1,\psi_2}$ corresponding to the irregular weight $(k,0)$. In the proofs, we split into further cases depending on the set $J$ obtained from Proposition \ref{thetatoirregshape} and Proposition \ref{weightoneshape}.

We start with an auxiliary lemma relating the extensions of Kisin modules corresponding to the regular weights $(k',0)$ and $(k^{\mu},l^{\mu})$.

\begin{lemma}[Map between $\ofM^{\mu}$ and $\ofM'$]
\label{muprimemap}
Let $\mu \in \tilde{M}$. Let $J'$ and $J^{\mu}$ be subsets of $\Sigma$ such that $\mu \in J' \cap J^{\mu}$. Moreover, assume 
\[
\tau \in J' \iff \tau \in J^{\mu}
\]
for all $\tau \not \in (J_0 \cap \text{Block}_{\mu})$ and that for all $\tau \in \text{Block}_{\mu} \cap J_0$ we have $\tau \in J'$ and $\tau \not \in J^{\mu}$.

Let $\ofM^{\mu}$ be an extension of Kisin modules as described in \eqref{eq: kisinmu}. Then there is an extension $\ofM'$ of the type \eqref{eq: kisinprime} such that 
\[
\Tkis(\ofM^{\mu} \otimes \ofM(g_0^{\mu}, \dots,g_{f-1}^{\mu};1)) \cong \Tkis(\ofM' \otimes \ofM(d_0,\dots,d_{f-1};1)),
\]
where 
\[
g^{\mu}_i= \begin{cases}
1 & \text{if } \tau_i \in \tilde{M} \setminus \{\mu \},\\
0 & \text{otherwise},    
\end{cases} \quad
d_i= \begin{cases}
1 & \text{if } i \in \tilde{M}, \\
0 & {otherwise},    
\end{cases}
\]
and $\ofM'$ is such that 
\[
\varphi(e'_{i-1})=(b^{\mu})_i u^{t'_i} e'_i, \quad \varphi(f'_{i-1})=(a^{\mu})_i u^{s'_i} f'_i + x^{\mu}_i e'_i.
\]
\end{lemma}
\begin{proof}
We show the conditions of Proposition \ref{prop523reg} are satisfied. Let $\ofM'$ be as in \eqref{eq: kisinirreg} and consider the extension $\ofM^{\mu}$ as in the statement.
Write
\[
\ofN^{\mu}_{twist}:=\ofN^{\mu} \otimes \ofM(g_0^{\mu}, \dots,g_{f-1}^{\mu};1) \text{ and } \ofP^{\mu}_{twist}:=\ofP^{\mu} \otimes \ofM(g_0^{\mu}, \dots,g_{f-1}^{\mu};1),
\]
and 
\[
\ofN'_{twist}:=\ofN' \otimes \ofM(d_0, \dots,d_{f-1};1) \text{ and } \ofP'_{twist}:=\ofP' \otimes \ofM(d_0, \dots,d_{f-1};1).
\]
Note these twists add $1$ to $s'_i,t'_i$ at embeddings in $i \in \tilde{M}$. Recall that in Section \ref{regularweightdescription} we added 1 to $s^{\mu}_i,t^{\mu}_i$ for $i$ such that $\tau_i=\mu$ in order to apply Theorem \ref{thm79} (by twisting the representation by a fundamental character). Combining this with the twist by $\ofM(g_0^{\mu}, \dots,g_{f-1}^{\mu};1)$ above, we have now also added $1$ to $s^{\mu}_i,t^{\mu}_i$ at all embeddings in $i \in \tilde{M}$.
This means by Lemma \ref{alphaprimemu} and Lemma \ref{lemma512} we have non-zero maps $\ofP^{\mu}_{twist} \to \ofP'_{twist}$ and $ \ofN'_{twist} \to  \ofN^{\mu}_{twist}$.
The claim about the existence of a basis with properties as in Proposition \ref{prop523reg} follows from \eqref{eq: kisinmu}. So we obtain some $\ofM'_{twist} \in \Ext^1(\ofN'_{twist},\ofP'_{twist})$ such that
\[
\Tkis(\ofM^{\mu} \otimes \ofM(g_0^{\mu}, \dots,g_{f-1}^{\mu};1)) \cong \Tkis(\ofM'_{twist})
\]
with parameters
\[
y'_{twist,i}=u^{\alpha_i(\ofP'_{twist}, \ofP^{\mu}_{twist})+p(\alpha_{i-1}(\ofN^{\mu}_{twist}, \ofN'_{twist}))} y^{\mu}_{twist, i},
\]
by Proposition \ref{prop523reg}. Here 
\[
y^{\mu}_{twist,i}= 
\begin{cases}
x_i^{\mu}, & i \not \in \tilde{M} \setminus \{\mu\}, \\
u x_i^{\mu}, & i \in \tilde{M} \setminus \{\mu\}.
\end{cases}
\]
By definition we have $x_i^{\mu}=0$ if $i \not \in J^{\mu}$, so then $y_{twist,i}^{\mu}=0$. This gives $y'_{twist, i}=0$. Moreover, we have assumed $\mu \in J' \cap J^{\mu}$ and $i \not \in J^{\mu}$ for all $i \in J_0 \cap \text{Block}_{\mu}$. So also for such $i$ we have $x_i^{\mu}=0$, so again we have $y'_{twist,i}=0$.

Suppose $x_i^{\mu} \neq 0$, so that $i \in J^{\mu}$. Then by Lemma \ref{alphaprimemu} we have 
\[
\alpha_i(\ofP'_{twist}, \ofP^{\mu}_{twist})+p(\alpha_{i-1}(\ofN^{\mu}_{twist}, \ofN'_{twist}))=
\begin{cases}
1, & i=i_{\mu}, \\
0, & \text{otherwise}.
\end{cases}
\]
Now by the above we find
\[
y'_{twist,i}= 
\begin{cases}
x_i^{\mu}, & i \not \in \tilde{M}, \\
u x_i^{\mu}, & i \in \tilde{M}.    
\end{cases}
\]
Finally it is easy to see that we have $ \Tkis(\ofM'_{twist}) \cong \ofM' \otimes \ofM(d_0, \dots,d_{f-1};1)$ for $\ofM'$ as in \eqref{eq: kisinprime} with $a'=a^{\mu}, b'=b^{\mu}$ and $x'_i=x^{\mu}_i$. 
So we have
$\Tkis(\ofM^{\mu} \otimes \ofM(g_0^{\mu}, \dots,g_{f-1}^{\mu};1)) \cong \Tkis(\ofM' \otimes \ofM(d_0,\dots,d_{f-1};1))$. 
\end{proof}

We can now prove the result if $\rho$ has lifts of weights $(k',l')$ and $\{k^{\mu},l^{\mu}\}_{\mu \in \tilde{M}}$.

\begin{proposition}
\label{intersectionresult}
Suppose $\rho: G_K \to \GL_2(\Fpb)$ has crystalline lifts of weights $(k',l')$ and $\{k^{\mu},l^{\mu}\}_{\mu \in \tilde{M}}$. Let $J'$ and $\{J^{\mu}\}_{\mu \in \tilde{M}}$ be sets as in Corollary \ref{corextclass} such that for each set $J^{\mu}$ the following holds: for all $\tau \in \Sigma$ such that $\tau \not \in \text{Block}_{\mu} \cap J_0$ we have 
\[
\tau \in J' \iff \tau \in J^{\mu}.
\]
Let $J$ be as in Proposition \ref{weightoneshape}.

For each $\mu \in \tilde{M}$, consider $L_{\psi_1,\psi_2}$, $L_{\psi'_1,\psi'_2}$ and $L_{\psi^{\mu}_1, \psi^{\mu}_2}$ with $\psi_1,\psi_2$, $\psi'_1,\psi'_2$ and  $\psi^{\mu}_1,\psi^{\mu}_2$ corresponding to $(k,0)$, $(k',l')$ and $(k^{\mu},l^{\mu})$ respectively, as in \eqref{eq: charpsi}. Suppose $\overline{\psi_1 \psi_2^{-1}}= \overline{\psi'_1} \overline{\psi'_2}^{-1}=\overline{\psi^{\mu}_1} \overline{\psi^{\mu}_2}^{-1}$ for each $\mu \in \tilde{M}$ such that $\mu \in J'$. Then we have
\[
c_{\rho} \in L_{\psi_1,\psi_2}. 
\]
\end{proposition}

\begin{proof}
By Lemma \ref{nusemimatch}, Lemma \ref{primenotexceptional} and Lemma \ref{munotexceptional} we are not in the exceptional case for any of the weights. So we can consider the Kisin modules as in Section \ref{regularweightdescription}. 

Suppose we have an element in $L_{\psi'_1,\psi'_2}$. Then by Theorem \ref{thm79}, we obtain $\ofM'$ as in \eqref{eq: kisinprime} with parameters $x'_i$ and some constants $a', b' \in \F_L$. Here $x'_i=0$ if $i \not \in J'$. In particular, if $\mu \not \in J'$ for some $\mu \in \tilde{M}$, then by Lemma \ref{nusemimatch} we have $x'_i=0$ for all $i \in J_0 \cap \text{Block}_{\mu}$. We will use Proposition \ref{isoregprime}. We note the condition on the set $J'$ therein follows from our assumptions on $J$ together with Lemma \ref{nusemimatch}. 

We write $c_{\rho} \in H^1(G_K,\Fpb(\overline{\psi}))$ for the extension class corresponding to $\rho$. Consider $c_{\rho}$ as an element of $L_{\psi'_1,\psi'_2}$ and the corresponding extension $\ofM'$ of Kisin modules as in \eqref{eq: kisinprime} with fixed $a',b' \in \F_L$. Then if $\mu \not \in J'$ for all $\mu \in \tilde{M}$ we have $x'_i=0$ for all $i \in J_0$ and so the result follows from Proposition \ref{isoregprime}.

Next suppose there is some $\mu \in \tilde{M}$ such that $\mu \in J'$. Consider $c_{\rho} \in L_{\psi^{\mu}_1,\psi^{\mu}_2}$ and the corresponding extension $\ofM^{\mu}$ of Kisin modules as in \eqref{eq: kisinmu}. Assume $a^{\mu}=a'$ and $b^{\mu}=b'$ (possibly after twisting the extension by $\ofM(0,\dots,0;c)$ where $c$ is such that $\Tkis(\ofM(0,\dots,0;c)) = (\overline{\psi'_1 \psi^{\mu}_1}^{-1})|_{G_{\infty}}=(\overline{\psi_2 \psi^{\mu}_2}^{-1})|_{G_{\infty}}$)). Write $\ofM'_{\mu}$ for the extension as in $\eqref{eq: kisinprime}$ above that we obtain by Lemma \ref{muprimemap}. Then letting
\[
\ofM'_{\text{twist},\mu}:= \ofM'_{\mu} \otimes \ofM(d_0,\dots,d_{f-1};1)),
\]
we have
\[
\Tkis(\ofM^{\mu} \otimes \ofM(g_0^{\mu}, \dots,g_{f-1}^{\mu};1)) \cong \Tkis(\ofM'_{\text{twist},\mu}),
\]
where
\[
g^{\mu}_i= \begin{cases}
1 & \text{if } \tau_i \in \tilde{M} \setminus \{\mu \},\\
0 & \text{otherwise},    
\end{cases} \quad
d_i= \begin{cases}
1 & \text{if } i \in \tilde{M}, \\
0 & \text{otherwise}.   
\end{cases}
\]
For the constants $a'_{\mu}, b'_{\mu} \in \F_L$ of $\ofM'_{\mu}$ (and thus $\ofM'_{\text{twist},\mu}$), we have $a'_{\mu}=a^{\mu}=a'$ and $b'_{\mu}=b^{\mu}=b'$. For the parameters, writing $x'_{\mu,i}$, we have $x'_{\mu,i}=x^{\mu}_i$ by Lemma \ref{muprimemap}. This gives $x'_{\mu,i}=0$ if $i \in J_0 \cap \text{Block}_{\mu}$ by Lemma \ref{nusemimatch}, and for $i \not \in J_0 \cap \text{Block}_{\mu}$ we have $x'_{\mu,i} = 0$ if $i \not \in J'$ by our assumptions. 

In this way, for each $\mu \in \tilde{M}$ such that $\mu \in J'$ we obtain an extension $\ofM'_{\text{twist},\mu}$ of two rank one Kisin modules
\begin{align}
\label{eq: classicIIIkisin}
0 \to \ofM(\tilde{\wbottom}'_0,\dots,\tilde{\wbottom}'_{f-1}; b') \to \ofM'_{\text{twist},\mu} \to \ofM(\tilde{\wtop}'_0,\dots,\tilde{\wtop}'_{f-1}; a') \to 0,
\end{align}
where 
\[
\tilde{\wbottom}'_i = 
\begin{cases}
\wbottom'_i +1, & \tau_i \in \tilde{M},\\  
\wbottom'_i, & \tau_i \not \in \tilde{M}, 
\end{cases}
\quad
\tilde{\wtop}'_i=
\begin{cases}
\wtop'_i +1, & \tau_i \in \tilde{M},  \\
\wtop'_i, & \tau_i \not \in \tilde{M}. 
\end{cases}
\]
We have an action of the Frobenius on the bases $\{e'_i,f'_i\}$ as follows:
\begin{align}
\label{eq: basisfortwistmu}    
\varphi(e'_{i-1})=
  (b' )_i u^{\tilde{\wbottom}'_i} e'_i, \quad 
\varphi(f'_{i-1})=
\begin{cases}    
(a' )_i u^{\tilde{\wtop}'_i} f'_i + x'_{\mu,i} e'_i, & i \not \in \tilde{M}, \\
(a' )_i u^{\tilde{\wtop}'_i} f'_i + u x'_{\mu,i} e'_i, & i \in \tilde{M}.
\end{cases}
\end{align}
Note only the $x'_{\mu,i}$ are dependent on the choice of $\mu \in \tilde{M}$.

For each $\mu \in \tilde{M}$ we have 
\[
\hat{T}(\hofM(g_0^{\mu}, \dots,g_{f-1}^{\mu};1)) \otimes \omega_{\mu} \cong  \hat{T}(\hofM(d_0,\dots, d_{f-1};1)).
\]
Since we have $c_{\rho} \in L_{\psi'_1,\psi'_2}$ and $ c_{\rho} \in L_{\psi_1^{\mu},\psi_2^{\mu}}$, again writing $\ofM'$ for the extension corresponding to $c_{\rho} \in L_{\psi'_1,\psi'_2}$, we find
\begin{align*}
\left( \rho \otimes \prod_{\tau \in \tilde{M}} \omega_{\tau} \right)_{\big\vert_{G_{\infty}}} & \cong \Tkis(\ofM' \otimes \ofM(d_0,\dots, d_{f-1};1)), \\
&\cong \Tkis (\ofM^{\mu} \otimes \ofM(g_0^{\mu}, \dots, g_{f-1}^{\mu}; 1), \\ 
&\cong \Tkis (\ofM'_{\text{twist},\mu}).
\end{align*}

By Lemma \ref{lem542} the above will give isomorphisms of $G_K$-representations $\hat{T}(\hat{\ofM'} \otimes \hofM(d_0,\dots, d_{f-1};1)) \cong \hat{T}(\ofM'_{\text{twist},\mu}) $

The extensions $\ofM'_{\text{twist},\mu}$ and $\ofM' \otimes \ofM(d_0,\dots, d_{f-1};1)$ are both as in \eqref{eq: classicIIIkisin} with the action of $\varphi$ exactly as in \eqref{eq: basisfortwistmu}, except that the parameters $x'_i$ corresponding to $\ofM' \otimes \ofM(d_0,\dots, d_{f-1};1)$ could be different than the $x'_{\mu,i}$ listed there. But by the isomorphisms above and Remark \ref{diffxi} we must have $x'_{\mu,i}= x'_i$ for all $i$. So in particular, we obtain $x'_i=0$ for $i \in J_0 \cap \text{Block}_{\mu}$.

Since $c_{\rho} \in \cap_{\mu \in \tilde{M}} L_{\psi_1^{\mu},\psi_2^{\mu}}$, we
repeat this for all $\mu \in \tilde{M}$ such that $\mu \in J'$. Note we already had $x'_i=0$ for all $i \in J_0 \cap \text{Block}_{\mu}$ for $\mu \in \tilde{M}$ such that $\mu \not \in J'$.
But then we must have that
\[
\text{$x'_i=0$ if $\tau_i \in J_0$ and $x'_i=0$ if $i \not \in J'$}.
\]
By Proposition \ref{isoregprime}, we are done.
\end{proof}

\subsubsection{Two weight version}
\label{twoweightsection}
We now return to the other set of regular weights, $(k',l')$ and $(k^{\Theta},l^{\Theta})$, occurring in the two weight version of the main theorem. We need some more auxiliary definitions and results. We first introduce one more space of extensions. This space is chosen so that we can compare the extensions corresponding to $(k',l')$ to those corresponding to $(k^{\Theta},l^{\Theta})$.

\begin{definition}[Auxiliary space]
\label{auxiliarydescription}
We fix subsets $J'$ and $J^{\Theta}$ of $\Sigma$ such that 
\[
\tau \in J' \iff \tau \in J^{\Theta}
\]
for all $\tau \not \in J_0$. Moreover, we assume that for all $\nu \in \tilde{M}$ we have 
\[
\tau_{i_{\nu}}, \tau_{i_{\nu}+1}, \dots, \tau_E\in J' \text{ and }\tau_{i_\nu+1},\dots, \tau_{i_E} \not \in J^{\Theta},
\]
or
\[
\tau_{i_{\nu}}, \tau_{i_{\nu}+1}, \dots, \tau_E \not \in J' \text{ and }\tau_{i_\nu+1},\dots, \tau_{i_E}  \in J^{\Theta}.
\]
where $\tau_{i_{\nu}}, \dots, \tau_E$ is the part of $\text{Block}_{\nu}$ starting at $\tau_{i_{\nu}}=\nu$ and where $\tau_E$ is the end of the block. 

We use the notation as in \eqref{eq: rexpressions} and \eqref{eq: sexpressions}. We let $\psi^{\gamma}_1,\psi^{\gamma}_2: G_K \to \overline{\Z}_p^{\times}$ be two crystalline characters such that
\begin{align}
\label{eq: gammaweights}
\HT(\psi^{\gamma}_1)_{i}:= \wtop^{\gamma}_i =
\begin{cases}
    \wtop_i', & i \in J', \\
    \wtop^{\Theta}_i, & i \not \in J',
\end{cases}
\quad
\HT(\psi^{\gamma}_2)_{i} := \wbottom^{\gamma}_i =
\begin{cases}
    \wbottom_i', & i \in J',\\
    \wbottom^{\Theta}_i, & i \not \in J'.
\end{cases}
\end{align}
These exist by Lemma \ref{cryslift} and are unique up to an unramified twist. Let $L_{\psi^{\gamma}_1,\psi^{\gamma}_2}$ be the corresponding subspace.
\end{definition}

We write down an example before we study the corresponding Kisin modules.

\begin{example}[Example \ref{explicitweight} continued]
As before we let $f=4$ and consider the irregular weight $(k,l)=((1,k_1,1,k_3),(0,0,0,0))$ with $k_1,k_3 \geq 3$.
Suppose $J'=\{\tau_1,\tau_2\}$ and $J^{\Theta}=\{\tau_0,\tau_1\}$. Then we find 
\[
(\wtop_0^{\gamma},\wtop_1^{\gamma},\wtop_2^{\gamma},\wtop_3^{\gamma})=(p,k_1-1,p,-1) \text{ and } (\wbottom_0^{\gamma},\wbottom_1^{\gamma},\wbottom_2^{\gamma},\wbottom_3^{\gamma})=(0,0,0,k_3-1).
\]
\end{example}

\noindent \textbf{Kisin modules coming from $L_{\psi^{\gamma}_1,\psi^{\gamma}_2}$.}
For each element in $L_{\psi^{\gamma}_1,\psi^{\gamma}_2}$ we consider the corresponding Galois representation and twist it by $\prod_{\mu \in \tilde{M}, \mu \not \in J'} \omega_{\mu}$. This removes the negative Hodge--Tate weights at these embeddings and we can now apply Theorem \ref{thm79}. 

We assume we are not in the exceptional case of Theorem \ref{thm79}. We write $\ofM^{\gamma}$ for the obtained extension of Kisin modules:
\begin{align}
\label{eq: kisinnew}
0 \to \ofM(\tilde{t}^{\gamma}_0,\dots,\tilde{t}^{\gamma}_{f-1};b^{\gamma}) \to \ofM^{\gamma} \to \ofM(\tilde{s}^{\gamma}_0,\dots,\tilde{s}^{\gamma}_{f-1};a^{\gamma}) \to 0,
\end{align}
where
\begin{align*}
    \tilde{s}^{\gamma}_i = 
    \begin{cases}
    \wtop_i', & i \in J', \\
    \wtop^{\Theta}_i, & i \not \in J',\\
    \wtop^{\Theta}_i+1, & i \not \in J', i \in \tilde{M},
\end{cases}
\quad
    \tilde{t}^{\gamma}_i =
    \begin{cases}
    \wbottom_i', & i \in J',\\
    \wbottom^{\Theta}_i, & i \not \in J',\\
    \wbottom^{\Theta}_i+1, & i \not \in J', i \in \tilde{M},
\end{cases}
\end{align*}
and where $a^{\gamma}, b^{\gamma} \in \F_L$. We obtain bases $\{e^{\gamma}_i,f^{\gamma}_i\}$ for $\ofM^{\gamma}_i$ such that
\begin{align*}
    \varphi(e^{\gamma}_{i-1}) = (b^{\gamma})_i u^{\tilde{t}^{\gamma}_i} e^{\gamma}_i, \quad    \varphi(f^{\gamma}_{i-1}) = (a^{\gamma})_i u^{\tilde{s}^{\gamma}_i} f^{\gamma}_i + x^{\gamma}_i e^{\gamma}_i,
\end{align*}
where $x^{\gamma}_i \in \F_L$ and $x^{\gamma}_i=0$ if both $\tau_i \not \in J'$ and $\tau_i \in \tilde{M}$, or if both $\tau_i \not \in J'$ and $\tau_i \not \in J_0 \cup \tilde{M}$.
In the following we write $\ofN^{\gamma}$ and $\ofP^{\gamma}$ for the rank one modules above, so that $\ofM^{\gamma} \in \Ext^1(\ofN^{\gamma},\ofP^{\gamma})$.

We will again use Proposition \ref{prop523reg} to relate the Kisin modules as above to our earlier regular extensions of Kisin modules.

\begin{lemma}[Map from $\ofM^{\Theta}$ to $\ofM^{\gamma}$]
\label{thetagammamap}
Let $J'$ and $J^{\Theta}$ be subsets of $\Sigma$ as in Definition \ref{auxiliarydescription}. 
Let $\ofM^{\Theta}$ be an extension of Kisin modules as described in \eqref{eq: kisintheta}. Then there is an extension $\ofM^{\gamma}$ of the type as described in \eqref{eq: kisinnew} such that 
\[
\Tkis(\ofM^{\Theta}) \cong \Tkis(\ofM^{\gamma} \otimes \ofM(g_0,\dots,g_{f-1};1)),\]
where $g_i=1$ if $\tau_i \in \tilde{M}, \tau_i \in J'$ and $0$ otherwise. Moreover, $\ofM^{\gamma}$ is such that 
\[
\varphi(e^{\gamma}_{i-1})=(b^{\Theta})_i u^{\tilde{t}^{\gamma}_i} e^{\gamma}_i, \quad \varphi(f^{\gamma}_{i-1})=(a^{\Theta})_i u^{\tilde{s}^{\gamma}_i} f^{\gamma}_i + x^{\Theta}_i e^{\gamma}_i.
\]
\end{lemma}
\begin{proof}
Consider the rank one modules $\ofN^{\gamma}$ and $\ofP^{\gamma}$ with $a^{\gamma}=a^{\Theta}$ and $b^{\gamma} = b^{\theta}$. Now set 
\[
\ofN^{\gamma}_{twist}:=\ofN^{\gamma} \otimes \ofM(g_0, \dots,g_{f-1};1) \text{ and }\ofP^{\gamma}_{twist}:=\ofP^{\gamma} \otimes \ofM(g_0, \dots,g_{f-1};1).
\]
Note that before applying Theorem \ref{thm79}, we twisted any representation in $L_{\psi_1^{\gamma},\psi_2^{\gamma}}$ by $\prod_{\mu \in \tilde{M},\mu \not \in J'} \omega_{\mu}$. Combined with the above twists, this means we have added 1 to $s^{\gamma}_i,t^{\gamma}_i$ at embeddings $\tau_i$ such that $\tau_i \in \tilde{M}$. Similarly, we have added $1$ to $s^{\Theta}_i,t^{\Theta}_i$ when $\tau_i \in \tilde{M}$ in Section \ref{regularweightdescription} by twisting by a product of fundamental characters.

Consider the rank one modules $\ofN^{\Theta}$ and $\ofP^{\Theta}$ as in \eqref{eq: kisintheta}. By Lemma \ref{alphaprimetheta} and Lemma \ref{lemma512}, we have non-zero maps $\ofN^{\gamma}_{twist} \to \ofN^{\Theta}$ and $\ofP^{\Theta} \to \ofP^{\gamma}_{twist}$. The condition on the basis in Proposition \ref{prop523reg} for $\ofM^{\Theta}$ follows from \eqref{eq: kisintheta}. 
So by Proposition \ref{prop523reg} we obtain some $\ofM^{\gamma}_{twist} \in \Ext^1(\ofN^{\gamma}_{twist},\ofP^{\gamma}_{twist})$ such that
\[
\Tkis(\ofM^{\gamma}_{twist})  \cong  \Tkis(\ofM^{\Theta}),
\]
with parameters $y^{\gamma}_{twist,i}$ such that
\[y^{\gamma}_{twist,i}=u^{\alpha_i(\ofP^{\Theta},\ofP^{\gamma}_{twist})+p(\alpha_{i-1}(\ofN^{\gamma}_{twist}, \ofN^{\Theta}))} y^{\Theta}_i,
\]
with $y^{\Theta}_i= x^{\Theta}_i$. We have $x_i^{\Theta}=0$ if $i \not \in J^{\Theta}$, so that in those cases $y^{\gamma}_{twist,i}=0$ and thus $x_i^{\gamma}=0$.

We now suppose $x_i^{\Theta} \neq 0$, so that we must have $i \in J^{\Theta}$, then by Lemma \ref{alphaprimetheta} we have
\[
\alpha_i(\ofP^{\Theta}, \ofP^{\gamma}_{twist})+p(\alpha_{i-1}(\ofN^{\gamma}_{ twist}, \ofN^{\Theta})) = 
\begin{cases}
1, & i \in \tilde{M,} \\
0, & \text{otherwise},
\end{cases}
\]
where for $i \in J_0$ we use that $i \in J^{\Theta}$ implies $i \not \in J'$ by our assumptions. By the above this gives 
\[
y^{\gamma}_{twist,i}= 
\begin{cases}
x_i^{\Theta}, & i \not \in \tilde{M},  \\
u x_i^{\Theta}, & i \in \tilde{M}.
\end{cases}
\]
We note that $x_i^{\Theta}=0$ if $i \not \in  J'$ for $i \in \tilde{M}$. It is then easy to see that we have $\Tkis(\ofM^{\gamma} \otimes \ofM(g_0, \dots,g_{f-1};1)) \cong \Tkis(\ofM^{\gamma}_{twist})$ for $\ofM^{\gamma}$ as in \eqref{eq: kisinnew} with $a'=a^{\Theta}, b'=b^{\Theta}$ and $x'_i=x^{\Theta}_i$ so the result follows. 
\end{proof}

\begin{lemma}[Map from $\ofM'$ to $\ofM^{\gamma}$]
\label{primegammamap}
Let $J'$ and $J^{\Theta}$ be subsets of $\Sigma$ as in Definition \ref{auxiliarydescription}.
Let $\ofM'$ be an extension of Kisin modules as described in \eqref{eq: kisinprime}. Then there is an extension $\ofM^{\gamma}$ as in \eqref{eq: kisinnew} such that 
\[
\Tkis(\ofM' \otimes \ofM(d_0,\dots,d_{f-1};1) ) \cong \Tkis(\ofM^{\gamma} \otimes \ofM(g_0,\dots,g_{f-1};1)),
\]
with 
\[
d_i = \begin{cases}
1, &\text{if } \tau_i \in \tilde{M},\\
0, & \text{otherwise},
\end{cases}
\text{ and }
g_i=\begin{cases}
1, & \text{if } \tau_i \in \tilde{M}, \tau_i \in J',\\
0, & \text{otherwise},
\end{cases}
\]
and $\ofM^{\gamma}$ is such that 
\[
\varphi(e'_{i-1})=(b')_i u^{\tilde{t}^{\gamma}_i} e'_i, \quad \varphi(f'_{i-1})=(a')_i u^{\tilde{s}^{\gamma}_i} f'_i + x'_i e'_i.
\]
\end{lemma}

\begin{proof}
Let $\ofM'$ be as in \eqref{eq: kisinprime}. Write 
\[
\ofM'_{twist}:=\ofM' \otimes \ofM(d_0, \dots,d_{f-1};1),
\]
and write $\ofN'_{twist}$ and $\ofP'_{twist}$ for the corresponding twisted rank one modules.  Consider the rank one modules $\ofN^{\gamma}$ and $\ofP^{\gamma}$ with $a^{\gamma}=a'$ and $b^{\gamma} = b'$. Now set 
\[
\ofN^{\gamma}_{twist}:=\ofN^{\gamma} \otimes \ofM(g_0, \dots,g_{f-1};1) \text{ and }\ofP^{\gamma}_{twist}:=\ofP^{\gamma} \otimes \ofM(g_0, \dots,g_{f-1};1).
\]
The twists add $1$ to $s'_i,t'_i$ at embeddings $\tau_i$ such that $\tau_i \in \tilde{M}$. Note that before applying Theorem \ref{thm79}, we twisted any representation in $L_{\psi_1^{\gamma},\psi_2^{\gamma}}$ by $\prod_{\mu \in \tilde{M},\mu \not \in J'} \omega_{\mu}$. Combined with the above twists, this means we have added 1 to $s^{\gamma}_i,t^{\gamma}_i$ at embeddings $\tau_i$ such that $\tau_i \in \tilde{M}$. By Lemma \ref{alphaprimetheta} and Lemma \ref{lemma512}, we have non-zero maps $\ofN^{\gamma}_{twist} \to \ofN'_{twist}$ and $\ofP'_{twist} \to \ofP^{\gamma}_{twist}$.

We check the remaining conditions of Proposition \ref{prop523reg}. The condition on the basis for $\ofM'_{twist}$ follows from \eqref{eq: kisinprime}. 
By Proposition \ref{prop523reg} we obtain some $\ofM^{\gamma}_{twist} \in \Ext^1(\ofN^{\gamma}_{twist},\ofP^{\gamma}_{twist})$ such that 
\[
\Tkis(\ofM^{\gamma}_{twist}) \cong \Tkis(\ofM' \otimes \ofM(d_0,\dots,d_{f-1};1)),
\]
with parameters $y^{\gamma}_{twist,i}$ such that
\[
y^{\gamma}_{twist,i}=u^{\alpha_i(\ofP'_{twist},\ofP^{\gamma}_{twist})+p(\alpha_{i-1}(\ofN^{\gamma}_{twist},\ofN'_{twist}))} y'_{twist,i}.
\]
Here 
\[
y'_{twist,i}= 
\begin{cases}
x'_i, & i \not \in \tilde{M},  \\
u x'_i, & i \in \tilde{M}. 
\end{cases}
\]
We have $x'_i=0$ if $i \not \in J'$, in which case clearly $y^{\gamma}_{twist,i}=0$. Now suppose $i \in J'$, then by Lemma \ref{alphaprimetheta} we find 
\[
\alpha_i(\ofP'_{twist}, \ofP^{\gamma}_{twist})+p(\alpha_{i-1}(\ofN^{\gamma}_{twist}, \ofN'_{twist}))=0.\]
This gives us
\[
y^{\gamma}_{twist,i}= 
\begin{cases}
x'_i, & i \not \in \tilde{M} \cap J', \\
u x'_i, & i \in \tilde{M} \cap J'.
\end{cases}
\]
It is easy to see that we have $ \Tkis(\ofM^{\gamma}_{twist}) \cong \ofM^{\gamma} \otimes \ofM(g_0, \dots,g_{f-1};1)$ for $\ofM^{\gamma}$ as in \eqref{eq: kisinnew} with $a^{\gamma}= a', b^{\gamma} = b', x^{\gamma}_i = x'_i$ as before and the result follows. 
\end{proof}

Now we are ready to prove the result we need for the two weight version. 
\begin{proposition}
\label{thetaresult}
Suppose $\rho: G_K \to \GL_2(\Fpb)$ is a reducible representation with crystalline lifts of weights $(k',0)$ and $(k^{\Theta},l^{\Theta})$. Let $J'$ and $J^{\Theta}$ be as in Corollary \ref{corextclass} and as in Definition \ref{auxiliarydescription}.

Let $J$ be as in Proposition \ref{thetatoirregshape} and consider $L_{\psi_1,\psi_2}$, $L_{\psi'_1, \psi'_2}$,  $L_{\psi^{\Theta}_1,\psi^{\Theta}_2}$ and $L_{\psi^{\gamma}_1,\psi^{\gamma}_2}$ as in \eqref{eq: charpsi}. Suppose $\overline{\psi_1 \psi_2^{-1}}=\overline{\psi^{\gamma}_1} \overline{\psi^{\gamma}_2}^{-1}=\overline{\psi^{\Theta}_1} \overline{\psi^{\Theta}_2}^{-1}=\overline{\psi'_1 \psi'_2}^{-1}$. Then we have
\[
c_{\rho} \in L_{\psi_1,\psi_2}. 
\]
\end{proposition}

\begin{proof}
By Lemma \ref{primenotexceptional}, Lemma \ref{thetanotexceptional} we are not in the exceptional case of Theorem \ref{thm79} for any of the weights $(k',0)$ and $(k^{\Theta},l^{\Theta})$. 

Suppose we have an element in $L_{\psi'_1,\psi'_2}$. Then by Theorem \ref{thm79}, we obtain $\ofM'$ as in \eqref{eq: kisinprime} with parameters $x'_i$ and constants $a'$ and $b' \in \F_L$. Here $x'_i=0$ if $i \not \in J'$. In particular, if $\mu \not \in J'$ for some $\mu \in \tilde{M}$, then by Lemma \ref{thetasemimatch} we have $x'_i=0$ for all $i \in J_0 \cap \text{Block}
_{\mu}$. We will use Proposition \ref{isoregprime}. We note the condition on the set $J'$ therein follows from our assumptions on $J$ together with Lemma \ref{thetasemimatch}. 

Consider $c_{\rho} \in L_{\psi'_1,\psi'_2}$ and the corresponding extension $\ofM'$ of Kisin modules as in \eqref{eq: kisinprime}. Suppose first that $\mu \not \in J'$ for all $\mu \in \tilde{M}$.  Then as above, by Lemma \ref{thetasemimatch} again, we have $x'_i=0$ for all $i \in J_0$. The result follows from Proposition \ref{isoregprime}.

Next suppose there is some $\mu \in \tilde{M}$ with $\mu \in J'$.
By Lemma \ref{primegammamap}, we have an extension $\ofM^{\gamma,1} \in \Ext^1(\ofN^{\gamma}, \ofP^{\gamma})$ as in $\eqref{eq: kisinnew}$ with $x_i^{\gamma}=x'_i$ such that 
\[
\Tkis(\ofM' \otimes \ofM(d_0,\dots, d_{f-1};1)) \cong \Tkis(\ofM^{\gamma,1} \otimes \ofM(g_0,\dots, g_{f-1};1)),
\]
where
\[
d_i = \begin{cases}
1, &\text{if } \tau_i \in \tilde{M},\\
0, & \text{otherwise},
\end{cases}
\text{ and }
g_i=\begin{cases}
1, & \text{if } \tau_i \in \tilde{M}, \tau_i \in J',\\
0, & \text{otherwise}.
\end{cases}
\]
Recall that for each $\mu \in \tilde{M}$, we have $x'_i=0$ if $i \in J_0 \cap \text{Block}_{\mu}$ if $\mu \not \in J'$ by Lemma \ref{thetasemimatch}.

Next consider $c_{\rho}$ as an element in $L_{\psi^{\Theta}_1,\psi^{\Theta}_2}$ and the corresponding Kisin extension $\ofM^{\Theta}$ as in \eqref{eq: kisintheta} (after twisting by $ \otimes_{\mu \in \tilde{M}} \omega_{\mu}$).  Assume $a^{\Theta}=a'$ and $b^{\Theta}=b'$ (possibly after twisting the extension by $\ofM(0,\dots,0;c)$ where $c$ is such that $\Tkis(\ofM(0,\dots,0;c)) = (\overline{\psi'_1}\overline{\psi^{\Theta}_1}^{-1})|_{G_{\infty}}=(\overline{\psi'_2}\overline{\psi^{\Theta}_2}^{-1})|_{G_{\infty}}$). By Lemma \ref{thetagammamap}, we then obtain an extension $\ofM^{\gamma,2}  \in \Ext^1(\ofN^{\gamma}, \ofP^{\gamma})$ with $x_i^{\gamma,2}=x^{\Theta}_i$ such that 
\[
\Tkis(\ofM^{\Theta}) \cong \Tkis(\ofM^{\gamma,2} \otimes \ofM(g_0,\dots, g_{f-1};1)),
\]
with $g_i$ as before. Recall that for each $\mu \in \tilde{M}$, we have $x^{\Theta}_i=0$ if $i \in J_0 \cap \text{Block}_{\mu}$ if $\mu \in J'$ by Lemma \ref{thetasemimatch}.

Since we have
\[
\hat{T}(\hofM(d_0,\dots,d_{f-1};1)) \cong \otimes_{\mu \in \tilde{M}} \omega_{\mu},
\]
we obtain
\[
 \Tkis(\ofM^{\gamma,1} \otimes \ofM(g_0,\dots, g_{f-1};1)) \cong  \Tkis(\ofM^{\gamma,2} \otimes \ofM(g_0,\dots, g_{f-1};1)).
\]
By Remark \ref{diffxi}, this implies $x_i^{\gamma,2}=x_i^{\gamma,1}$. Since we also have $x_i^{\gamma,1}=x'_i$ and $x_i^{\gamma,2}=x^{\Theta}_i$, this implies $x'_i=x^{\Theta}_i$ for all $i$. For this to hold we need $x_i^{\gamma,1}=x_i^{\gamma,2}=0$ for all $i \in J_0$.

This means the element $\ofM'$ corresponding to $c_\rho$ has to satisfy $x'_i=0$ for all $i \in J_0$.  This means the Galois representation associated to $c_\rho$ is in $L'_{sub}$. The result follows by Proposition \ref{isoregprime}.
\end{proof}


\subsection{Wrapping up the proof}

We can now prove both versions of our main result in the reducible case.

\begin{theorem}
\label{twoproof}
Let $K$ be an unramified finite extension of $\Qp$. Let $(k,0) \in \Z^{\Sigma} \times \Z^{\Sigma}$ be an irregular weight satisfying $1 \leq k_{\tau} \leq p$ for all $\tau \in \Sigma$ and such that $k_{\tau} \neq 2$ if $k_{\Frob^{-1} \circ \tau}=1$.  Suppose further $k \neq 1$. Suppose $\rho: G_K \to \GL_2(\Fpb)$ is a reducible representation. Then $\rho$ has a crystalline lift of weight $({k},{0})$ if and only if it has crystalline lifts of weights $(k',l')$ and $(k^{\Theta},l^{\Theta})$.
\end{theorem}

\begin{proof}
We have already proved this for $\rho$ semisimple in see Proposition \ref{semisimpleproof}, so assume $\rho$ is not semisimple. 

First assume $\rho: G_K \to \GL_2(\Fpb)$ has a crystalline lift of weight $({k},{0})$. By Proposition \ref{regularshape}, we know that $\rho|_{I_K}$ can be written to have the shape corresponding to the weights $(k',l')$ and $(k^{\Theta},l^{\Theta})$ for some sets $J'$ and $J^{\Theta}$. It remains to show that $c_{\rho} \in L_{\psi_1',\psi_2'}$ and $c_{\rho} \in L_{\psi_1^{\Theta},\psi_2^{\Theta}}$ for these sets. This follows from Proposition \ref{isoregprime} and Proposition \ref{isoregtheta}.

Next assume $\rho$ has crystalline lifts of weights $(k',l')$ and $(k^{\Theta},l^{\Theta})$. We choose sets $J'$ and $J^{\Theta}$ as in the statement of Proposition \ref{thetaresult} (which we can do by Lemma \ref{minimalthetamatch} and Lemma \ref{thetasemimatch}). By Proposition \ref{thetatoirregshape} we can rewrite $\rho|_{I_K}$ to have the shape corresponding to the weight $(k,0)$ for a set $J$ defined in this proposition. It remains to show that $c_{\rho} \in L_{\psi_1,\psi_2}$, the space corresponding to the weight $(k,0)$ and the subset $J$.
This follows from Proposition \ref{thetaresult}.
\end{proof}

\begin{theorem}
\label{reduciblethm}
Let $K$ be an unramified finite extension of $\Qp$. Let $(k,0) \in \Z^{\Sigma} \times \Z^{\Sigma}$ be an irregular weight satisfying $1 \leq k_{\tau} \leq p$ for all $\tau \in \Sigma$ and such that $k_{\tau} \neq 2$ if $k_{\Frob^{-1} \circ \tau}=1$.  Suppose further $k \neq 1$. Suppose $\rho: G_K \to \GL_2(\Fpb)$ is a reducible representation.
Then $\rho$ has a crystalline lift of weight $(k,0)$ if and only if it has crystalline lifts of weights $(k',l')$ and $(k^{\mu},l^{\mu})$ for each $\mu \in \tilde{M}$. 
\end{theorem}

\begin{proof}
By Proposition \ref{semisimpleproof} we can assume $\rho$ is not semisimple. 

First assume $\rho: G_K \to \GL_2(\Fpb)$ has a crystalline lift of weight $(k,0)$.
By Proposition \ref{regularshape}, we know that $\rho|_{I_K}$ can be written to have the shape corresponding to the weights $(k',l')$ and $(k^{\mu},l^{\mu})$  for each $\mu \in \tilde{M}$ for some sets $J'$ and $\{J^{\mu}\}_{\mu \in \tilde{M}}$. It remains to show that $c_{\rho} \in L_{\psi'_1,\psi'_2}$ and $L_{\psi^{\mu}_1,\psi^{\mu}_2}$ for these sets. This follows from Proposition \ref{isoregprime} and Proposition \ref{isoregmu}. 

Next assume $\rho$ has crystalline lifts of weights $(k',l')$ and $(k^{\mu},l^{\mu})$ for each $\mu \in \tilde{M}$. We choose sets $J'$ and $\{J^{\mu}\}_{\mu \in \tilde{M}}$ as in the statement of Proposition \ref{intersectionresult} (which we can do by Lemma \ref{minimalmumatch}). By Proposition \ref{weightoneshape} we know that $\rho|_{I_K}$ can be written to have the shape corresponding to the weight $(k,0)$ for some set $J$. It remains to show that the extension class $c_{\rho}$ lies in $L_{\psi_1,\psi_2}$ for this $J$.
This follows from Proposition \ref{intersectionresult} so that $\rho$ indeed has a crystalline lift of weight $(k,0)$.
\end{proof}

\section{The irreducible case}

Finally, we tackle the irreducible case. We separate the two directions. Recall that we write $K'$ for the quadratic unramified extension of $K$, $\F'$ for the residue field and $\Sigma'=\{ \tau: \F' \hookrightarrow \Fpb \}$. Finally let $\pi: \Sigma' \to \Sigma$ be defined by $\sigma \mapsto \sigma|_{\F}$.

\subsection{Obtaining regular lifts}

We start with the easier direction from irregular weights to regular weights.

\begin{proposition}
\label{irreducibletoregular}
Let $K$ be an unramified finite extension of $\Qp$. Let $(k,0) \in \Z^{\Sigma} \times \Z^{\Sigma}$ be an irregular weight satisfying $1 \leq k_{\tau} \leq p$ for all $\tau \in \Sigma$ and such that $k_{\tau} \neq 2$ if $k_{\Frob^{-1} \circ \tau}=1$.  Suppose further $k \neq 1$. If an irreducible representation $\rho: G_K \to \GL_2(\Fpb)$ has a crystalline lift of weights $(k,0)$, then $\rho$ also has crystalline lifts of weights $(k',0)$, $(k^{\Theta},l^{\Theta})$, and $(k^{\mu},l^{\mu})$ for all $\mu \in \tilde{M}$.   
\end{proposition}

\begin{proof}
Note first the weight $(k',0)$ is obtained from $(k,0)$ by adding the vectors $h_{\tau}= p e_{\Frob^{-1} \circ \tau}-e_{\tau}$ to $k \in \Z^{\Sigma}$ for each $\tau \in M$. The other weights are obtained similarly, as in Section \ref{regularweightdescription}.

By Theorem \ref{thm101} we have $\rho \sim \Ind_{G_{K'}}^{G_K} \xi$ with
    \[
    \xi|_{I_{K'}}=\prod_{\sigma \in J} \omega_{\sigma}^{k_{\pi(\sigma)}-1},
    \]
for some balanced subset $J \subset \Sigma'$.

Suppose $\nu \in \tilde{M}$. Choose $\sigma_0 \in \Sigma'$ such that $\sigma_0 \in J$ and $\pi(\sigma_0)=\nu$. Then rewrite $\xi|_{I_{K'}}$ as follows:
\[
    \xi|_{I_{K'}}=\omega_{\sigma_0}^{-1} \omega_{\Frob^{-1} \circ \sigma_0}^p \cdot \prod_{\sigma \in J} \omega_{\sigma}^{k_{\pi(\sigma)}-1}.
\]
We see $\xi$ has a lift $\tilde{\xi}$ such that $\Ind_{G_{K'}}^{G_K} \tilde{\xi}$ is a crystalline lift with weight $(k,0)+(h_{\nu},0)$. The new Hodge--Tate weights at $\pi(\sigma_0)$ are $\{k_{\pi(\sigma_0)}-2,0\}$ and at $\pi(\Frob^{-1} \circ \sigma_0)$ they are $
\{p,0\}$. If $\Frob^{-1} \circ \sigma_0 \in M$, we repeat this step and write:
\[
    \xi|_{I_{K'}}=\omega_{\Frob^{-1} \circ \sigma_0}^{-1} \omega_{\Frob^{-2} \circ \sigma_0}^p \cdot \left( \omega_{\sigma_0}^{-1} \omega_{\Frob^{-1} \circ \sigma_0}^p \prod_{\sigma \in J} \omega_{\sigma}^{k_{\pi(\sigma)}-1} \right).
\]
We repeat this until $\Frob^{-s} \circ \sigma_0 \not \in M$. We do this for each $\nu \in \tilde{M}$ to obtain $(k',0)$.

To obtain $(k^{\mu},l^{\mu})$, we do the above for all $\nu \in M \setminus \{\mu \}$. Then choose $\sigma_0 \in \Sigma'$ such that $\sigma_0 \not \in J$ and $\pi(\sigma_0)=\mu$. Then proceed as above for this choice of $\sigma_0$.

Finally to obtain $(k^{\Theta},l^{\Theta})$, we repeat the step we did to obtain $(k^{\mu},l^{\mu})$ for $\mu \in \tilde{M}$, but now for all $\nu \in \tilde{M}$.
\end{proof}

\subsection{Obtaining irregular lifts}
Fix an irregular weight $(k,0)$ such that $1 \leq k_{\tau} \leq p$, $k_{\tau} \neq 2$ if $k_{\Frob^{-1} \circ \tau}=1$, and $k \neq 1$. 

Consider $(k',0)$ as in Section \ref{regularweightdescription} and as before, write $\{b'_{\tau,1},b'_{\tau,2}\}_{\tau}$ for the corresponding Hodge--Tate weights. Similarly, for each $\mu \in \tilde{M}$, consider $(k^{\mu},l^{\mu})$ and write $\{b^{\mu}_{\tau,1},b^{\mu}_{\tau,2}\}_{\tau}$ for the corresponding Hodge--Tate weights and write $\{b^{\Theta}_{\tau,1}, b^{\Theta}_{\tau,2}\}_{\tau}$ for the weight $(k^{\Theta},l^{\Theta})$.

Fix subsets $J', \{J^{\mu}\}_{\mu \in \tilde{M}}$ and $J^{\Theta}$ of $\Sigma'$. We set
\begin{align}
\label{eq: irrrexpressions}    
\wtop'_{\sigma}=\begin{cases}
    b'_{\pi(\sigma),1} & \sigma \in J', \\
    b'_{\pi(\sigma),2} & \sigma \not \in J',     
\end{cases}
\quad \quad
\wtop^{\mu}_{\sigma}=\begin{cases}
    b^{\mu}_{\pi(\sigma),1} & \sigma \in J^{\mu}, \\
    b^{\mu}_{\pi(\sigma),2} & \sigma \not \in J^{\mu},    
\end{cases}
\quad \quad
\wtop^{\Theta}_{\sigma}=\begin{cases}
    b^{\Theta}_{\pi(\sigma),1} & \sigma \in J^{\Theta}, \\
    b^{\Theta}_{\pi(\sigma),2} & \sigma \not \in J^{\Theta},    
\end{cases}
\end{align}
and similarly
\begin{align}
\label{eq: irrsexpressions}    
\wbottom'_{\sigma}=\begin{cases}
    b'_{\pi(\sigma),2} & \sigma \in J', \\
    b'_{\pi(\sigma),1} & \sigma \not \in J',     
\end{cases}
\quad \quad
\wbottom^{\mu}_{\sigma}=\begin{cases}
    b^{\mu}_{\pi(\sigma),2} & \sigma \in J^{\mu}, \\
    b^{\mu}_{\pi(\sigma),1} & \sigma \not \in J^{\mu},    
\end{cases}
\quad \quad
\wbottom^{\Theta}_{\sigma}=\begin{cases}
    b^{\Theta}_{\pi(\sigma),2} & \sigma \in J^{\Theta}, \\
    b^{\Theta}_{\pi(\sigma),1} & \sigma \not \in J^{\Theta}.    
\end{cases}
\end{align}

Recall we have labelled the embeddings $\tau$ in $\Sigma$ letting $\tau_{i+1}=\Frob^{-1} \circ \tau_{i}$ for some $\tau_0 \in \Sigma$. Next we label the embeddings in $\Sigma' = \{ \sigma_0, \dots, \sigma_{2f-1}\}$ where $\pi(\sigma_i)=\pi(\sigma_{f+i})=\tau_i$. We will write $s'_i$ and $t'_i$ to mean $s'_{\sigma_i}$ and $t'_{\sigma_i}$ and similarly for the other weights.

We are going to continue to work in \emph{blocks}. We say a string $\sigma_i,\dots, \sigma_{i+\singleparam}$ is a block in $\Sigma'$ if $\pi(\sigma_i),\dots, \pi(\sigma_{i+\singleparam})$ is a block in the sense of Definition \ref{defblock}.

Assume an irreducible Galois representation $\rho: G_K \to \Fpb$ has crystalline lifts of weights $(k',0)$ and $\{(k^{\mu},l^{\mu})\}_{\mu \in \tilde{M}}$.
By Theorem \ref{thm101} we have balanced subsets $J'$ and $J^{\mu}$ for each $\mu \in \tilde{M}$ such that 



\begin{align}
\label{eq: irrcongstring}
    \sum_{i \in \{0, \dots 2f-1\}} \wtop'_{i} \cdot p^{2f-i-1} \equiv     \sum_{i \in \{0, \dots 2f-1\}} \wtop^{\mu}_{i} \cdot p^{2f-i-1} \mod (p^{2f}-1),
\end{align}
and 
\begin{align}
\label{eq: sirrcongstring}
    \sum_{i \in \{0, \dots 2f-1\}} \wbottom'_{i} \cdot p^{2f-i-1} \equiv     \sum_{i \in \{0, \dots 2f-1\}} \wbottom^{\mu}_{i} \cdot p^{2f-i-1} \mod (p^{2f}-1),
\end{align}
for each $\mu \in \tilde{M}$.

Similarly, if $\rho$ has lifts of weights $(k',0)$ and $(k^{\Theta},l^{\Theta})$, then
\begin{align}
\label{eq: thetacongstring}
    \sum_{i \in \{0, \dots 2f-1\}} \wtop'_{i} \cdot p^{2f-i-1} \equiv     \sum_{i \in \{0, \dots 2f-1\}} \wtop^{\Theta}_{i} \cdot p^{2f-i-1} \mod (p^{2f}-1),
\end{align}
and 
\begin{align}
\label{eq: thetabottomcongstring}
    \sum_{i \in \{0, \dots 2f-1\}} \wbottom'_{i} \cdot p^{2f-i-1} \equiv     \sum_{i \in \{0, \dots 2f-1\}} \wbottom^{\Theta}_{i} \cdot p^{2f-i-1} \mod (p^{2f}-1).
\end{align}

We first need some auxiliary results.

\begin{proposition}
\label{irrj0emb}
Suppose an irreducible representation $\rho: G_K \to \GL_2(\Fpb)$ has crystalline lifts of weights $(k',0)$ and $(k^{\mu},l^{\mu})$ for some $\mu \in \tilde{M}$. Let $\sigma_i \in \Sigma'$ be such that $\pi(\sigma_i)=\mu \in \tilde{M}$. 
Then either
\[
\sigma_i \in J' \text{ and } \sigma_{i+1}, \dots, \sigma_{i+\singleparam} \in J',
\]
or
\[
\sigma_i \not \in J' \text{ and } \sigma_{i+1}, \dots, \sigma_{i+\singleparam} \not \in J',
\]
where $\pi(\sigma_{i+1}), \dots, \pi(\sigma_{i+\singleparam}) \in \text{Block}_{\mu} \cap J_0$ and $\pi(\sigma_{i+\singleparam})$ is the last embedding in the block.
\end{proposition}
\begin{proof}
We obtain balanced subsets $J'$ and $J^{\mu}$ in $\Sigma'$ from Theorem \ref{thm101}. For $\sigma_i=\mu$, we suppose $\sigma_i \in J'$ and $\sigma_i \in J^{\mu}$. Then we have
\[
\wtop'_{\sigma_i}-\wtop^{\mu}_{\sigma_i}=(k_{\mu}-2) - (k_{\mu}-1) = -1.
\]
Now by \eqref{eq: irrcongstring} and Lemma \ref{lem71} we must have that
\begin{align}
\label{eq: irrequation}
(\wtop'_{\sigma_{i+1}}-\wtop_{\sigma_{i+1}}^{\mu}, \dots, 
\wtop'_{\sigma_{i+n-1}}-\wtop_{\sigma_{i+n-1}}^{\mu}, \wtop'_{\sigma_{i+n}}-\wtop_{ \sigma_{i+n}}^{\mu}) \\ = (p-1,\dots,p-1,p), \nonumber
\end{align}
which forces $\sigma_{i+1}, \dots, \sigma_{i+n} \in J'$ by \eqref{eq: bprime} and \eqref{eq: bmu}.


Now suppose $\sigma_i \in J'$ and $\sigma_i \not \in J^{\mu}$. Then we have
\[
\wtop'_{\sigma_i}-\wtop^{\mu}_{\sigma_i}=(k_{\mu}-2)-(-1) = k_{\mu}-1 
\]
which forces $k_{\mu}=p$ by \eqref{eq: congstring} and Lemma \ref{lem71} (since $k_{\mu} \neq 2$ by assumption). Again by Lemma \ref{lem71}, \eqref{eq: irrequation} must hold, 
which forces $\sigma_{i+1}, \dots, \sigma_{i+n} \in J'$.

The cases where $\sigma_i \not \in J' \cap J^{\mu}$ and $\sigma_i \not \in J'$ but $\sigma_i \in J^{\mu}$ follow by Remark \ref{symmrandh}.
\end{proof}

\begin{proposition}
\label{irrj0embtheta}
Suppose an irreducible representation $\rho: G_K \to \GL_2(\Fpb)$ has crystalline lifts of weights $(k',0)$ and $(k^{\Theta},l^{\Theta})$. Let $\mu \in \tilde{M}$. Let $\sigma_i \in \Sigma'$ be such that $\pi(\sigma_i)=\mu \in \tilde{M}$. 
Then either
\[
\sigma_i \in J' \text{ and } \sigma_{i+1}, \dots, \sigma_{i+\singleparam} \in J',
\]
or
\[
\sigma_i \not \in J' \text{ and } \sigma_{i+1}, \dots, \sigma_{i+\singleparam} \not \in J',
\]
where $\pi(\sigma_{i+1}), \dots, \pi(\sigma_{i+\singleparam}) \in \text{Block}_{\mu} \cap J_0$ and $\pi(\sigma_{i+\singleparam})$ is the last embedding in the block.
\end{proposition}
\begin{proof}
This follows immediately from the proof of Proposition \ref{irrj0emb} given that the weights $(k^{\mu},l^{\mu})$ and $(k^{\Theta},l^{\Theta})$ are identical for all embeddings in $\text{Block}_{\mu}$.
\end{proof}

We can now prove that if $\rho$ has a lift of the set of regular weights, that then $\rho$ also has a lift of irregular weight $(k,0)$.

\begin{proposition}
\label{irreducibletoirregular}
Let $K$ be an unramified finite extension of $\Qp$. Let $(k,0) \in \Z^{\Sigma} \times \Z^{\Sigma}$ be an irregular weight satisfying $1 \leq k_{\tau} \leq p$ for all $\tau \in \Sigma$ and such that $k_{\tau} \neq 2$ if $k_{\Frob^{-1} \circ \tau}=1$.  Suppose further $k \neq 1$.
Suppose an irreducible representation $\rho: G_K \to \GL_2(\Fpb)$ either has crystalline lifts of weights $(k',0)$ and $(k^{\mu},l^{\mu})$ for all $\mu \in \tilde{M}$, or has crystalline lifts of weights $(k',0)$ and $(k^{\Theta},l^{\Theta})$, then $\rho$ has a crystalline lift of weight $(k,0)$.
\end{proposition}
\begin{proof}
We show that there is a set $J \subset \Sigma'$ so that we have
\[
    \rho|_{I_{K}} \sim \mtwo{\prod_{\sigma \in J} \omega_{\sigma}^{b_{\pi(\sigma),1}} \prod_{\sigma \not \in J} \omega_{\sigma}^{b_{\pi(\sigma),2}}}{0}{0}{\prod_{\sigma \in J} \omega_{\sigma}^{b_{\pi(\sigma),1}} \prod_{\sigma \not \in J} \omega_{\sigma}^{b_{\pi(\sigma),2}}}, 
\]
where we write $\{b_{\pi(\sigma),1},b_{\pi(\sigma),2}\}$ for the Hodge--Tate weights corresponding to $(k,0)$.

We first assume our Galois representation $\rho$ has lifts $(k',0)$ and $(k^{\mu},l^{\mu})$ for all $\mu \in \tilde{M}$. So we find balanced subsets $J'$ and $\{J^{\mu}\}_{\mu \in \tilde{M}}$ in $\Sigma'$ by Theorem \ref{thm101}. We will choose our set $J$ so that the above holds. Note that for the weight $(k,0)$ we have $b_{\tau,1},b_{\tau,2}=0$ for all $\tau \in J_0$. So we only need to decide for $\Sigma' \setminus \pi^{-1}(J_0)$.

Recall that for all $\sigma$ such that $\pi(\sigma) \in \Sigma \setminus(J_0 \cup \tilde{M})$ we have
\[
{b_{\pi(\sigma),1}}={b'_{\pi(\sigma),1}}={b^{\mu}_{\pi(\sigma),1}},
\]
and
\[
{b_{\pi(\sigma),2}}={b'_{\pi(\sigma),2}}={b^{\mu}_{\pi(\sigma),2}}.
\]
For such embeddings we let $\sigma \in J \iff \sigma \in J'$.

We proceed per blocks. Say $\sigma_0, \dots, \sigma_{E}$ is a block, so that $\pi(\sigma_0), \dots \pi(\sigma_E)$ is a block as in Definition \ref{defblock}.  
Let $\nu$ be the unique embedding in this block such that $\nu \in \tilde{M}$. Suppose we have $\pi(\sigma_{i_{\nu}})=\nu$ for $0 \leq i_{\nu} \leq E$.

Now suppose $\sigma_{i_{\nu}} \in J'$, then we let $\sigma_{i_{\nu}} \in J$. 
We write
\[
\wtop_{i}= \begin{cases}
    b_{\pi(\sigma_i),1} & \sigma_i \in J, \\
    b_{\pi(\sigma_i),2} & \sigma_i \not \in J, \end{cases}
    \quad \text{ and } \quad 
\wbottom_{i}= \begin{cases}
b_{\pi(\sigma_i),1} & \sigma_i \not \in J, \\
b_{\pi(\sigma_i),2} & \sigma_i \in J. \\
\end{cases}
\]
We obtain the following congruences
\[
\sum_{i \in \{0, \dots, E \}} \wtop_i p^{2f-i-1} \equiv \sum_{i \in \{0, \dots, E \}} \wtop'_i p^{2f-i-1}
\]
since $\wtop_i=\wtop'_i$ for $i=0,  \dots, i_{\nu}-1$ and $\wtop_{i_{\nu}}=k_{\nu}-1$ and $\wtop'_{i_{\nu}}=k_{\nu}-2$ with $\wtop'_{i_{\nu}+1}, \dots, \wtop'_E = p-1, \dots, p-1,p$ by Proposition \ref{irrj0emb} (with $\wtop_{i_{\nu}+1}, \dots, \wtop_E = 0$). 
If $\sigma_{i_{\nu}} \not \in J'$, then we set $\sigma_{i_{\nu}} \not \in J$, in which case $\wtop_i=\wtop'_i=0$ for all $i \geq i_{\nu}$ within this block, and $\wtop_i=\wtop'_i$ for $i=0,  \dots, i_{\nu}-1$ so that the congruence also holds in this case.

This argument works exactly the same for each and any other block. Now putting the congruences together we indeed find:
\[
\sum_{i \in \{0, \dots, 2f-1 \}} \wtop_i p^{2f-i-1} \equiv \sum_{i \in \{0, \dots, 2f-1 \}} \wtop'_i p^{2f-i-1},
\]
and similarly
\[
\sum_{i \in \{0, \dots, 2f-1 \}} \wbottom_i p^{2f-i-1} \equiv \sum_{i \in \{0, \dots, 2f-1 \}} \wbottom'_i p^{2f-i-1}.
\]
This gives 
\[
    \rho|_{I_{K}} \sim \mtwo{\prod_{\sigma \in J} \omega_{\sigma}^{b_{\pi(\sigma),1}} \prod_{\sigma \not \in J} \omega_{\sigma}^{b_{\pi(\sigma),2}}}{0}{0}{\prod_{\sigma \in J} \omega_{\sigma}^{b_{\pi(\sigma),1}} \prod_{\sigma \not \in J} \omega_{\sigma}^{b_{\pi(\sigma),2}}}. 
\]
and the result now follows from lifting $\xi = \prod_{\sigma \in J} \omega_{\sigma}^{b_{\pi(\sigma),1}} \prod_{\sigma \not \in J} \omega_{\sigma}^{b_{\pi(\sigma),2}}$ to a character $\tilde{\xi}$ such that $\Ind_{G_{K'}}^{G_K} \tilde{\xi}$ is the lift of $\rho$ with Hodge--Tate weights $(k,0)$ (which we can do by Lemma \ref{cryslift}).
The case where $\rho$ has lifts of weights $(k',l')$ and $(k^{\Theta},l^{\Theta})$ is identical upon replacing Proposition \ref{irrj0emb} by Proposition \ref{irrj0embtheta} in above argument.
\end{proof}

\appendix

\section{Appendix}

This appendix consists of three parts. The first part has to do with the expression $\alpha_i$ as in Definition \ref{defalpha}. We need to compute these values to determine the existence of maps between rank one Kisin modules. The second part computes these values as well, but the focus of that section is the matching results. In the final part we prove some auxiliary results we need to justify assumptions made in Section \ref{nonsemisimplesection}.

\subsection{Computations of $\alpha_i$}
We have the following auxiliary lemma: 

\begin{lemma}
\label{genalphacalc}
Suppose we have a string of integers $r_0, \dots r_{f-1}$ in the range $[-p,p]$ that, considered as a cyclic list, can be broken up into strings of the form $\pm (-1,p-1,\dots,p-1,p)$ (where there may not be any occurrences of $p-1$) and strings of the form $(0,\dots, 0)$. Then
\[
\alpha'_i(r_0,\dots, r_{f-1}) : = \frac{1}{p^f-1} \sum_{j=1}^f p^{f-j} (r_{j+i}) = \begin{cases}
    \pm 1 & \text{if } r_{i+1} \in \pm \{p-1,p\}, \\
    0 & \text{ otherwise}.
\end{cases} 
\]
\end{lemma}

\begin{proof}
We first note that 
\[
p^{m}\cdot (-1)+p^{(m-1)} \cdot (p-1)+ \dots + p^{m-(s-1)} \cdot (p-1) + p^{m-s} \cdot p = 0,
\]
for any $m > s \geq 1$. This gives $\alpha'_i(r_0,\dots, r_{f-1})=0$ if $r_{i+1}=0, -1$. Next suppose $r_{i+1}=p-1$ or $p$. It is easy to see that in these cases we have
\begin{align*}  
\sum_{j=1}^f p^{f-j} (r_{j+i}) = p^f-1, 
\end{align*}
and so the result follows.
\end{proof}
We want to apply this to the regular and irregular weights we introduced. We define the following notation:
\begin{definition}
 Let $(c_n)_n=c_0,\dots, c_{f-1}$, $(d_n)_n=d_0, \dots, d_{f-1}$ be two sequences of $f$ integers, indexes considered in $\Z/(f-1)\Z$. Then define
 \[
 \alpha_i(c_n),(d_n)) := \frac{1}{p^f-1} \sum_{j=1}^f p^{f-j} ((c_n)_{j+i}-(d_n)_{j+i}).
 \]
\end{definition}
We now go back to our setting. We again use the notation from Section \ref{weightnotation}. Then, for example, for fixed $J$ and $J'$, we have
\begin{align*}   
\alpha_i((s'_n), (s_n))&= \alpha'_i(s'_0-s_0,\dots, s'_{f-1}-s_{f-1}) \\
&= \alpha_i(\ofM'(s'_0,\dots, s'_{f-1};a'),\ofM(s_0,\dots,s_{f-1};a)), 
\end{align*}
where the latter is the value defined in the statement of Lemma \ref{lemma512}.

\begin{lemma}
\label{alpharegtoirreg}
Fix a subset $J \in \Sigma$, and let $J'$, $J^{\Theta}$, $\{J^{\mu}\}_{\mu \in \tilde{M}}$ be subsets as in the proof of Proposition \ref{regularshape}. Then we have
\begin{align*}  
\alpha_i((s'_n), (s_n))= 
\begin{cases}
    1 & i \in \tilde{M} \text{ and } i \in J', \\
    1 & i \in J_0, i \in J' \text{ and } i+1 \in J_0, \\
    0 & \text{otherwise, } 
\end{cases}
\end{align*}
and
\begin{align*}    
\alpha_i((t'_n), (t_n))= 
\begin{cases}
    1 & i \in \tilde{M} \text{ and } i \not \in J', \\
    1 & i \in J_0, i \not \in J' \text{ and } i+1 \in J_0, \\
    0 & \text{otherwise, } 
\end{cases}
\end{align*}
and
\begin{align*}  
\alpha_i((s^{\mu}_n), (s_n))= 
\begin{cases}
    1 & i \in \tilde{M} \setminus \{ \mu \} \text{ and } i \in J^{\mu}, \\
    1 & \tau_i=\mu \text{ and } i \not \in J^{\mu}, \\
    1 & i \in J_0, i \in J^{\mu} \text{ and } i+1 \in J_0, \\
    0 & \text{otherwise, } 
\end{cases}
\end{align*}
and
\begin{align*}  
\alpha_i((t^{\mu}_n), (t_n))= 
\begin{cases}
    1 & i \in \tilde{M}  \setminus \{ \mu \}  \text{ and } i \not \in J^{\mu}, \\
    1 & \tau_i=\mu \text{ and } i \in J^{\mu}, \\
    1 & i \in J_0, i \not \in J^{\mu} \text{ and } i+1 \in J_0, \\
    0 & \text{otherwise, } 
\end{cases}
\end{align*}
and 
\begin{align*}  
\alpha_i((s^{\Theta}_n), (s_n))= 
\begin{cases}
    1 & i \in \tilde{M} \text{ and } i \not \in J^{\Theta}, \\
    1 & i \in J_0, i \in J^{\Theta} \text{ and } i+1 \in J_0, \\
    0 & \text{otherwise, } 
\end{cases}
\end{align*}
and
\begin{align*}   
\alpha_i((t^{\Theta}_n), (t_n))= 
\begin{cases}
    1 & i \in \tilde{M}  \text{ and } i \in J^{\Theta}, \\
    1 & i \in J_0, i \not \in J^{\Theta} \text{ and } i+1 \in J_0, \\
    0 & \text{otherwise. } 
\end{cases}
\end{align*}
\end{lemma}

\begin{proof}
The equalities all follow from Lemma \ref{genalphacalc} and the proof of Proposition \ref{regularshape}, along with Section \ref{weightnotation}.
\end{proof}

\subsection{Comparing regular lifts}
We now prove some matching results. This is fairly straightforward, unless we are in some exceptional situation, as we shall see.

\begin{lemma}
\label{minimalmumatch}
Suppose $\rho: G_K \to \GL_2(\Fpb)$ has crystalline lifts of weights $(k',l')$ and $\{(k^{\mu},l^{\mu})\}_{\mu \in \tilde{M}}$. Let $\mu \in \tilde{M}$. Then there are $J'$ and $J^{\mu}$ as in Corollary \ref{corextclass} such that for all $\tau \in \Sigma$ such that $\tau \not \in \text{Block}_{\mu} \cap J_0$ we have $\tau \in J' \iff \tau \in J^{\mu}$.
\end{lemma}

\begin{proof}
Note that \eqref{eq: congstring} needs to hold for all elements of $\tilde{M}$. For all $\tau \in \Sigma$ such that $\tau \not \in \text{Block}_{\mu}$ we have
\[
b'_{\tau,1}= b^{\mu}_{\tau,1} \text{ and } b'_{\tau,2}= b^{\mu}_{\tau,2}=0.
\]
Now suppose for some $\tau \not \in \text{Block}_{\mu}$ we have $\tau \in J'$ and $\tau \not \in J^{\mu}$, then we obtain
\[
\wtop'_{\tau}-\wtop^{\mu}_{\tau}= b'_{\tau,1} \neq 0.
\]
By Lemma \ref{lem71} we obtain 
\[
(\wtop'_j-\wtop^{\mu}_j, \wtop'_{j+1}-\wtop^{\mu}_{j+1}, \dots, \wtop'_{j+\singleparam-1}-\wtop^{\mu}_{j+\singleparam-1}, \wtop'_{j+\singleparam}-\wtop^{\mu}_{j+\singleparam})= (-1,p-1, \dots, p-1,p),
\]
with $\tau \in \tau_j, \dots, \tau_{j+\singleparam}$. We note that by \eqref{eq: bprime} and \eqref{eq: bmu} all of these embeddings must be in the same block and $\tau_{j+\singleparam}$ must be the last embedding in this block.

To achieve $\wtop'_j-\wtop^{\mu}_j=-1$, we must have $j \not \in J_0$. Then we have $\tau_{j+i} \in \tilde{M}$ for some $j+i$. If $\tilde{M}=\{\mu\}$ we are done, so suppose not. 

In this case we find $\wtop'_{j+i}-\wtop^{\mu}_{j+i}=\pm (k_{\tau_{j+i}}-2)$ or $0$ (since $\tau_{j+i} \not \in \text{Block}_{\mu}$). In order to have  $\wtop'_{j+i}-\wtop^{\mu}_{j+i} \in \{-1,p-1\}$ we find $k_{\tau_{j+i}}=3$ and $\tau_{j+i} \not \in J'$, but $\tau_{j+i}\in J^{\mu}$.

Let $\nu = \tau_{j+i}$ and consider $J^{\nu}$ as given by Corollary \ref{corextclass}. Then we have 
\[
\{b^{\nu}_{j+i,1}, b^{\nu}_{j+i,2}\}=\{k_{\tau_{j+i}}-1, -1\},
\]
so that 
\[
\wtop'_{j+i}-\wtop^{\nu}_{j+i}= \begin{cases} -(k_{\tau_{j+i}}-1) & j+i \in J^{\nu}, \\
1 &j+i  \not \in J^{\nu},
\end{cases} \quad
\wtop^{\mu}_{j+1}-\wtop^{\nu}_{j+i}= \begin{cases}
-1 & j+i \in J^{\nu}, \\
(k_{\tau_{j+i}}-1)  & j+i  \not \in J^{\nu},
\end{cases}
\]
since $k_{\tau_{j+i}}=3$ this gives a contradiction with \eqref{eq: congstring} and Lemma \ref{lem71}.

It remains to prove the claim for an embedding $\tau \in \text{Block}_{\mu}$ with $\tau \not \in J_0$. Suppose $\tau \in J'$ and $\tau \not \in J^{\mu}$. Suppose first $\tau = \mu$, then 
\[
\wtop'_{\mu}-\wtop^{\mu}_{\mu}=(k_{\mu}-2)-(-1)=k_{\mu}-1.
\] 
Since $k_{\mu}\neq 2$, by Lemma \ref{lem71} we must have $k_{\mu}=p$. We again find a string of the form $(-1,p-1,\dots,p-1,p)$ with $p$ corresponding to the last embedding of the block $\tau_E$. To obtain $\wtop'_j-\wtop^{\mu}_j=-1$ we must have  $j \not \in J'$ and $j \in J^{\mu}$ so that $\wtop'_j-\wtop^{\mu}_j=- (k_{\tau_{j}}-1) = - 1$ forcing $k_{\tau_j}=2$. 

This further forces $k_{\tau_{j+i}}=p$ and $\tau_{j+1} \in J'$ and $\tau_{j+1} \not \in J^{\mu}$ for all $j+i \in {j+1, \dots, i_{\mu}}$, as well as $\tau_{i_{\mu}+1}, \dots, \tau_E \in J'$ but $\tau_{i_{\mu}+1}, \dots, \tau_E \not \in J^{\mu}$.

Suppose we are in this exceptional situation. Recall first that the set $J^{\mu}$ obtained from Corollary \ref{corextclass} is not unique. 

We twist our representation $\rho$ by $\omega_{\mu}$. Then the twist of $\rho$ has a lift with Hodge--Tate weights 
\begin{align*} 
&\{b^{\mu}_{\tau,1}, b^{\mu}_{\tau,2}\}, & \text{ if } \tau \neq \mu, \\
& \{b^{\mu}_{\tau,1}+1, b^{\mu}_{\tau,2}+1\}, & \text{ if } \tau = \mu,
\end{align*}
which by \eqref{eq: bmu} means that the weights are of the form $\{r_\tau,0\}_{\tau \in \Sigma}$ for $r_{\tau} \in [0,p]$.
In particular, we have
\[
(r^{\mu}_j,\dots, r^{\mu}_{E})=(1,p-1,\dots,p-1,p),
\]
so that by Definition \ref{jmax} our current set $J^{\mu}$ is not $J_{\text{max}}$. We obtain $J_{\text{max}}$ by removing $\tau_j$ from $J^{\mu}$ and adding $\tau_{j+1},\dots, \tau_E$. 
By Proposition \ref{prop88}, adapting $J^{\mu}$ this way results in an isomorphic Galois representation.
This means the conclusion of the lemma holds for this set (after untwisting).

Finally suppose again $\tau \in \text{Block}_{\mu}$ with $\tau \not \in J_0$ and $\tau \in J'$ and $\tau \not \in J^{\mu}$. Suppose now that $\tau \neq \mu$, so that $\wtop'_{\tau}-\wtop^{\mu}_{\tau}=k_{\tau}-1$. To achieve the string $(-1,p-1,\dots,p-1,p)$ obtained from Lemma \ref{lem71} we must include $\tau_E$ and therefore also $\mu$ with $\mu \in J'$ and $\mu \not \in J^{\mu}$. So the proof for this case follows from the case above where $\tau = \mu$.
\end{proof}

\begin{lemma}
\label{alphaprimemu}
Let $J'$ and $\{J^{\mu}\}_{\mu \in \tilde{M}}$ be as in the previous lemma. Suppose $\mu \in \tilde{M}$ is such that $\mu \in J'$, and assume that for all $\tau_i \in \text{Block}_{\mu} \cap J_0$ we have $\tau_i \in J'$ and $\tau_i \not \in J^{\mu}$. Then we have
\begin{align*}  
\alpha_i((\wtop'_n), (\wtop^{\mu}_n))= 
\begin{cases}
    1 & \tau_i=\mu, \\
    1 & i \in J_0 \cap \text{Block}_{\mu}, \text{ and } i+1 \in J_0, \\
    0 & \text{otherwise, } 
\end{cases}
\end{align*}
and
\begin{align*}
\alpha_i((\wbottom^{\mu}_n), (\wbottom'_n))= 
\begin{cases}
    1 & \tau_i=\mu, \\
    1 & i \in J_0 \cap \text{Block}_{\mu}, \text{ and } i+1 \in J_0, \\
    0 & \text{otherwise. } 
\end{cases}
\end{align*}
\end{lemma}

\begin{proof}
This follows from Lemma \ref{genalphacalc}, Section \ref{weightnotation} and our assumptions on the sets $J'$ and $J^{\mu}$.
\end{proof}

\begin{lemma}
\label{minimalthetamatch}
Suppose $\rho: G_K \to \GL_2(\Fpb)$ has crystalline lifts of weights $(k',l')$ and $(k^{\Theta},l^{\Theta})$. Then there are $J'$ and $J^{\Theta}$ as in Corollary \ref{corextclass} such that for all $\tau_i \in \Sigma$ such that $\tau_i \not \in  J_0$ we have $\tau_i \in J' \iff \tau_i \in J^{\Theta}$.
\end{lemma}
\begin{proof}
This follows from the proof of Lemma \ref{minimalmumatch} applied to each weight $(k^{\mu},l^{\mu})$ as follows. Consider a  block in $(k^{\Theta},l^{\Theta})$, say $\nu \in \tilde{M}$ is the unique embedding of $\tilde{M}$ in this block. Then this block is identical to the corresponding block in the weight $(k^{\nu},l^{\nu})$. This means we can follow the same steps, except in order to make use of $J_{\max}$ in the same way we need to twist and untwist by $\prod_{\mu \in \tilde{M}} \omega_{\mu}$ to eliminate negative weights.
\end{proof}

\begin{lemma}
\label{alphaprimetheta}
Let $J'$ and $J^{\Theta}$ be sets as in Definition \ref{auxiliarydescription}. 
Then
\begin{align*}  
\alpha_i((\wtop^{\gamma}_n), (\wtop^{\Theta}_n))= 
\begin{cases}
    1 & i \in J', i+1 \in J_0,\\
    0 & \text{otherwise, } 
\end{cases}
\end{align*}
and
\begin{align*}
\alpha_i((\wbottom^{\Theta}_n), (\wbottom^{\gamma}_n))= 
\begin{cases}
    1 & i  \in J', i+1 \in J_0,\\
    0 & \text{otherwise. } 
\end{cases}
\end{align*}
\begin{align*}  
\alpha_i((\wtop^{\gamma}_n), (\wtop'_n))= 
\begin{cases}
    1 & i \not \in J', i+1 \in J_0,\\
    0 & \text{otherwise, } 
\end{cases}
\end{align*}
and
\begin{align*}
\alpha_i((\wbottom'_n), (\wbottom^{\gamma}_n))= 
\begin{cases}
    1 & i \not \in J', i+1 \in J_0,\\
    0 & \text{otherwise. } 
\end{cases}
\end{align*}
\end{lemma}

\begin{proof}
This follows immediately from Lemma \ref{genalphacalc} and the weight descriptions in Section \ref{weightnotation}.
\end{proof}

\subsection{Auxiliary results}

In this short section we address the exceptional case of Theorem \ref{thm79} as explained in Remark \ref{remarknotexceptional}. As before, we fix an irregular weight $(k,0)$ such that $1 \leq k_{\tau} \leq p$ for all $\tau \in \Sigma$ and that there is no $\tau \in \Sigma$ such that $k_{\tau} \neq 2$ if $k_{\Frob^{-1} \circ \tau} = 1$. Suppose moreover that $k \neq 1$.

\begin{lemma}
\label{notexceptional}
Let $\rho:G_K \to \GL_2(\Fpb)$ be a reducible Galois representation with a crystalline lift of irregular weight $(k,0)$. Then we are not in the exceptional case of Theorem \ref{thm79}.
\end{lemma}

\begin{proof}
Suppose we are in the exceptional case of Theorem \ref{thm79}. We must have $(r_0,\dots, r_{f-1}) \in \mathcal{P}'$ so that $r_i \in \{0,1,p-1,p\}$. This implies $k_i \in \{1,2,p,p+1\}$. 

Since we have $k_i \leq p$, we obtain $k_i \neq p+1$ for all $i$.  Since we have $k \neq 1$, we then must have $k_i \in \{2,p\}$ for some $i$. Suppose first $k_i=2$, then looking at the conditions on the set $\mathcal{P}'$, we must have $k_{i+1}=p$ and therefore $k_{i+s}=p$ for all $s \geq 0$ which is obviously a contradiction. Similarly, if we have $k_i=p$, we arrive at the same conclusion.
\end{proof}

\begin{lemma}
\label{primenotexceptional}
Let $\rho:G_K \to \GL_2(\Fpb)$ be a reducible Galois representation with a crystalline lift of regular weight $(k',0)$. Suppose $J'$ is such that for all $\mu \in \tilde{M}$ we have
\[
\mu \in J' \text{ and } \tau_i \in J' \text{ for all } \tau_i \in  \Block_{\mu} \cap J_0, 
\]
or
\[
\mu \not \in J' \text{ and } \tau_i \not \in J' \text{ for all } \tau_i \in  \Block_{\mu} \cap J_0. 
\]
Then we are not in the exceptional case of Theorem \ref{thm79}.
\end{lemma}

\begin{proof}
Suppose we are in the exceptional case of Theorem \ref{thm79}. Then we must have $k'_i \in \{1,2,p,p+1\}$. By our assumptions we have $k'_i=p+1$ if and only if $\tau_i \in J_0$ and $\tau_i \not \in M$. 

Let $\mu \in \tilde{M}$, then we have $k'_{i_\mu}=k_{i_\mu}-1$ so we must have $k_{i_{\mu}}=3$ by the above since $3 \leq k_{i_{\mu}} \leq p$.

Using the notation as in Theorem \ref{thm79}, this gives us $r'_{i_\mu}=1$ so we must have $i_\mu \not \in J'$ since we are in the exceptional case of Theorem \ref{thm79}. But if this is the case then we must have $\text{Block}_\mu \cap J_0 \not \in J'$ by our assumptions. But this gives a contradiction with the conditions for the exceptional case as for these embeddings we have $r'_i \in \{p-1,p\}$.
\end{proof}

\begin{lemma}
\label{munotexceptional}
Let $\rho:G_K \to \GL_2(\Fpb)$ be a reducible Galois representation with a crystalline lift of regular weight $(k^{\mu},l^{\mu})$ for some $\mu \in \tilde{M}$. Suppose $J^{\mu}$ is such that
\[
\mu \in J^{\mu} \text{ and } \tau_i \not \in J^{\mu} \text{ for all } \tau_i \in  \Block_{\mu} \cap J_0, 
\]
or
\[
\mu \not \in J^{\mu} \text{ and } \tau_i \in J^{\mu} \text{ for all } \tau_i \in  \Block_{\mu} \cap J_0. 
\]
Then we are not in the exceptional case of Theorem \ref{thm79}.
\end{lemma}

\begin{proof}
We first twist $\rho$ by $\omega_{\mu}$ as in Section \ref{regularweightdescription}. 

Now we can apply Theorem \ref{thm79}. Suppose we are in the exceptional case.  Using the notation as in Theorem \ref{thm79}, we must have $r^{\mu}_i \in \{0,1,p-1,p\}$. 

Suppose ${\mu}=\tau_{i_{\mu}}$. Then we have $r^{\mu}_{i_\mu}=k^{\mu}_{i_{\mu}}$. So we obtain $r^{\mu}_{i_\mu} \in \{p-1,p\}$. So we must have $r^{\mu}_{i_\mu} \in J^{\mu}$ to satisfy the conditions for the exceptional case, but by our assumptions we have $r^{\mu}_{i_\mu+1} \not \in J^{\mu}$. This gives a contradiction with the conditions for the exceptional case as we have $r^{\mu}_{i_\mu+1} \in \{p-1,p\}$.
\end{proof}

\begin{lemma}
\label{thetanotexceptional}
Let $\rho:G_K \to \GL_2(\Fpb)$ be a reducible Galois representation with a crystalline lift of regular weight $(k^{\Theta},l^{\Theta})$. Suppose $J^{\Theta}$ is such that for all $\mu \in \tilde{M}$ we have
\[
\mu \in J^{\Theta} \text{ and } \tau_i \not \in J^{\Theta} \text{ for all } \tau_i \in  \Block_{\mu} \cap J_0, 
\]
or
\[
\mu \not \in J^{\Theta} \text{ and } \tau_i \in J^{\Theta} \text{ for all } \tau_i \in  \Block_{\mu} \cap J_0. 
\]
Then we are not in the exceptional case of Theorem \ref{thm79}.
\end{lemma}

\begin{proof}
We first twist $\rho$ by $\prod_{\mu \in \tilde{M}} \omega_{\mu}$ as in Section \ref{regularweightdescription}. Then we can proceed exactly as in the previous lemma, for any choice of $\mu \in \tilde{M}$.
\end{proof}

\bibliographystyle{plain}

\bibliography{references}

\end{document}